\documentclass[12pt]{article}

\usepackage{amsmath}
\usepackage{amsfonts}
\usepackage{amssymb}
\usepackage[draft]{todonotes}   % notes showed

\usepackage{cases}
\usepackage{xcolor}

\newcommand{\dt}{\partial_t}

\newcommand{\dr}{\partial_r}
\newcommand{\dx}{\partial_x}

\newcommand{\abs}[2]{\bigl| #1 \bigr|^{#2}}
\newcommand{\norm}[2]{\bigl\Arrowvert #1 \bigr\Arrowvert_{#2}}
\newcommand{\supnorm}{L^{\infty}}

\newtheorem{thm}{Theorem}[section]
\newtheorem{lm}{Lemma}

\newenvironment{pf}{\paragraph{Proof}}{\hfill$\square$}
\newenvironment{spf}{\paragraph{Sketch of proof}}{\hfill  $\square$}

\numberwithin{equation}{section}
\allowdisplaybreaks
\title{Global Solutions to Compressible Navier-Stokes Equations with Spherically Symmetric Motion and Free Boundary}

\author{Xin Liu}

\begin{document}

\maketitle

\begin{abstract}
This work is devoted to study the global existence of strong and classical solutions to compressible Navier-Stokes equations with or without density jump on the moving boundary for spherically symmetric motion. We establish a unified method to track the propagation of regularity of strong and classical solutions which works for the cases when density connects to vacuum continuously and with a jump simultaneously. The result we obtain is able to deal with both strong solutions with physical vacuum for which the sound speed is $1/2$-H\"older continuous across the boundary, and classical solutions with physical vacuum when $ 1 < \gamma < 3 $. In contrast to the previous results of global weak solutions, we track the regularity globally-in-time up to the symmetric center and the moving boundary. In particular, the free boundary can be traced.
\end{abstract}

\tableofcontents

\section{Introduction}

\subsection{Description and Background}
The motion of a viscous barotropic gas(or fluid) can be described by the isentropic compressible Navier-Stokes equations. In particular, the following system with constant viscosities governs the spherical motion in three dimensional space,

\begin{equation}\label{SphNS}
	\begin{cases}
		\dt(r^2 \rho) + \dr (r^2 \rho u) = 0 & r \in (0, R(t)), \\
	\dt(r^2 \rho u ) +\dr(r^2 \rho u^2)+ r^2 \dr P = (2\mu + \lambda)r^2\dr \biggl( \dfrac{\dr (r^2 u)}{r^2} \biggr) & r \in (0, R(t)),
	\end{cases}
\end{equation}
where $ \rho, u, R(t) $ represent the density, the radial velocity and the radius of the boundary respectively. The Lam\'e constants $ \mu, \lambda $ denoting the viscosity coefficients would satisfy the relation $ \mu > 0, 3\lambda+2\mu \geq 0 $. Moreover, the pressure potential $ P $ is assumed to depend only on the density. For simplicity in this work, the equation of state is taken as $ P = \rho^\gamma $ with $ \gamma > 1 $. Also, we will work on the  Navier-Stokes system \eqref{SphNS} complemented with the following free boundary conditions,
\begin{equation}\label{FreeBoundary}
	\begin{aligned}
		& [P - (2\mu r^2 \dr u + \lambda\dr (r^2 u))](R(t),t)=0,\\
		& u(0,t) = 0, ~ \dt R(t) = u(R(t),t),
	\end{aligned}
\end{equation} 
where the first boundary condition represents the balance of stress tensor across the gas-vacuum interface. 
Also, the initial data is taken to be
\begin{equation}
%	\begin{aligned}
		R(0) = R_0 , ~ u(r,0) = u_0(r), ~ \rho(r,0) = \rho_0(r), ~ r \in (0,R_0).
%	\end{aligned}
\end{equation}
Without loss of generality, it is assumed $ R_0 = 1 $ in the following. 
Meanwhile, we do not impose any boundary profile on the initial density $ \rho_0 $. In fact, $ \rho_0 $ can connect to the vacuum with or without a jump. 

%with boundary condition
%\begin{equation}
%\begin{cases}
%[P - (2\mu r^2 \dr u + \lambda\dr (r^2 u))](R(t),t)=0\\
%u(0,t) = 0\\
%\dt R(t) = u(R(t),t)
%\end{cases}
%\end{equation}
%Where $ P = \rho^\gamma $. Without lost of generality, we assume $ R(0) = 1 $. 

In particular, $ \rho_0 $ can approach the vacuum continuously across the boundary. Such gas-vacuum interface problem has appeared in plenty of physical scenarios such as astrophysics, shallow water waves etc. For example, the configuration of a non-rotating gaseous star would admit the physical vacuum boundary, i.e.
\begin{equation}\label{PV}
	-\infty < \nabla_n c^2 \leq -C < 0, \rho = 0, ~~~~ \text{on the boundary},
\end{equation}
where $ n $ denotes the normal direction and $ c^2 = P'(\rho) $ is the square of the sound speed. Indeed, the physical vacuum boundary indicates that the sound speed $ c $ is only $1/2$-H\"older continuous instead of Lipschitz continuous across the boundary, which is quite troublesome (see \cite{Liu1996}). Only recently, some local-in-time well-posedness of the smooth solutions for such problems is available for the inviscid flows \cite{Coutand2010,Coutand2012,Coutand2011a,Gu2012,Gu2015,Jang2009b,Jang2015,LuoXinZeng2014} with or without self-gravitation and for the viscous flows \cite{Jang2010} with self-gravitation. As for the global dynamic of flows with physical boundary \eqref{PV}, Luo, Xin, Zeng \cite{LuoXinZeng2016} have shown that with small perturbation of the Lane-Emden solutions, the strong solution to the Navier-Stokes-Poissen system exists globally and converges to the equilibrium state. See \cite{LuoXinZeng2015} for the case with degenerate viscosities.  Meanwhile, Zeng has established the global regularity of the compressible Navier-Stokes equation in one dimensional setting which includes the case of physical vacuum in \cite{ZengHH2015}. Zeng's work extends the one in \cite{Luo2000a}, in which the authors have shown the global existence of the solutions to the one dimensional problem with constant viscosities but higher regularity for the density. 

When the density connects to vacuum with a jump, a global weak solution to the problem with density-dependent viscosities and spherically symmetric motion is obtained in \cite{Guo2012a} by Guo, Li, Xin. Moreover the solution is shown to be smooth away from the centre. Recently, such problem is studied in the setting of spherical symmetry in two dimensional space and $ \mu = constant $, $ \lambda = \lambda(\rho) = \rho^\beta $ with some $ \beta > 1 $ by Li, Zhang \cite{Li2016}. Working in both Lagrangian and Eulerian coordinates, the authors show the global existence of strong solutions. Another noticeable result is from Yeung and Yuen \cite{Yeung2009}, in which the authors have constructed analytic solutions in the case with density-dependent viscosities. Similar results were further studied in \cite{Guo2012b} with or without a density jump across the boundary. Such solutions indicate that the domain of the gas(fluid) would expand as time grows up, and the density would decreases to zero everywhere including the centre. 

Notice, for a spherically symmetric motion, the regularity of solutions at the centre has been only obtained in the case when it is with small moves, induced either by the local-in-time motion \cite{LuoXinZeng2014} or by the small perturbation and stability of the equilibrium state \cite{LuoXinZeng2016,LuoXinZeng2015}. However, when such structures do not exist, it is in general not clear how to perform similar estimates. In this work, we consider the global regularity of \eqref{SphNS} with the free boundary \eqref{FreeBoundary}. For one thing, in the case when the density connects to vacuum continuously, the degeneracy of the system makes the problem challenging as mentioned in \cite{Jang2010}. Indeed, the classical techniques \cite{Zadrzynska2001,Zajaczkowski1999,Zajaczkowski1994,Zajaczkowski1993} does not work in such a situation. Besides the degeneracy on the boundary, we will focus on he regularity at the centre. As mentioned above, the moving domain is expected to expand due to the absence of equilibrium state and large time. In particular, there is no a priori bound on the flow trajectories which causes additional difficulties. 

Before moving onto our strategy to overcome the difficulties listed above, we briefly review some classical works concerning the compressible Navier-Stokes system. 

There are rich literatures studying the Cauchy and first initial boundary value problem. In the absence of vacuum ($\rho \geq \underline\rho > 0 $), the local and global well-posedness of classical solutions have been investigated widely. To name a few, Serrin \cite{Serrin1959} considered the uniqueness of both viscous and inviscid compressible flows. Itaya \cite{Itaya1971} and Tani \cite{Tani1977} considered the Cauchy and first initial boundary problems. See also \cite{Tani1986}. Moreover, the pioneering works of Matsumura and Nishida \cite{Matsumura1980,Matsumura1983} showed the global well-posedness of classical solutions to Navier-Stokes equations with small perturbation of a uniform non-vacuum state. When the vacuum appears, some singular behaviour may occur. For instant, as pointed out by Xin in \cite{zpxin1998}, the classical solution to the Navier-Stokes equations may blow up in finite time. See also in \cite{XinYan2013,Cho2006}. Nevertheless, the local well-posedness theory was developed by Cho and Kim \cite{Cho2006c,Cho2006a} for barotropic and heat-conductive flows. See also \cite{Zhang2007a}. Unfortunately, the solutions from Cho and Kim is not in the same functional space as that in \cite{zpxin1998}. In particular, it can not track the entropy in the vacuum area. With small initial energy, Huang, Li, Xin \cite{HuangLiXin2012} show the global existence of classical solutions but with large oscillations to the isentropic compressible Navier-Stokes equations. However, as pointed out in \cite{Liu1998a}, such Navier-Stokes system might not develop a satisfactory solution at the vacuum state. Therefore, the authors introduced the problem with density-dependent viscosities and showed the well-posedness of the free boundary problem locally. 

As for the free boundary problems, the local well-posedness theory can be tracked back to Solonnikov and Tani \cite{Solonnikov1992}, Zadrzy\'nska and Zaj\c{a}czkowski \cite{Zadrzynska2001,Zadrzynska1994,Zajaczkowski1993,Zajaczkowski1995}. The global well-posedness problem were studied in \cite{Zajaczkowski1999,Zajaczkowski1993,Zajaczkowski1994}. Among these works, the stress tensor on the moving surface is balanced by a force induced by the surface tension or an external pressure. In particular, it admits a uniform state and the global well-posedness was achieved with small perturbation of such constant state. We refer other free boundary problems to \cite{Zhang2009,Jang2011,Okada1989,Okada2004,Ou2015,Zhang2009c} and the references therein. 

We will work in the Lagrangian coordinates induced by the flow trajectories, defined in \eqref{LgCor}. Comparing to the classical choice of Lagrangian mass coordinates, such coordinates would enable us to track the particle path both at the centre and near the moving surface as mentioned in \cite{ZengHH2015}. Indeed, as it will explain itself, the Lagrangian unknown $r$ has ready represented the density $\rho $ and the velocity $ u $ \eqref{LgNSDen}, and $ \rho_0 $ would appear as a degenerate coefficient as in \cite{LuoXinZeng2014}. The coordinate singularity at the center can be understood as follows. Let $ \mathcal U = (V_1,V_2,V_3) = V(r) \cdot \vec{y}/r \in \mathbb R^3 $ be a smooth radial vector field in $ \mathbb R^3 $ where $ r = |\vec{y}| $. Then $ |\nabla \mathcal U|^2 = V'^2 + 2 (\frac{V}{r})^2 $. In the Lagrangian coordinates system, such quantity is $ (\frac{V_x}{r_x})^2 + 2 (\frac{V}{r})^2 $. Thus it can be seen, in order to obtain the desirable bound on the vector field $ v(x,t) = u(r(x,t),t) $, the reasonable quantities to consider are $ \frac{v_x}{r_x} $ and $ \frac{v}{r} $. Such structure has already been studied in \cite{LuoXinZeng2014,LuoXinZeng2015,LuoXinZeng2016}. In the meantime, $ \frac{r}{x} $ a priorly admits only lower bound but not upper bound in contrast to the works by Luo, Xin, Zeng. In other word, we only expect the uniform (in time) bound of $ \frac{x}{r}$ and $ \frac{1}{r_x} $. While in \cite{ZengHH2015}, for the one dimensional problem, $ r_x $ can be estimated through a careful point-wise estimate, similar structure no longer exists for the spherical motion. In fact, the corresponding quantity is now $ \frac{r^2r_x}{x^2}$. However, the structure for this new quantity is not determinant enough (see \eqref{uee003}). In particular, it contains an extra boundary term $ \biggl. \frac{v}{r} \biggr|_{x=1} $. Luckily, %by assuming the initial distribution of density "scatter" enough (\eqref{DensityAsm}) and Remark \ref{rm:scatteringden}), 
the radius $ R(t) $ priorly admits uniform (in time) lower bound. Such structure would enable a $ L_t^2 L_x^\infty $ estimates of $ \frac{v_x}{r_x} $ and $ \frac{v}{r} $. Our $ L_t^\infty L_x^\infty $ estimates are achieved through some point-wise estimates, which are performed by integrating the equation from the boundary and from the centre. In fact, such calculations separate $ \frac{v_x}{r_x} $ and $ \frac{v}{r} $ from the viscosity tensor (see \eqref{ueept001} and \eqref{ueept002}). The benefit of such programs is that it avoids the integral multipliers used in \cite{LuoXinZeng2016} and the Hardy's inequality in \cite{ZengHH2015}, and therefore we are able to manipulate the case when $ \rho_0 $ is with general profiles. Eventually, the bound of $ \frac{v_x}{r_x}, \frac{v}{r}, \frac{x^2}{r^2 r_x}$ can be represented by the a prior bound \eqref{PrAssum} and the initial energy $ \mathcal E_0, \mathcal E_1, \mathcal E_2 $ defined in \eqref{InitialEne}. Then by choosing initial energy small, we shall close the a prior estimate. In addition, we study the relation between the a prior bound and the initial energy, from which it is possible to determine the a prior bound as a map of the initial energy, provided the initial energy small enough. We further discuss the regularity propagated by the solution, which would give the global regularity.  

The rest of the present work would be organised as follows. In the next section, we would introduce the basic notations and the main results in the Lagrangian coordinates. Then we shall start the a priori estimates. Under the a prior assumption, some standard energy estimates would be listed in Section \ref{sec:energyestimates}. The key estimates are studied in Section \ref{sec:uniformestimate}. We establish the uniform point-wise estimates discussed above. Furthermore, the smallness constraints on the initial energy is designed to close the estimates. Section \ref{sec:regularity} and Section \ref{sec:highregularity} are devoted to study the propagation of regularity for strong and classical solutions respectively. In the end, we present the suitable functional frame-work for the local well-posedness problem and sketch the proof of the main theorem.

\subsection{Lagrangian Reformulation and Main Results}

We define the Lagrangian coordinate as follows. For $ x \in (0,1) = (0,R_0) $, let the particle path $ r $ be a function defined by the following ordinary differential equation,
\begin{equation}\label{LgCor}
	\begin{cases}
	\dfrac{d}{dt} r(x,t) = u (r(x,t),t), \\
	r(x,0) = x.
	\end{cases}
\end{equation}
The Lagrangian unknowns are denoted by
\begin{equation*}
f(x,t) : = \rho(r(x,t),t), ~v(x,t) : = u(r(x,t),t).
\end{equation*}
Then the system \eqref{SphNS} can be written as
\begin{equation*}
	\begin{cases}
		(r^2 f)_t + (r^2 f) \cdot \dfrac{v_x}{r_x} = 0,\\
		r^2 f v_t + r^2 \dfrac{P_x}{r_x} = (2\mu + \lambda) \dfrac{r^2}{r_x} \dx \biggl( \dfrac{(r^2 v)_x}{r^2 r_x} \biggr).
	\end{cases}
\end{equation*}
Or equivalently, 
\begin{align}
	& \rho(r(x,t),t) = f(x,t) = \dfrac{x^2\rho_0}{r^2r_x}, ~ u(r(x,t),t) = v(x,t) = \dt r, {\label{LgNSDen}} \\
	& \biggl(\dfrac{x}{r}\biggr)^2 \rho_0 v_t + \biggl\lbrack\biggl(\dfrac{x^2 \rho_0}{r^2 r_x}\biggr)^\gamma \biggr\rbrack_x = (2\mu + \lambda)\biggl\lbrack \dfrac{(r^2 v)_x}{r^2 r_x} \biggr\rbrack_x= \mathfrak{B}_x + 4\mu \biggl( \dfrac{v}{r} \biggr)_x {\label{LgNS}},
\end{align}
where
\begin{align*}
P = \biggl(\dfrac{x^2 \rho_0}{r^2 r_x} \biggr)^\gamma, ~ \mathfrak{B} = (2\mu+\lambda) \dfrac{v_x}{r_x} + 2\lambda\dfrac{v}{r}.
\end{align*}
The equation \eqref{LgNS} is complemented with the boundary condition
\begin{equation} \label{LgNSBC}
	\begin{aligned}
	& \bigl.\biggl(\dfrac{x^2 \rho_0}{r^2 r_x} \biggr)^\gamma - \biggl((2\mu+\lambda) \dfrac{v_x}{r_x} + 2\lambda\dfrac{v}{r} \biggr) \bigr|_{x=1} = 0, \\
	& \bigl. r \bigr|_{x=0} = 0, ~\bigl. v \bigr|_{x=0} = 0,
	\end{aligned}
\end{equation}
and initial data
\begin{equation}\label{LgNSIN}
	r(x,0) = x, ~ r_t(x,0) = u_0(x). 
\end{equation}
Moreover, the following assumption is imposed on the initial density $ \rho_0 $, and the viscosity coefficients $ \mu, \lambda $, 
\begin{align}
&  \rho_0 > 0 ~ \text{for} ~ x \in [0,1), ~ \max_{0\leq x \leq 1} \rho_0 \leq \bar\rho, \label{InitialAssum} \\ % ~ \max_{0\leq x\leq 1} \biggl\lbrace \biggl|\dfrac{u_0}{x} \biggr|, \biggl| u_{0x} \biggr| \biggr\rbrace \leq \beta. 
& \bigl\Arrowvert (\rho_0^\gamma)_{x}\bigr\Arrowvert_{L_x^2(0,1)} < \infty, ~ \mu, \lambda > 0 . \label{InitialAssum2}
\end{align}
Now we define the initial energy we shall use in this work. Denote 
\begin{equation}\label{InitialEne}
\begin{aligned}
& \mathcal{E}_0 = \dfrac{1}{2} \int_0^1 x^2 \rho_0 u_0^2\,dx	+ \dfrac{1}{\gamma-1} \int_0^1 x^2 \rho_0^\gamma \,dx,\\
& \mathcal{E}_1 = \dfrac{1}{2} \int_0^1 x^2 \rho_0 u_{1}^2\,dx, ~~~~ \mathcal{E}_2 = \dfrac{1}{2} \int_0^1 \rho_0  u_{1}^2\,dx, \\
& \mathcal E_3 = \dfrac{1}{2}\int_0^1 x^2 \rho_0 u_{2}^2\,dx, ~~~~ \mathcal E_4 = \dfrac{1}{2}\int_0^1 \rho_0 u_{2}^2\,dx.
\end{aligned}
\end{equation}
where
\begin{align}
	& u_1 = \dfrac{1}{\rho_0}\biggl\lbrace (2\mu + \lambda) \biggl( \dfrac{(x^2 u_{0})_x}{x^2}\biggr)_x - \biggl( \rho_0^\gamma \biggr)_x \biggr\rbrace, \label{intut1}\\
	& u_{2} = \dfrac{1}{\rho_0}\biggl\lbrace (2\mu + \lambda) \biggl\lbrack\dfrac{(x^2u_1)_x}{x^2}  \biggr\rbrack_x + \gamma \biggl\lbrack\rho_0^\gamma \dfrac{(x^2 u_0)_x}{x^2} \biggr\rbrack_x  \biggr.   \nonumber\\
	& ~~~~ \biggl. - (2\mu + \lambda) \biggl\lbrack  u_{0,x}^2 + 2 \dfrac{u_0^2}{x^2} \biggr\rbrack \biggr\rbrace + 2 \dfrac{u_0 u_1}{x}.
\end{align}
Also, let $ M > 0 $ be a constant such that \begin{equation}\label{UpBdInitial}  \max_{0\leq x \leq 1}{\biggl\lbrace \biggl| \dfrac{u_0}{x} \biggr|, \biggl| u_{0,x} \biggr|  \biggr\rbrace} \leq M, \end{equation}
and the following compatible condition on the boundary is imposed,
\begin{equation}\label{InitialBoundary}
\biggl. \rho_0^\gamma - \biggl( (2\mu + \lambda) u_{0,x} + 2 \lambda \dfrac{u_0}{x} \biggr) \biggr|_{x=1} = 0.
\end{equation}

%In addition, it is assumed,
%\begin{equation}\label{DensityAsm}
%	0 < \eta \leq 3^{1/3} \mathcal E_{0}^{-\frac{1}{3(\gamma-1)}} \biggl( \int_0^1 s^2 \rho_0 \,ds \biggr)^{\frac{\gamma}{3(\gamma-1)}} \leq 1.
%\end{equation}
%
%\begin{rmk}\label{rm:scatteringden}
%	The assumption \eqref{DensityAsm} can be understood as follows. 
%\end{rmk}

% ===================================

Then our main theorem is stated as follows
\begin{thm}[Global existence]
Consider the initial boundary problem \eqref{LgNS}, \eqref{LgNSBC}, \eqref{LgNSIN} satisfying the assumption \eqref{InitialAssum}, \eqref{InitialAssum2}, \eqref{UpBdInitial}, \eqref{InitialBoundary}, %\eqref{DensityAsm}. 
There is an $ \bar \epsilon > 0 $ depending on $ \bar\rho_0, \mu,\lambda, M $ such that if
\begin{equation}\label{InitialAssum4}
\mathcal E_0, \mathcal E_1, \mathcal E_2 < \bar\epsilon \end{equation}
it admits a globally defined strong solution $ (r(x,t),v(x,t)) $ to \eqref{LgNS}. Moreover, there are positive constants $ \bar\alpha $ and $ \bar\beta > M $ such that
\begin{equation}
	0 < \dfrac{x^2}{r^2r_x} < \bar\alpha^3, ~~ 0 \leq \biggl| \dfrac{v}{r} \biggr|, \biggl|\dfrac{v_x}{r_x}\biggr| < \bar\beta.
\end{equation}
And the following regularity holds for any $ 0<T<\infty $, 
\begin{equation}\label{regularity1}
	\begin{cases}
		x\sqrt{\rho_0} v, x \sqrt{\rho_0} v_t, \sqrt{\rho_0}v_t \in L_t^\infty((0,T),L_x^2(0,1)),\\
		v, v_x, v_{xx}, r, r_x, r_{xx}, \dfrac{v}{x},  \biggl(\dfrac{v}{x}\biggr)_x, \biggl(\dfrac{r}{x}\biggr)_x \in L_t^\infty((0,T),L_x^2(0,1)), \\
		x v_x , x v_{xt},  v , v_t, \dfrac{v_t}{x}, v_{xt} \in L_t^2((0,T),L_x^2(0,1)).
	\end{cases}
\end{equation}
If, in addition,
\begin{equation}\label{InitialAssum3}
	\mathcal E_3, \mathcal E_4 < \infty, ~ \biggl\Arrowvert (\rho_0)_x \biggr\Arrowvert_{L_x^2(0,1)}, \biggl\Arrowvert (\rho_0^\gamma)_{xx} \biggr\Arrowvert_{L_x^2(0,1)} < \infty,
\end{equation}
the global strong solution is smooth, and satisfies, along with \eqref{regularity1}
\begin{equation}\label{regularity2}
	\begin{cases}
		x \sqrt{\rho_0} v_{tt}, \sqrt{\rho_0} v_{tt}, v_{xxt}, \biggl(\dfrac{v_t}{x}\biggr)_x \in L_t^\infty((0,T),L_x^2(0,1)), \\
		r_{xxx}, \biggl(\dfrac{r}{x}\biggr)_{xx}, v_{xxx}, \biggl(\dfrac{v}{x}\biggr)_{xx} \in L_t^\infty((0,T),L_x^2(0,1)), \\ 
		x v_{xtt}, v_{tt}, v_{xtt}, \dfrac{v_{tt}}{x} \in L_t^2((0,T),L_x^2(0,1)),
	\end{cases}
\end{equation}
for $ 0 < T < \infty $.

\end{thm}

\section{A Prior Estimates}\label{sec:estimates}

Through this section, it is assumed $ (r,v) = (r(x,t),v(x,t)) $ is a smooth solution to \eqref{LgNS} with \eqref{LgNSBC}, \eqref{LgNSIN} such that all the integrations by parts in the following hold. Moreover, it is a priorly assumed
\begin{equation}\label{PrAssum}
	0 < \dfrac{x^2}{r^2r_x} \leq \alpha^3 , ~ 0 \leq \biggl|\dfrac{v}{r}\biggr|, \biggl|\dfrac{v_x}{r_x}\biggr| \leq \beta
\end{equation}
for some $ \alpha > 1 $ and $ \beta \geq M $ ($ M $ is defined in \eqref{UpBdInitial}). 
%The goal is to show
%\begin{equation}
%	0 < \dfrac{x^2}{r^2r_x} \leq 2, 0 \leq \biggl|\dfrac{v}{r}\biggr|, \biggl|\dfrac{v_x}{r_x}\biggr| \leq 2\beta.
%\end{equation}
From \eqref{PrAssum} it can be derived
\begin{equation}\label{PrAssum2}
	\begin{gathered}
	r^3 = 3 \int_0^x r^2 r_x\,dx \geq 3 \alpha^{-3} \int_0^x x^2 \,dx = \alpha^{-3} x^3, ~ \text{or} ~ 0 < \dfrac{x}{r} \leq \alpha; \\ \text{in particular} ~ R(t) \geq \alpha^{-1} .
	\end{gathered}
\end{equation}
We will justify the a prior bound \eqref{PrAssum} in the end of Section \ref{sec:uniformestimate}. In the following, unless specified, it is denoted by
\begin{displaymath}
	\int \cdot \,dx = \int_0^1 \cdot \,dx, ~~ \int \cdot \,dt = \int_0^T \cdot \,dt ~ \text{for} ~ T>0.
\end{displaymath}

%===================================

\subsection{Basic Energy Estimates}\label{sec:energyestimates}
In this section, we start by some basic estimates under the a priori assumption \eqref{PrAssum}. Moreover, the estimates in this section are independent of time. 
For a constant $ C_0>0 $, determined in \eqref{lm2}, denote
\begin{equation} \label{TEnergy0}
%\begin{aligned}
	 E_0 = \mathcal{E}_0, ~ E_1 = \mathcal{E}_1 + C_0(\alpha^{6\gamma} + \beta^2) \mathcal{E}_{0}, ~ E_2 = \mathcal{E}_2.
%\end{aligned}
\end{equation}
The first lemma is concerning the kinetic and potential energy of \eqref{LgNS}. 
\begin{lm}[Basic Energy Estimate]\label{lm:basicest}
For a smooth solution to \eqref{LgNS}, the kinetic and potential energies satisfy the following identity. 
	\begin{equation}\label{lm1}
	\begin{aligned}
		& \dfrac{1}{2} \int x^2\rho_0 v^2\,dx + \dfrac{1}{\gamma-1} \int r^2 r_x \biggl( \dfrac{x^2 \rho_0}{r^2r_x}\biggr)^\gamma\,dx \\
		& ~~  + 2\mu \int \int \biggl( r^2 \dfrac{v_x^2}{r_x} +  2 r_x v^2\biggr) \,dx  \,dt  + \lambda \int \int  r^2 r_x \biggl( \dfrac{v_x}{r_x} + 2 \dfrac{v}{r} \biggr)^2\,dx  \,dt = E_0.
	\end{aligned}
	\end{equation}
\end{lm}

\begin{pf}
	Multiply \eqref{LgNS} with $ r^2 v $ and integrate the resulting in spatial variable, 
	\begin{equation}\label{ee001}
		\int x^2 \rho_0 v_t v \,dx + \int \biggl\lbrack\biggl(\dfrac{x^2 \rho_0}{r^2r_x}\biggr)^\gamma\biggr\rbrack_x r^2 v\,dx = \int \mathfrak{B}_x r^2 v\,dx + \int 4\mu \biggl(\dfrac{v}{r} \biggr)_x r^2 v\,dx.
	\end{equation}
	Integration by parts yields
	\begin{displaymath}
		\begin{aligned}
			& \int \biggl\lbrack\biggl(\dfrac{x^2 \rho_0}{r^2r_x}\biggr)^\gamma\biggr\rbrack_x r^2 v\,dx - \int \mathfrak{B}_x r^2 v\,dx = - \int \biggl(\dfrac{x^2 \rho_0}{r^2r_x} \biggr)^\gamma (r^2 v)_x \,dx\\
			& ~~~ + \int \mathfrak{B} (r^2 v)_x\,dx = \dfrac{d}{dt} \biggl\lbrace \dfrac{1}{\gamma-1} \int r^2 r_x \biggl(\dfrac{x^2 \rho_0}{r^2 r_x} \biggr)^\gamma\,dx \biggr\rbrace + \int \mathfrak{B} (r^2 v)_x\,dx.
		\end{aligned}
	\end{displaymath}
	Moreover, 
	\begin{displaymath}
		\begin{aligned}
			& \int \mathfrak{B} (r^2 v)_x\,dx - \int 4\mu \biggl(\dfrac{v}{r}\biggr)_x r^2 v\,dx \\
			& ~~~~~~~ = \int (2\mu + \lambda) r^2 \dfrac{v_x^2}{r_x} + 4\lambda r v v_x + (4\mu + 4\lambda) r_x v^2\,dx\\
			& ~~~~~~~ = 2\mu \int \biggl( r^2 \dfrac{v_x^2}{r_x} + 2 r_x v^2\biggr) \,dx + \lambda \int r^2 r_x \biggl( \dfrac{v_x}{r_x} + 2 \dfrac{v}{r} \biggr)^2\,dx.
		\end{aligned}
	\end{displaymath}
	Therefore, \eqref{ee001} can be written as
	\begin{equation}\label{ee002}
	\begin{aligned}
		& \dfrac{d}{dt}\biggl\lbrace \dfrac 1 2 \int x^2 \rho_0 v^2 \,dx + \dfrac{1}{\gamma-1} \int r^2 r_x \biggl( \dfrac{x^2 \rho_0}{r^2r_x}\biggr)^\gamma\,dx \biggr\rbrace \\
		& ~~~~~~~~ + 2\mu \int \biggl( r^2 \dfrac{v_x^2}{r_x} + 2 r_x v^2\biggr) \,dx + \lambda \int r^2 r_x \biggl( \dfrac{v_x}{r_x} + 2 \dfrac{v}{r} \biggr)^2\,dx = 0.
	\end{aligned}
	\end{equation}
	Integrating over temporal variable yields \eqref{lm1}.
\end{pf}

%=========================

The next lemma concerns the time derivative of \eqref{LgNS}.

\begin{lm} There is a constant $ C_0 > 0 $ depending on $ \mu, \lambda, \bar\rho_0 $ such that
	\begin{equation}\label{lm2}
	\begin{aligned}
		& \dfrac{1}{2} \int x^2 \rho_0 v_t^2 \,dx +  \mu \int \int \biggl( r^2 \dfrac{v_{xt}^2}{r_x} + 2 r_x v_t^2 \biggr) \,dx  \\& ~~~~~~~~ %+ \lambda \int \int r^2 r_x \biggl( \dfrac{v_{xt}}{r_x} + 2 \dfrac{v_t}{r} \biggr)^2\,dx \,dt 
		\leq \mathcal E_1 + C_0 \biggl(\alpha^{6\gamma}+\beta^2\biggr) \mathcal E_0 = E_1.
	\end{aligned}
	\end{equation}
\end{lm}

\begin{pf}
	Multiply \eqref{LgNS} with $ r^2 $ and take the time derivative of the resulting equation. It follows,
	\begin{equation}\label{LgNS1}
		\begin{aligned}
			& x^2 \rho_0 v_{tt} +  r^2 \biggl\lbrack \biggl( \dfrac{x^2 \rho_0}{r^2 r_x} \biggr)^\gamma \biggr\rbrack_{xt} = r^2 \mathfrak{B}_{xt} + 4\mu r^2 \biggl(\dfrac{v}{r}\biggr)_{xt} \\
			& ~~~~~~~ + 2 r v \biggl( \mathfrak{B}_x + 4\mu \biggl(\dfrac{v}{r}\biggr)_x \biggr) - 2 r v \biggl\lbrack \biggl(\dfrac{x^2\rho_0}{r^2r_x}\biggr)^\gamma \biggr\rbrack_x.
		\end{aligned}
	\end{equation}
	\eqref{LgNS1} would be complemented with additional boundary conditions
	\begin{equation}
		\biggl. \biggl\lbrace \biggl\lbrack \biggl( \dfrac{x^2 \rho_0}{r^2 r_x} \biggr)^\gamma \biggr\rbrack_{t} - \mathfrak{B}_{t} \biggr\rbrace \biggr|_{x=1} = 0, ~ v_t(0,t) = 0.
	\end{equation}
	Multiply \eqref{LgNS1} with $ v_t $ and integrate the resulting in spatial variable. After integration by parts, it follows,
	\begin{equation}\label{ee003}
		\begin{aligned}
			& \int x^2 \rho_0 v_{tt}v_t \,dx  - \int \biggl\lbrack \biggl(\dfrac{x^2\rho_0}{r^2r_x}\biggr)^\gamma \biggr\rbrack_t (r^2 v_t)_x\,dx + \int \mathfrak{B}_t (r^2 v_t)_x\,dx \\ & ~~~~ - \int 4 \mu r^2 v_t \biggl(\dfrac{v}{r} \biggr)_{xt}\,dx = 2 \int \biggl(\dfrac{x^2 \rho_0}{r^2r_x}\biggr)^\gamma (r v v_t)_x \,dx - 2 \int \mathfrak{B} (r v v_t )_x\,dx \\ & ~~~~ + 8\mu \int r v v_t \biggl(\dfrac{v}{r}\biggr)_x\,dx := L_1 + L_2 + L_3.
		\end{aligned}
	\end{equation}
	Notice, by using \eqref{PrAssum}
	\begin{displaymath}
		\begin{aligned}
			& \biggl| \int \biggl\lbrack \biggl(\dfrac{x^2\rho_0}{r^2r_x}\biggr)^\gamma \biggr\rbrack_t (r^2 v_t)_x\,dx \biggr| = \biggl| \gamma \int \biggl(\dfrac{ x^2 \rho_0 }{r^2 r_x} \biggr)^\gamma \dfrac{(r^2r_x)_t}{r^2 r_x} (r^2 v_t)_x\,dx \biggr| \\
			& ~~~~~~ \leq C \alpha^{3\gamma} \int \biggl|\biggl( \dfrac{v_x}{r_x} + 2 \dfrac{v}{r} \biggr) (r^2 v_t )_x \biggr| \,dx \\
			& ~~~~~~ \leq \ C \alpha^{3\gamma} \int \biggl|4r_x v v_t\biggr| + \biggl|2 r v v_{xt}\biggr| + \biggl|2 r v_t v_x\biggr| + \biggl|\dfrac{r^2}{r_x} v_x v_{xt}\biggr| \,dx\\
			& ~~~~~~ \leq \delta \biggl\lbrace \int r^2 \dfrac{v_{xt}^2}{r_x}\,dx + \int r_x v_t^2\,dx \biggr\rbrace + C_{\delta} \alpha^{6\gamma} \biggl\lbrace \int r_x v^2\,dx + \int r^2 \dfrac{v_x^2}{r_x} \,dx   \biggr\rbrace. 
	\end{aligned}\end{displaymath} 
	In the meantime, 
	\begin{displaymath}
		\begin{aligned}
			& \int \mathfrak{B}_t (r^2 v_t)_x\,dx - \int 4 \mu r^2 v_t \biggl(\dfrac{v}{r} \biggr)_{xt}\,dx \\
			& ~~~~~~ = 2 \mu \int \biggl( r^2 \dfrac{v_{xt}^2}{r_x} + 2 r_x v_t^2 \biggr) \,dx + \lambda\int r^2 r_x \biggl( \dfrac{v_{xt}}{r_x} + 2 \dfrac{v_t}{r} \biggr)^2\,dx - L_4,
		\end{aligned}
	\end{displaymath}
	where
	\begin{displaymath}
	\begin{aligned}
		& L_4 = (4\mu + 2 \lambda) \int \biggl(\dfrac{v_x}{r_x}\biggr) r v_t v_x \,dx + (2\mu + \lambda) \int \biggl(\dfrac{v_x}{r_x}\biggr) \dfrac{r^2 v_x v_{xt}}{r_x}\,dx \\
		& ~~~~ + (8\mu + 4\lambda) \int \biggl(\dfrac{v}{r} \biggr) r_x v v_t \,dx + 2\lambda \int \biggl(\dfrac{v}{r} \biggr) r v v_{xt} \,dx - 8\mu \int \biggl(\dfrac{v}{r} \biggr) r v_x v_t\,dx.
	\end{aligned}
	\end{displaymath}
	Again, using \eqref{PrAssum} and applying Cauchy's inequality, it follows,
	\begin{displaymath}
	\begin{aligned}
		& L_1,L_2,L_3,L_4 \leq C (\alpha^{3\gamma}+\beta) \biggl\lbrace \int r v_t v_x \,dx + \int r_x v v_t\,dx \biggr.\\
		& \biggl. + \int r v v_{xt}\,dx + \int \dfrac{r^2 v_x v_{xt}}{r_x}\,dx  \biggr\rbrace 		\leq \delta \biggl\lbrace \int r^2 \dfrac{v_{xt}^2}{r_x}\,dx + \int r_x v_t^2 \,dx\biggr\rbrace \\
		& ~~~~~~ + C_{\delta} (\alpha^{6\gamma} + \beta^2) \biggl\lbrace \int r^2 \dfrac{v_x^2}{r_x} \,dx + \int r_x v^2\,dx \biggr\rbrace.
	\end{aligned}
	\end{displaymath}
	Summing up the inequality above, it follows from \eqref{ee003},
	\begin{equation}\label{ee0031}
		\begin{aligned}
			& \dfrac{d}{dt} \biggl\lbrace \dfrac{1}{2} \int x^2 \rho_0 v_t^2 \,dx \biggr\rbrace + \int 2\mu \biggl( r^2 \dfrac{v_{xt}^2}{r_x} + 2 r_x v_t^2 \biggr) + \lambda r^2 r_x \biggl( \dfrac{v_{xt}}{r_x} + 2 \dfrac{v_t}{r} \biggr)^2\,dx \\
			& \leq 5 \delta \biggl\lbrace \int \biggl( r^2 \dfrac{v_{xt}^2}{r_x} + 2 r_x v_t^2 \,dx\biggr) \biggr\rbrace + C_{\delta}(\alpha^{6\gamma} +\beta^2) \biggl\lbrace \int r^2 \dfrac{v_x^2}{r_x} \,dx + \int r_x v^2\,dx \biggr\rbrace.
		\end{aligned}
	\end{equation}
	We shall choose $ \delta $ small enough. Then 
	\eqref{lm2} follows from integration in temporal variable of \eqref{ee0031} together with \eqref{lm1}.
\end{pf}

\subsection{Uniform Estimates}\label{sec:uniformestimate}
The aim in this section is to show some uniform (in time) estimates. With these estimates, it would be able to design the restriction on initial data,  with which under the a prior assumption \eqref{PrAssum}, it can be shown,
\begin{equation}
	\dfrac{x^2}{r^2r_x} < \alpha^3, \biggl| \dfrac{v_x}{r_x} \biggr|, \biggl|\dfrac{v}{r}\biggr| < \beta
\end{equation}
for a smooth solution to \eqref{LgNS}. In particular, the a prior bounds in \eqref{PrAssum} can be justified. Moreover, the point-wise bounds of $$ \dfrac{x^2}{r^2r_x},~ \dfrac{v_x}{r_x}, ~ \dfrac{v}{r}  $$ are independent of time, which would be important ingredients to establish the propagation of regularity in the next section. In this and the following sections, we agree that the constant $ C > 0 $ is a universal constant which might be different from line to line and depends only on $ \bar \rho, \mu, \lambda, \gamma $.

A direct consequence of \eqref{PrAssum2} is $ R(t) = r(1,t) \geq \alpha^{-1} $. 
%With \eqref{lm303} in hand, 
With such a prior lower bound of the radius, the value of $ v $ on the boundary($ x=1 $) admits the following estimates.
\begin{lm}\label{lm:boundary}   %[$ \biggl.v \biggr|_{x=1}$]
Under the same assumptions as in Lemma \ref{lm:basicest}, there is a constant $ C>0 $ such that,
	\begin{align}
	& 	\int  \biggl.v^2 \biggr|_{x=1} \,dt \leq C \alpha E_{0}, \label{lm401} \\
	& \int \biggl. v_t^2 \biggr|_{x=1} \, dt \leq C \alpha E_1, \label{lm403} \\
	&  \biggl.v^2 \biggr|_{x=1} \leq C \alpha E_{0}  + C \alpha E_1. \label{lm402}
	\end{align}
\end{lm}

\begin{pf}
	As we have mentioned in \eqref{PrAssum2}, $ \exists 0 < \sigma <  \alpha^{-1}/2 $, satisfying $ R(t) - \sigma > \sigma $. For instant, we shall take $ \sigma = \alpha^{-1}/4 $. Then we have, by applying the mean value theorem, fundamental theorem of calculus and Cauchy's inequality,
	\begin{equation*}
	\begin{aligned}
		& \biggl.v^2\biggr|_{x=1} = u^2(R(t),t) \leq \sigma^{-1} \int_{R(t) - \sigma}^{R(t)} u^2\,dr + \int_{R(t) - \sigma}^{R(t)} \abs{(u^2)_r}{}\,dr\\
		& \leq C \sigma^{-1} \int_{R(t)-\sigma}^{R(t)} u^2 \,dr + C \sigma \int_{R(t)-\sigma}^{R(t)} u_r^2\,dr \leq C \sigma^{-1} \int_{R(t) - \sigma}^{R(t)} u^2\,dr \\& ~~~~ + C \sigma (R(t) - \sigma)^{-2}  \int_{R(t) - \sigma}^{R(t)} r^2 u_r^2\,dr. % = C  \int r_x v^2 \,dx + C \int r^2 \dfrac{v_x^2}{r_x} \,dx
	\end{aligned}
	\end{equation*}
	Therefore,
	\begin{equation}\label{uee001}
		\begin{aligned}
			& \biggl.v^2\biggr|_{x=1} \leq C \sigma^{-1} \biggl( \int_0^{R(t)} u^2 \,dr + \int_0^{R(t)} r^2 u_r^2 \,dr \biggr) \\
			& ~~~~~~~~ = C \alpha \biggl(\int r_x v^2 \,dx + \int r^2 \dfrac{v_x^2}{r_x} \,dx \biggr),
		\end{aligned}
	\end{equation}
	where we have used the fact $ (R(t)-\sigma)^{-2} < \sigma^{-2} $. 
	Then from \eqref{lm1} it holds
	\begin{equation}
		\int  \biggl.v^2\biggr|_{x=1}\,dt \leq C \alpha E_{0}.
	\end{equation}
	Meanwhile, by \eqref{lm2}, similar argument as in \eqref{uee001} yields
	\begin{equation}\label{uee0011}
	\begin{aligned}
		& \int \biggl. v_t^2\biggr|_{x=1} \,dt  \leq C \alpha \int \biggl( \int r_x v_t^2 \,dx  + \int  r^2 \dfrac{v_{xt}^2}{r_x} \,dx \biggr) \,dt  \leq C \alpha E_1.
	\end{aligned}
	\end{equation}
	Consequently, 
	\begin{equation}\label{uee002}
		\begin{aligned}
			& \biggl.v^2\biggr|_{x=1} \leq % \biggl. v^2\biggr|_{x=1,t=0} + 
			C \int \biggl( \biggl. v^2\biggr|_{x=1} + \biggl. v_t^2\biggr|_{x=1} \biggr) \,dt  \leq C \alpha E_{0}+ C \alpha E_1.
%			& ~~~~~~ + \int \int \biggl( r_x v_t^2 + r^2 \dfrac{v_{xt}^2}{r_x} \biggr)\,dx \,dt	 \leq \initial
%			 C \int \biggl(\int r_x v^2 \,dx + \int r^2 \dfrac{v_x^2}{r_x} \,dx  \biggr) \,dt\\
%			& ~~~~~~ + C \int \biggl( \int \dt(r_x v^2) \,dx +  \int \dt \biggl( r^2 \dfrac{v_x^2}{r_x} \biggr) \,dx  \biggr)\,dt
		\end{aligned}
	\end{equation}
%	from \eqref{lm1} and \eqref{lm2}.
%	Notice
%	\begin{displaymath}
%		\begin{aligned}
%			& \dt (r_x v^2 ) = v_x v^2 + 2 r_x v v_t = \dfrac{v_x}{r_x} r_x v^2 + 2 r_x v v_t\\
%			& \dt\biggl( r^2 \dfrac{v_x^2}{r_x} \biggr) = 2 r v \dfrac{v_x^2}{r_x} + 2 r^2 \dfrac{v_x v_{xt}}{r_x} - r^2 \dfrac{v_x^2}{r_x^2} v_x = \dfrac{v}{r} 2 r^2 \dfrac{v_x^2}{r_x} + 2 r^2 \dfrac{v_xv_{xt}}{r_x} - \dfrac{v_x}{r_x} r^2 \dfrac{v_x^2}{r_x}
%		\end{aligned}
%	\end{displaymath}
%	Then by applying H\"older inequality to \eqref{uee002}, \eqref{lm1} and \eqref{lm2} yield
%	\begin{equation}
%		\begin{aligned}
%			& \biggl.v^2\biggr|_{x=1} \leq C \int \biggl(\int r_x v^2 \,dx + \int r^2 \dfrac{v_x^2}{r_x} \,dx  \biggr) \,dt\\
%			& ~~~~~~ + C \int \biggl(\int r_x v_t^2 \,dx + \int r^2 \dfrac{v_{xt}^2}{r_x} \,dx  \biggr) \,dt \leq \initial.
%		\end{aligned}
%	\end{equation}

\end{pf}

%============

The following lemma is the key ingredient in this work. It shows the time-integrability of $$ \biggl(\dfrac{x^2\rho_0}{r^2r_x}\biggr)^{2\gamma}, \biggl|\dfrac{v_x}{r_x}\biggr|^2, \biggl|\dfrac{v}{r}\biggr|^2. $$

\begin{lm} Under the same assumptions as in Lemma \ref{lm:basicest}, there exists a constant   $ C > 0 $  % depending on $ \mu, \lambda, \eta, \bar\rho $ 
such that
	\begin{align}
		& \biggl\Arrowvert\biggl( \dfrac{x^2 \rho_0}{r^2r_x} \biggr)^\gamma \biggr\Arrowvert_{L^\infty_x}  + C \int\biggl\Arrowvert\biggl( \dfrac{x^2 \rho_0}{r^2r_x} \biggr)^{2\gamma} \biggr\Arrowvert_{L^\infty_x} \,dt + C \int \biggl\Arrowvert \dfrac{v_x}{r_x} + 2 \dfrac{v}{r} \biggr\Arrowvert_{L^\infty_x}^2 \,dt {\nonumber} \\
		& ~~~~~~~~ \leq \bar\rho^\gamma + C \alpha^5 E_1 + C \alpha^3 E_{0}, \label{lm501} \\
		& \int \biggl\Arrowvert\dfrac{v_x}{r_x}-\dfrac{v}{r} \biggr\Arrowvert_{L^\infty_x}^2 \,dt \leq C \bar\rho^\gamma + C \alpha^{5} E_1 + C \alpha^3 E_{0}, \label{lm502}\\
		& \int \biggl\Arrowvert \dfrac{v_x}{r_x} \biggr\Arrowvert_{L^\infty_x}^2 + \biggl\Arrowvert \dfrac{v}{r} \biggr\Arrowvert_{L^\infty_x}^2 \,dt \leq  C \bar\rho^\gamma + C \alpha^{5} E_1 + C \alpha^3 E_{0}.   \label{lm503}
	\end{align}
\end{lm}

\begin{pf}
Integrate \eqref{LgNS} over $ (x,1) $ in spatial variable for $ 0 < x < 1 $. It follows,
\begin{equation}\label{uee003}
	- (2 \mu + \lambda ) \dfrac{(r^2 v)_x}{r^2 r_x} + \biggl( \dfrac{x^2 \rho_0}{r^2 r_x} \biggr)^\gamma = \int_x^1 \biggl(\dfrac{x}{r}\biggr)^2 \rho_0 v_t\,dx - 4\mu \biggl.\dfrac{v}{r} \biggr|_{x=1}.
\end{equation}
Notice
$$ (r^2 v)_x = \dfrac{1}{3} (\dt r^3 )_x = ( r^2 r_x )_t. $$
Taking square on both sides of \eqref{uee003} yields, 
\begin{equation}\label{uee004}
	\begin{aligned}
		& - 2 ( 2 \mu + \lambda ) \dfrac{(r^2r_x)_t}{r^2r_x} \cdot \biggl( \dfrac{x^2 \rho_0}{r^2r_x} \biggr)^\gamma + \biggl( \dfrac{x^2 \rho_0}{r^2r_x} \biggr)^{2\gamma} + \biggl(( 2 \mu + \lambda ) \dfrac{(r^2r_x)_t}{r^2r_x} \biggr)^2\\
		& ~~~~~~ \leq 2 \biggl( \int_x^1 \biggl(\dfrac{x}{r}\biggr)^2 \rho_0 v_t\,dx\biggr)^2 + 32 \mu^2 \biggl( \biggl.\dfrac{v}{r}\biggr|_{x=1} \biggr)^2 \\
		& ~~~~~~ \leq 2  \int_x^1 \biggl(\dfrac{x}{r}\biggr)^2 \rho_0 v_t^2 \,dx \cdot \int_x^1 \biggl(\dfrac{x}{r}\biggr)^2 \rho_0 \,dx  + C \biggl. \dfrac{v^2}{r^2} \biggr|_{x=1}.
	\end{aligned}
\end{equation}
Notice, from \eqref{PrAssum}, \eqref{PrAssum2} and \eqref{InitialAssum}
\begin{displaymath}
\begin{aligned}
	& - \dfrac{(r^2r_x)_t}{r^2r_x} \cdot \biggl( \dfrac{x^2 \rho_0}{r^2r_x} \biggr)^\gamma = \dfrac{1}{\gamma}\dfrac{\partial}{\partial t}\biggl( \dfrac{x^2\rho_0}{r^2r_x}\biggr)^\gamma, \\
	& % \biggl( \int_x^1 \biggl(\dfrac{x}{r}\biggr)^2 \rho_0 v_t\,dx\biggr)^2 \leq
	 \int_x^1 \biggl(\dfrac{x}{r}\biggr)^2 \rho_0 v_t^2 \,dx \cdot \int_x^1 \biggl(\dfrac{x}{r}\biggr)^2 \rho_0 \,dx \leq C \alpha^5 \int r_x v_t^2\,dx.
\end{aligned}
\end{displaymath}
%where we have used, from \eqref{InitialAssum}, \eqref{PrAssum}, \eqref{PrAssum2} in the last inequality.
%$ \biggl(\frac{x}{r}\biggr)^2 = \dfrac{x^2}{r^2 r_x} r_x \leq C r_x $ and $ \frac{x}{r} \leq C $.\\
Integration over $ ( 0, t) $ in temporal variable of \eqref{uee004} then implies
\begin{equation}\label{06-July-001}
	\begin{aligned}
		& \norm{\biggl( \dfrac{x^2 \rho_0}{r^2r_x} \biggr)^\gamma}{L^\infty_x} + C  \int\biggl( \dfrac{x^2 \rho_0}{r^2r_x} \biggr)^{2\gamma} \,dt + C \int \biggl(\dfrac{(r^2r_x)_t}{r^2r_x} \biggr)^2 \,dt  \\
		& ~~~~~~ \leq \rho_0^\gamma + C \alpha^5 \int \,dt \int r_x v_t^2\,dx + C \int \biggl.\dfrac{v^2}{r^2}\biggr|_{x=1} \,dt \\
		& ~~~~~~ \leq\bar\rho^\gamma + C \alpha^5 E_1 + C \alpha^3 E_{0},
	\end{aligned}
\end{equation}
where the last inequality follows from \eqref{PrAssum2}, \eqref{lm2} and \eqref{lm401}%, \eqref{lm303}
.  %The estimate of $  \norm{\biggl( \dfrac{x^2 \rho_0}{r^2r_x} \biggr)^\gamma}{L^\infty_x} $ in \eqref{lm501} follows. 
On the other hand, considering \eqref{uee004} at the point $ (x_1(t),t) $ where $$ \biggl(\dfrac{x^2 \rho_0}{r^2 r_x}\biggr) (x_1(t),t) = \norm{\dfrac{x^2 \rho_0}{r^2 r_x}}{L^\infty_x}(t). $$ Then, with $ \dfrac{d}{dt} := \partial_t + \dt x_1(t) \partial_x $ representing the trajectory differential along $ (x_1(t),t) $, we have
\begin{gather*}
- \dfrac{(r^2r_x)_t}{r^2r_x} \cdot \biggl( \dfrac{x^2 \rho_0}{r^2r_x} \biggr)^\gamma\bigg|_{(x_1(t),t)} = \dfrac{1}{\gamma}\dfrac{d}{dt} \biggl( \dfrac{x^2\rho_0}{r^2r_x}\biggr)^\gamma \bigg|_{(x_1(t),t)} \\
 - \dfrac{1}{\gamma} \dfrac{\partial}{\partial x}\biggl( \dfrac{x^2\rho_0}{r^2r_x}\biggr)^\gamma \bigg|_{(x_1(t),t)} \times \dt x_1(t) = \dfrac{1}{\gamma}\dfrac{d}{dt} \biggl( \dfrac{x^2\rho_0}{r^2r_x}\biggr)^\gamma \bigg|_{(x_1(t),t)}.
\end{gather*}
Thus integrating \eqref{uee004} in temporal variable along $ (x_1(t),t) $ yields, with the rest estimated similarly as above,
\begin{equation}\label{06-July-002}
	\begin{gathered}
	  \int \biggl\Arrowvert \biggl( \dfrac{x^2 \rho_0}{r^2r_x} \biggr)^{2\gamma} \biggr\Arrowvert_{L^\infty_x} \,dt \leq \dfrac{2(2\mu+\lambda)}{\gamma}\biggl( \dfrac{x^2\rho_0}{r^2r_x}\biggr)^\gamma \bigg|_{(x_1(t),t)}  + \int \biggl\Arrowvert \biggl( \dfrac{x^2 \rho_0}{r^2r_x} \biggr)^{2\gamma} \biggr\Arrowvert_{L^\infty_x} \,dt \\
	  \leq  \bar\rho^\gamma + C \alpha^5 E_1 + C \alpha^3 E_{0}.
	  \end{gathered}
\end{equation}
Moreover, considering \eqref{uee004} at the point $ (x_2(t),t) $ where $$ \biggl(\dfrac{(r^2r_x)_t}{r^2r_x} \biggr)(x_2(t),t) = \norm{\dfrac{(r^2r_x)_t}{r^2r_x} }{L^\infty}(t). $$ Applying Cauchy's inequality to the resultant inequality yields
\begin{gather*}
	\dfrac{(2\mu+\lambda)^2}{2} \biggl(\dfrac{(r^2r_x)_t}{r^2r_x} \biggr)^2(x_2(t),t) \leq  C  \int_{x_2(t)}^1 \biggl(\dfrac{x}{r}\biggr)^2 \rho_0 v_t^2 \,dx \cdot \int_{x_2(t)}^1 \biggl(\dfrac{x}{r}\biggr)^2 \rho_0 \,dx \\
	  + C \biggl. \dfrac{v^2}{r^2} \biggr|_{x=1} + \norm{\biggl( \dfrac{x^2 \rho_0}{r^2r_x} \biggr)^{2\gamma}}{L^\infty_x}.
\end{gather*}
Integrating in the temporal variable along $ (x_2(t),t) $ yields,  with the right hand side estimated similarly as above, 
\begin{equation}\label{06-July-003}
	\int \norm{\dfrac{(r^2r_x)_t}{r^2r_x}}{L^\infty_x}^2 \,dt \leq C \int \norm{\biggl( \dfrac{x^2 \rho_0}{r^2r_x} \biggr)^{2\gamma}}{L^\infty_x} \,dt + C \bar\rho^\gamma + C \alpha^5 E_1 + C \alpha^3 E_{0}.
\end{equation}
After noticing
\begin{displaymath}
	\dfrac{(r^2r_x)_t}{r^2r_x} = \dfrac{v_x}{r_x} + 2\dfrac{v}{r},
\end{displaymath}
\eqref{lm501} follows from \eqref{06-July-001}, \eqref{06-July-002}, and \eqref{06-July-003}.

On the other hand, multiply \eqref{LgNS} with $ r^3 $ and integrate the resulting over $ (0,x) $ in spatial variable for $ 0 < x < 1 $. After integration by parts, it holds, 
\begin{equation}\label{uee005}
\begin{aligned}
	& \int_0^x x^2 r \rho_0 v_t \,dx + \biggl(\dfrac{x^2\rho_0}{r^2r_x}\biggr)^\gamma r^3 - \int_0^x\biggl(\dfrac{x^2\rho_0}{r^2r_x}\biggr)^\gamma (r^3)_x\,dx\\
	& ~~~~~~ = (2\mu+\lambda) \biggl( \dfrac{(r^2v)_x}{r^2r_x} r^3 -  \int_0^x \dfrac{(r^2v)_x}{r^2r_x} (r^3)_x\,dx \biggr).
\end{aligned}
\end{equation}
Direct calculation yields,
\begin{displaymath}
	\begin{aligned}
	& \dfrac{(r^2v)_x}{r^2r_x} r^3 -  \int_0^x \dfrac{(r^2v)_x}{r^2r_x} (r^3)_x\,dx = \biggl( \dfrac{v_x}{r_x} - \dfrac{v}{r} \biggr) \cdot r^3.\\
	\end{aligned}
\end{displaymath}
Moreover, similar as before, from \eqref{InitialAssum} and \eqref{PrAssum2}
\begin{displaymath}
	\begin{aligned}
	& \int_0^x x^2 r \rho_0 v_t \,dx \leq \biggl( \int\biggl(\dfrac{x}{r}\biggr)^2 \rho_0 v_t^2\,dx \biggr)^{1/2} \biggl(\int_0^x \rho_0 x^2 r^4 \,dx\biggr)^{1/2} \\
	& ~~~~\leq C  \biggl( \int \biggl(\dfrac{x}{r}\biggr)^2 \rho_0 v_t^2\,dx \biggr)^{1/2} \biggl( x \dfrac{x^2}{r^2} r^6 \biggr)^{1/2} \\
	& ~~~~\leq C \alpha \biggl( \int\biggl(\dfrac{x}{r}\biggr)^2 \rho_0 v_t^2\,dx \biggr)^{1/2} r^3 , \\
	& \int_0^x\biggl(\dfrac{x^2\rho_0}{r^2r_x}\biggr)^\gamma (r^3)_x\,dx \leq \biggl\Arrowvert\biggl( \dfrac{x^2\rho_0}{r^2r_x} \biggr)^\gamma\biggr\Arrowvert_{L^\infty_x} \int_0^x (r^3)_x \,dx \leq \biggl\Arrowvert\biggl( \dfrac{x^2\rho_0}{r^2r_x} \biggr)^\gamma\biggr\Arrowvert_{L^\infty_x} r^3.
	\end{aligned}
\end{displaymath}
Therefore, \eqref{uee005} implies
\begin{equation*}
	\biggl| \dfrac{v_x}{r_x} - \dfrac{v}{r} \biggr| \cdot r^3 \leq C \alpha \biggl( \int\biggl(\dfrac{x}{r}\biggr)^2 \rho_0 v_t^2\,dx \biggr)^{1/2}  r^{3} + C \biggl\Arrowvert\biggl( \dfrac{x^2\rho_0}{r^2r_x} \biggr)^\gamma\biggr\Arrowvert_{L^\infty_x} r^3,
\end{equation*}
or
\begin{equation}\label{uee0051}
	\begin{aligned}
	& \biggl| \dfrac{v_x}{r_x} - \dfrac{v}{r} \biggr|^2 \leq C \alpha^2  \int\biggl(\dfrac{x}{r}\biggr)^2 \rho_0 v_t^2\,dx + C \biggl\Arrowvert\biggl( \dfrac{x^2\rho_0}{r^2r_x} \biggr)^{2\gamma}\biggr\Arrowvert_{L^\infty_x} \\
	& ~~~~~~ \leq C \alpha^{5} \int r_x v_t^2 \,dx + C \biggl\Arrowvert\biggl( \dfrac{x^2\rho_0}{r^2r_x} \biggr)^{2\gamma}\biggr\Arrowvert_{L^\infty_x},
	\end{aligned}
\end{equation}
where \eqref{InitialAssum} and \eqref{PrAssum} is applied in the last inequality. Integrating the above inequality over temporal variable, together with \eqref{lm2} and \eqref{lm501} then yields \eqref{lm502}.
\eqref{lm503} follows from \eqref{lm501} and \eqref{lm502}.
\end{pf}

%===========
We shall need the following lemmas concerning the interior energy inspired by \cite{LuoXinZeng2016}, which will claim extra integrability of the unknowns. 

\begin{lm}
Under the same assumptions as in Lemma \ref{lm:basicest}, there exists a constant $ C > 0 $ % depending on $ \mu, \lambda, \eta, \bar\rho $ 
such that,
\begin{equation}\label{lm6}
\begin{aligned}
& \int \dfrac{v_x^2}{r_x}\,dx + \int r_x \biggl(\dfrac{v}{r}\biggr)^2\,dx  \leq  C \alpha^2 E_0 + C( \alpha^2 + \Gamma) E_1 + C \alpha^2 \Gamma^{2\gamma-1},
\end{aligned}	
\end{equation}
where 
\begin{equation}\label{lm602}
	\Gamma =: \min \biggl\lbrace \bar\rho \alpha^3, (\bar\rho^\gamma + C \alpha^5 E_1 + C \alpha^3 E_0)^{1/\gamma}\biggr\rbrace.
\end{equation}
As a consequence of \eqref{PrAssum}, \eqref{InitialAssum} and \eqref{lm501},
\begin{equation}\label{lm601}
	\max_{0\leq x \leq 1} \dfrac{x^2\rho_0}{r^2r_x} \leq \Gamma.
\end{equation}

\end{lm}

\begin{pf}

Multiply \eqref{LgNS} with $ v $ and integrate the resulting in spatial variable. Then it holds
\begin{equation}\label{uee006}
	\int \biggl(\dfrac{x}{r}\biggr)^2 \rho_0 v_t v\,dx - \int \biggl(\dfrac{x^2 \rho_0}{r^2 r_x} \biggr)^\gamma v_x \,dx = - (2\mu + \lambda)\int \dfrac{(r^2 v)_x}{r^2 r_x } v_x\,dx + 4 \mu \biggl.\dfrac{v^2}{r}\biggr|_{x=1}.
\end{equation}
where the integration by parts is applied with the boundary condition \eqref{LgNSBC}. 
Notice, integration by parts again yields, %integration by parts yields, by \eqref{LgNSBC}
\begin{displaymath}
	\begin{aligned}
		& (2\mu + \lambda)\int \dfrac{(r^2 v)_x}{r^2 r_x } v_x\,dx = (2\mu + \lambda) \int \dfrac{v_x^2}{r_x}\,dx + (2\mu + \lambda) \int \dfrac{2 v v_x}{r}\,dx \\
		& ~~~~~~ = (2\mu + \lambda) \int \dfrac{v_x^2}{r_x}\,dx + (2\mu + \lambda) \int r_x\dfrac{v^2}{r^2}\,dx + (2\mu + \lambda) \biggl. \dfrac{v^2}{r} \biggr|_{x=1}.
	\end{aligned}
\end{displaymath}
Moreover, from \eqref{InitialAssum}, \eqref{PrAssum2}, \eqref{lm601} and Cauchy's inequality,
\begin{displaymath}
	\begin{aligned}
		& \int \biggl(\dfrac{x}{r}\biggr)^2 \rho_0 v_t v\,dx \leq \delta\int r_x \biggl(\dfrac{v}{r}\biggr)^2\,dx + C_\delta \int \dfrac{x^2\rho_0}{r^2r_x} x^2 \rho_0 v_t^2 \,dx \\
		& ~~~~~~ \leq \delta \int r_x \biggl(\dfrac{v}{r}\biggr)^2\,dx + C_{\delta} \Gamma \int x^2 \rho_0 v_t^2 \,dx,\\
		& \int \biggl(\dfrac{x^2 \rho_0}{r^2 r_x} \biggr)^\gamma v_x \,dx  \leq \delta \int \dfrac{v_x^2}{r_x}\,dx + C_\delta \int \dfrac{x^2 \rho_0}{r^2}     \biggl(\dfrac{x^2 \rho_0}{r^2 r_x} \biggr)^{2\gamma-1} \,dx\\
		& ~~~~~~ \leq \delta \int \dfrac{v_x^2}{r_x}\,dx + C_\delta \alpha^{2} \Gamma^{2\gamma-1}.
	\end{aligned}
\end{displaymath}
Therefore, from \eqref{uee006} it follows
\begin{equation}
\begin{aligned}
	& (2\mu + \lambda) \biggl( \int \dfrac{v_x^2}{r_x}\,dx + \int r_x \biggl(\dfrac{v}{r}\biggr)^2\,dx \biggr)  \leq 2 \delta \biggl( \int r_x\biggl(\dfrac{v}{r}\biggr)^2\,dx + \int \dfrac{v_x^2}{r_x}\,dx \biggr) \\& ~~+ C_{\delta}\Gamma \int x^2 \rho_0 v_t^2 \,dx + C_\delta \alpha^{2} \Gamma^{2\gamma-1} + C \biggl. \dfrac{v^2}{r}\biggr|_{x=1}\\
	& \leq 2 \delta \biggl( \int r_x\biggl(\dfrac{v}{r}\biggr)^2\,dx + \int \dfrac{v_x^2}{r_x}\,dx \biggr) + C_{\delta}\Gamma E_1 + C_\delta \alpha^2 \Gamma^{2\gamma-1} \\
	& ~~~~~~~ + C \alpha^2 E_{0} + C \alpha^2 E_1,  
\end{aligned}
\end{equation}
where the last inequality follows from \eqref{PrAssum2}, \eqref{lm2}, \eqref{lm402}.
By choosing $ \delta $ small enough, \eqref{lm6} follows.
\end{pf}

\begin{lm}
%For $ E_0, E_1, E_2 $ small enough
%There exists a constant $ C > 0 $ % depending on $ \mu, \lambda, \eta, \bar\rho $ 
%such that
Under the same assumptions as in Lemma \ref{lm:basicest}, there exists a constant $ C > 0 $ such that, 
	\begin{equation}\label{lm801}
	\begin{aligned}
	& \dfrac{1}{2}\int \biggl(\dfrac{x}{r} \biggr)^2\rho_0 v_t^2\,dx + \dfrac{2\mu+\lambda}{2} \int  \int \biggl( r_x \dfrac{v_t^2}{r^2} + \dfrac{v_{xt}^2}{r_x}\biggr)\,dx \,dt \\
	& \leq E_2 + C (\alpha^2 E_0 + (\alpha^2 +\Gamma)E_1 + \alpha^{2}\Gamma^{2\gamma-1})\times ( \alpha^3 E_0 + \alpha^5 E_1 + 1). 
	\end{aligned}
	\end{equation}
\end{lm}

\begin{pf}
	Taking the time derivative of \eqref{LgNS},
	\begin{equation} \label{LgNS01}
		\biggl( \dfrac{x}{r} \biggr)^2\rho_0 v_{tt} + \biggl\lbrack \biggl( \dfrac{x^2\rho_0}{r^2r_x} \biggr)^\gamma \biggr\rbrack_{xt} = (2\mu+\lambda) \biggl\lbrack \dfrac{(r^2 v)_x}{r^2r_x} \biggr\rbrack_{xt} + 2 \dfrac{x^2}{r^2}\dfrac{v}{r} \rho_0 v_t.
	\end{equation}
%	Define an interior cut-off function $ \chi \in C^{\infty}(0,1) $ satisfying
%	\begin{equation}
%		\chi(x) = \begin{cases}
%			1, & 0 < x < 0.5 \\
%			0, & 0.75 < x < 1 \\
%			\text{decreasing}, & 0.5 < x < 0.75
%		\end{cases}
%	\end{equation}
%	with 
%	\begin{equation}
%		\chi'(x) \leq 6.
%	\end{equation}
	Multiply this equation with $ v_t $ and integrate the resulting in spatial variable. 
	\begin{equation}
		\begin{aligned}
		& \int \biggl(\dfrac{x}{r}\biggr)^2 \rho_0 v_{tt} v_t\,dx - \int v_{xt} \biggl\lbrack \biggl( \dfrac{x^2\rho_0}{r^2r_x}\biggr)^\gamma \biggr\rbrack_t \,dx \\
		& ~~~~~~ = - (2\mu + \lambda) \int v_{xt} \biggl\lbrack \dfrac{(r^2v)_x}{r^2r_x} \biggr\rbrack_t \,dx + 4 \mu \biggl. v_t \biggl( \dfrac{v}{r} \biggr)_{t} \biggr|_{x=1} + 2 \int\dfrac{x^2}{r^2}\dfrac{v}{r} \rho_0 v_t^2\,dx.
		\end{aligned}
	\end{equation}
	Direct calculation yields,
	\begin{equation}\label{uee009}
		\begin{aligned}
			& \dfrac{d}{dt} \biggl\lbrace \dfrac{1}{2} \int \biggl( \dfrac{x}{r}\biggr)^2 \rho_0 v_t^2 \biggr\rbrace  + (2\mu+\lambda) \int \biggl( r_x \dfrac{v_t^2}{r^2} + \dfrac{v_{xt}^2}{r_x} \biggr)\,dx \\
			& ~~~~~~ = \int \dfrac{x^2\rho_0}{r^2}\dfrac{v}{r} v_t^2\,dx + (2\mu + \lambda) \int \biggl(2 \dfrac{v^2}{r^2} + \dfrac{v_x^2}{r_x^2}\biggr) v_{xt}\,dx\\
%			& ~~~~~~~~ - (2\mu + \lambda) \int \chi' \biggl( \dfrac{v_t^2}{r} + \dfrac{v_{xt}v_t}{r_x} - 2 \dfrac{v^2 v_t}{r^2} - \dfrac{v_x^2 v_t}{r_x^2} \biggr) \,dx\\
			& ~~~~~~~~ - \gamma \int \biggl( \dfrac{x^2\rho_0}{r^2 r_x} \biggr)^\gamma \biggl( \dfrac{v_x v_{xt}}{r_x} + 2 \dfrac{v v_{xt}}{r} \biggr) \,dx \\
			& ~~~~~~~~  + \biggl. \biggl\lbrack-(2\mu+\lambda) \dfrac{v_t^2}{r} + 4\mu v_t \biggl(\dfrac{v}{r}\biggr)_t\biggr\rbrack\biggr|_{x=1} : = I_1 + I_2 + I_3 + I_4.
		\end{aligned}
	\end{equation}
	
	We estimate $ I_1,  \cdots,  I_4 $ as follows. By noticing \eqref{PrAssum2}, \eqref{lm402} and \eqref{lm6}, Cauchy's inequality yields
	\begin{displaymath}
		\begin{aligned}
%			& I_1 \leq C {AB} \int r_x v_t^2 \,dx \\
			& I_2 \leq \delta \int \dfrac{v_{xt}^2}{r_x}\,dx + C_\delta \biggl\Arrowvert \dfrac{v_x}{r_x} \biggr\Arrowvert_{L^\infty_x}^2 \int \dfrac{v_x^2}{r_x}\,dx + C_\delta \biggl\Arrowvert \dfrac{v}{r} \biggr\Arrowvert_{L^\infty_x}^2\int r_x \dfrac{v^2}{r^2} \,dx\\
			& ~~~~ \leq \delta \int \dfrac{v_{xt}^2}{r_x}\,dx + C_\delta ( \alpha^2 E_0 + (\alpha^2 +\Gamma)E_1 + \alpha^{2}\Gamma^{2\gamma-1}) \\ & ~~~~~~~~ \times \biggl\lbrace \biggl\Arrowvert \dfrac{v_x}{r_x} \biggr\Arrowvert_{L^\infty_x}^2 + \biggl\Arrowvert \dfrac{v}{r} \biggr\Arrowvert_{L^\infty_x}^2 \biggr\rbrace, \\
			& I_3 \leq \delta \int \dfrac{v_{xt}^2}{r_x}\,dx + C_\delta \biggl\Arrowvert \biggl(\dfrac{x^2\rho_0}{r^2r_x} \biggr)^{2\gamma} \biggr\Arrowvert_{L^\infty_x} \cdot \biggl\lbrace \int \dfrac{v_x^2}{r_x}\,dx + \int r_x \dfrac{v^2}{r^2}\,dx \biggr\rbrace \\
			& ~~~~ \leq \delta \int \dfrac{v_{xt}^2}{r_x}\,dx + C_\delta (\alpha^2 E_0 + (\alpha^2 +\Gamma)E_1 + \alpha^{2}\Gamma^{2\gamma-1}) \cdot  \biggl\Arrowvert \biggl(\dfrac{x^2\rho_0}{r^2r_x} \biggr)^{2\gamma} \biggr\Arrowvert_{L^\infty_x} ,\\
			& I_4 \leq C \biggl.( \alpha v_t^2 + \alpha^3 v^4)\biggr|_{x=1} \leq  C  \alpha \biggl. v_t^2\biggr|_{x=1} + C \alpha^4 (E_0+E_1)\biggl. v^2\biggr|_{x=1} .
			\end{aligned}
	\end{displaymath}
Moreover, after plugging \eqref{LgNS} into $ I_1 $, applying integration by parts to the resulting expression then yields, together with \eqref{PrAssum2},
\begin{displaymath}
	\begin{aligned}
		& I_1 = \int \dfrac{v}{r}v_t \biggl( (2\mu+\lambda) \biggl(\dfrac{(r^2v)_x}{r^2r_x} \biggr)_x - \biggl\lbrack\biggl(\dfrac{x^2\rho_0}{r^2r_x}\biggr)^\gamma\biggr\rbrack_x \biggr)\,dx\\
		& ~~~~ = \int \biggl( \dfrac{v}{r} v_t \biggr)_x \biggl( \biggl(\dfrac{x^2\rho_0}{r^2r_x}\biggr)^\gamma - (2\mu + \lambda) \dfrac{(r^2 v)_x}{r^2r_x} \biggr) \,dx + 4\mu \biggl.\dfrac{v^2}{r^2}v_t \biggr|_{x=1}\\
		& \leq \delta \biggl\lbrace \int r_x \dfrac{v_t^2}{r^2}\,dx + \int \dfrac{v_{xt}^2}{r_x}\, dx \biggr\rbrace + C_{\delta} \biggl\lbrace \biggl\Arrowvert \biggl( \dfrac{x^2\rho_0}{r^2r_x}\biggr)^{2\gamma} \biggr\Arrowvert_{L_x^\infty} \biggr.\\ & ~~ \biggl. + \biggl\Arrowvert \dfrac{v_x}{r_x} \biggr\Arrowvert_{L^\infty_x}^2 + \biggl\Arrowvert \dfrac{v}{r} \biggr\Arrowvert_{L^\infty_x}^2  \biggr\rbrace  \cdot   \biggl\lbrace \int \dfrac{v_x^2}{r_x}\,dx + \int r_x\dfrac{v^2}{r^2}\,dx  \biggr\rbrace  + C  ( \alpha \biggl.v_t^2 + \alpha^3 v^4 )\biggr|_{x=1}.
	\end{aligned}
\end{displaymath}
Therefore, as consequences of \eqref{lm402} and \eqref{lm6},
\begin{displaymath} \begin{aligned}
		& I_1 \leq \delta \biggl\lbrace  \int r_x \dfrac{v_t^2}{r^2}\,dx + \int \dfrac{v_{xt}^2}{r_x}\, dx \biggr\rbrace + C_{\delta}( \alpha^2 E_0 + ( \alpha^2 +\Gamma)E_1 + \alpha^{2}\Gamma^{2\gamma-1}) \\& ~~~~~~ \times \biggl\lbrace \biggl\Arrowvert \biggl( \dfrac{x^2\rho_0}{r^2r_x}\biggr)^{2\gamma} \biggr\Arrowvert_{L_x^\infty} + \biggl\Arrowvert \dfrac{v_x}{r_x} \biggr\Arrowvert_{L^\infty_x}^2 + \biggl\Arrowvert \dfrac{v}{r} \biggr\Arrowvert_{L^\infty_x}^2  \biggr\rbrace\\
		& ~~~~~~ + C \alpha v_t^2 \bigr|_{x=1} + C \alpha^4 (E_0 + E_1) v^2\bigr|_{x=1}.
	\end{aligned}
\end{displaymath}

Consequently, by choosing $ \delta $ small enough, integration of \eqref{uee009} in temporal variable yields, 
\begin{equation}
	\begin{aligned}
		& \dfrac{1}{2} \int \biggl(\dfrac{x}{r}\biggr)^2 \rho_0 v_t^2\,dx + \dfrac{2\mu+\lambda}{2} \int  \int \biggl(r_x \dfrac{v_t^2}{r^2}+ \dfrac{v_{xt}^2}{r_x} \biggr)\,dx \,dt \\
		& \leq \dfrac{1}{2} \int \rho_0 u_1^2\,dx + C (\alpha^2 E_0 + (\alpha^2 +\Gamma)E_1 + \alpha^{2}\Gamma^{2\gamma-1}) \\&  ~~~~~~ \times \int \biggl\lbrace \biggl\Arrowvert \biggl( \dfrac{x^2\rho_0}{r^2r_x}\biggr)^{2\gamma} \biggr\Arrowvert_{L_x^\infty} + \biggl\Arrowvert \dfrac{v_x}{r_x} \biggr\Arrowvert_{L^\infty_x}^2 + \biggl\Arrowvert \dfrac{v}{r} \biggr\Arrowvert_{L^\infty_x}^2  \biggr\rbrace \,dt \\
		& ~~~~~~  + C \alpha \int \biggl. v_t^2 \biggr|_{x=1} \,dt + C \alpha^4 (E_0 + E_1) \int\biggl.v^2\biggr|_{x=1} \,dt \\
		& \leq C \alpha^5 E_0 (E_0+E_1)+ C \alpha^2 E_1 + E_2 \\
		& ~~~~ + C (\alpha^2 E_0 + (\alpha^2 +\Gamma)E_1 + \alpha^{2}\Gamma^{2\gamma-1})\times ( \alpha^3 E_0 + \alpha^5 E_1 + \bar\rho^\gamma),
	\end{aligned}
\end{equation}
where the last inequality is from \eqref{lm401}, \eqref{lm403},  \eqref{lm501}  and \eqref{lm503}. Then \eqref{lm801} follows by noticing $ \alpha > 1 $. 
\end{pf}

% ========
Now we have enough materials to show the point-wise bounds of $$ \dfrac{x^2}{r^2r_x}, \dfrac{v_x}{r_x}, \dfrac{v}{r}, $$ under the a prior assumption \eqref{PrAssum}. In particular, the bounds are independent of time. More precisely, 
\begin{lm}\label{lm:pointwisebound} Under the same assumptions as in Lemma \ref{lm:basicest}, there exists $ C > 0 $ % depending on $ \mu,\lambda,\eta,\bar\rho $ 
such that,
\begin{align}
	& \biggl\Arrowvert \dfrac{v_x}{r_x} + 2 \dfrac{v}{r} \biggr\Arrowvert_{L_x^\infty}^2 \leq C \alpha^2 ( ( \alpha^2 E_0 + (\alpha^2 +\Gamma)E_1 + \alpha^{2}\Gamma^{2\gamma-1}) {\nonumber} \\
	& ~~~~~~~\times ( \alpha^3 E_0 + \alpha^5 E_1 + 1) + E_2 ) + \Gamma^{2\gamma}, \label{lm901} \\
	& \biggl\Arrowvert  \dfrac{v_x}{r_x} - \dfrac{v}{r}   \biggr\Arrowvert_{L_x^\infty}^2  \leq C \alpha^2 ( ( \alpha^2 E_0 + (\alpha^2 +\Gamma)E_1 + \alpha^{2}\Gamma^{2\gamma-1}) {\nonumber} \\
	& ~~~~~~~\times ( \alpha^3 E_0 + \alpha^5 E_1 + 1) + E_2 ) + \Gamma^{2\gamma}, \label{lm902}\\
	& \biggl\Arrowvert \dfrac{x^2}{r^2r_x} \biggr\Arrowvert_{L_x^\infty} \leq \exp \bigl\lbrace C\alpha^{5/2} ( (1 + \beta )E_0 + E_1 )^{1/2}  {\nonumber} \\
	& ~~~~~~~  + C \alpha^{5/2} ( \alpha^3 E_0 + \alpha^5 E_1 + 1)^{1/2} E_1^{1/2}\bigr\rbrace.  \label{lm9021}
\end{align}

\end{lm}

\begin{pf}
	From \eqref{uee004}, \eqref{InitialAssum}, \eqref{PrAssum2}, 
	\begin{equation}\label{ueept001}
		\begin{aligned}
			& (2\mu+\lambda)^2 \biggl( \dfrac{(r^2r_x)_t}{r^2r_x} \biggr)^2 \leq \int_x^1 \biggl(\dfrac{x}{r}\biggr)^2 \rho_0 v_t^2 \,dx \cdot \int_x^1 \biggl(\dfrac{x}{r}\biggr)^2 \rho_0 \,dx  + C \biggl. \dfrac{v^2}{r^2} \biggr|_{x=1} \\% \biggl( \int\biggl(\dfrac{x}{r}\biggr)^2 \rho_0v_t\,dx\biggr)^2 + 8 \mu^2 \biggl(\biggl. \dfrac{v}{r}\biggr|_{x=1} \biggr)^2\\
			& ~~~~~~~~ - \biggl( \dfrac{x^2\rho_0}{r^2r_x} \biggr)^{2\gamma} + 2 (2\mu +\lambda) \dfrac{(r^2 r_x)_t}{r^2r_x}\cdot \biggl( \dfrac{x^2\rho_0}{r^2 r_x} \biggr)^\gamma\\
			& ~~~~~~ \leq \delta \biggl( \dfrac{(r^2r_x)_t}{r^2r_x} \biggr)^2 + C \alpha^2  \int \biggl(\dfrac{x}{r}\biggr)^2\rho_0 v_t^2 \,dx  + C \alpha^2 \biggl.v^2\biggr|_{x=1} + C_\delta \biggl( \dfrac{x^2\rho_0}{r^2r_x} \biggr)^{2\gamma}.
		\end{aligned}
	\end{equation}
	After choosing $ \delta  $ small enough, \eqref{InitialAssum}, \eqref{PrAssum}, \eqref{lm402}, \eqref{lm601} and \eqref{lm801} imply
	\begin{equation}
		\begin{aligned}
			& \biggl( \dfrac{(r^2r_x)_t}{r^2r_x} \biggr)^2 \leq C \alpha^2 ( ( \alpha^2 E_0 + (\alpha^2 +\Gamma)E_1 + \alpha^{2}\Gamma^{2\gamma-1})\times ( \alpha^3 E_0 + \alpha^5 E_1 + 1) \\& ~~~~~~  + E_2 ) + C \alpha^3 ( E_0 + E_1 ) + \Gamma^{2\gamma}.
		\end{aligned}
	\end{equation}
	Then \eqref{lm901} follows after noticing $ \alpha> 1 $ and the fact
	\begin{displaymath}
		\dfrac{(r^2r_x)_t}{r^2r_x} = \dfrac{v_x}{r_x} + 2\dfrac{v}{r}.
	\end{displaymath}
	
	To show \eqref{lm902}, from \eqref{uee0051}, \eqref{lm801}, \eqref{lm601} %, and \eqref{lm802},
	\begin{equation}\label{ueept002}
	\begin{aligned}
	&	\biggl| \dfrac{v_x}{r_x} - \dfrac{v}{r} \biggr|^2 \leq C \alpha^2  \int\biggl(\dfrac{x}{r}\biggr)^2 \rho_0 v_t^2\,dx + C \biggl\Arrowvert\biggl( \dfrac{x^2\rho_0}{r^2r_x} \biggr)^{2\gamma}\biggr\Arrowvert_{L^\infty_x}\\
	& \leq C \alpha^2 ( ( \alpha^2 E_0 + (\alpha^2 +\Gamma)E_1 + \alpha^{2}\Gamma^{2\gamma-1})\times ( \alpha^3 E_0 + \alpha^5 E_1 + 1 ) \\& ~~~~~~ + E_2 ) + C\Gamma^{2\gamma}.
	\end{aligned}
	\end{equation}
	
	In order to show \eqref{lm9021}, from \eqref{uee003}
	\begin{equation*}
		\begin{aligned}
		& -(2\mu + \lambda) \dfrac{d}{dt} \ln{r^2r_x} \leq -(2\mu + \lambda) \dfrac{(r^2r_x)_t}{r^2r_x} + \biggl(\dfrac{x^2 \rho_0}{r^2r_x}\biggr)^\gamma = \int_x^1 \dfrac{x^2}{r^2}\rho_0 v_t\,dx \\ & ~~ - 4\mu \dfrac{d}{dt}\ln{R(t)} = \dfrac{d}{dt} \biggl\lbrace \int_x^1 \dfrac{x^2}{r^2}\rho_0v\,dx -4\mu \ln {R(t)} \biggr\rbrace + 2 \int_x^1 \dfrac{x^2}{r^2} \dfrac{v}{r} \rho_0 v\,dx,
		\end{aligned}
	\end{equation*}
	where \begin{displaymath}
		R(t) = r(x = 1,t) .
	\end{displaymath}
	Integration in temporal variable yields, noticing $ R(0) = 1$ and $ r(x,t=0) = x $,
	\begin{equation}\label{uee010}
		\begin{aligned}
			& -(2\mu+\lambda) \ln { \dfrac{r^2r_x}{x^2}} \leq  - 4\mu \ln{R(t)} + \int_x^1 \dfrac{x^2}{r^2}\rho_0 v\,dx - \int_x^1 \rho_0 u_0\,dx \\
			& ~~~~ + 2 \int \int_x^1 \dfrac{x^2}{r^2}\dfrac{v}{r}\rho_0 v\,dx\,dt.
		\end{aligned}
	\end{equation}
%	where it has been applied the fact $ R(0) = 1$ and $ \frac{r^2 r_x}{x^2}(x,t=0) = 1 $.
	Denote
	\begin{displaymath}
		h = \int_x^1 \dfrac{x^2}{r^2}\rho_0 v\,dx - \int_x^1 \rho_0 u_0\,dx + 2 \int \int_x^1 \dfrac{x^2}{r^2}\dfrac{v}{r}\rho_0 v\,dx\,dt.
	\end{displaymath}
%	From \eqref{PrAssum}, \eqref{PrAssum2}, \eqref{InitialAssum}, \eqref{lm1}, and \eqref{lm503},
	Then, \eqref{uee010} can be written as
	\begin{equation}\label{uee011}
		\begin{aligned}
			& \dfrac{x^2}{r^2r_x} \leq R(t)^{-\frac{4\mu}{2\mu + \lambda}} \exp{\frac{h}{2\mu+\lambda} }. % \leq R(t)^{-\frac{4\mu}{2\mu+\lambda}} \exp \frac{C\alpha \beta + C \alpha(\alpha^\gamma + \alpha^{2} E_1 +  E_{0})}{2\mu + \lambda}
		\end{aligned}
	\end{equation}
	Denote $$\mathcal X = \max_{0\leq x\leq 1} \dfrac{x^2}{r^2r_x}. $$ Then 
	\begin{displaymath}
		R^3(t) = 3 \int_0^1 r^2r_x \,dx \geq 3 {\mathcal X}^{-1} \int_0^1 x^2 \,dx = {\mathcal X}^{-1}.
	\end{displaymath}
	Therefore, from \eqref{uee011}, it holds
	\begin{equation}
		\mathcal{X}^{\frac{2\mu+3\lambda}{3(2\mu+\lambda)}} \leq \exp (C \sup_{0\leq x\leq 1} h). 
	\end{equation}
	Meanwhile, by applying H\"older inequality, \eqref{InitialAssum}
	\begin{equation}\label{uee0101}
		\begin{aligned}
			& h\leq C  \biggl( \int \dfrac{x^2}{r^2} \,dx \biggr)^{1/2} \biggl( \int \dfrac{x^2}{r^2} v^2\,dx \biggr)^{1/2}  - \int_x^1 \rho_0u_0 \,dx \\
			& ~~~~~~ + C \biggl(\int \int \dfrac{x^2}{r^2} v^2 \,dx\,dt\biggr)^{1/2}\biggl(\int \int \dfrac{x^2}{r^2} \biggl(\dfrac{v}{r}\biggr)^2\,dx\,dt\biggr)^{1/2}. %\leq C \alpha \beta + C \int \alpha \biggl\lbrace \int r_xv^2 \,dx + \biggl\Arrowvert \dfrac{v}{r}\biggr\Arrowvert_{L_x^\infty}^2 \biggr\rbrace \,dt\\
		%	& \leq C\alpha \beta + C \alpha(\alpha^\gamma + \alpha^{2} E_1 +  E_{0}).
		\end{aligned}
	\end{equation}
	In addition, from \eqref{InitialAssum}, \eqref{PrAssum}, \eqref{lm1}, \eqref{lm2},
%	\begin{equation}
%		\begin{aligned}
%			& \int \dfrac{x^2}{r^2} v\,dx \leq C \biggl( \int \dfrac{x^2}{r^2} \,dx \biggr)^{1/2} \biggl( \int \dfrac{x^2}{r^2} v^2\,dx \biggr)^{1/2} \leq C \biggl( \int r_x v^2 \,dx \biggr)^{1/2}\\
%			& ~~~~ \leq C ( (1 + \beta )E_0 + E_1 )^{1/2},\\
%%%			& \int \int \dfrac{x^2}{r^2} v^2 \,dx\,dt \leq C  \int \int r_x v^2 \,dx\,dt \leq C E_0,
%		\end{aligned}
%	\end{equation}
%	where we have used Sobolev embedding in temporal variable,
	\begin{displaymath}
		\begin{aligned}
			& \int \dfrac{x^2}{r^2} v^2 \,dx \leq C \alpha^3 \int r_x v^2\,dx \leq  C \alpha^3 \int\,dt \int r_x v^2\,dx \\
			& ~~~~~~ + C \alpha^3\int \,dt \int \biggl(r_x v^2\biggr)_t\,dx \leq C \alpha^3 \biggl(1 + \biggl\Arrowvert \dfrac{v_x}{r_x}\biggr\Arrowvert_{L_x^\infty} \biggr) \int\,dt\int r_x v^2 \,dx \\
			& ~~~~~~ +C \alpha^3 \int\,dt\int r_x v_t^2 \,dx \leq C \alpha^3 ( (1 + \beta )E_0 + E_1 ) ,\\
			& \biggl| \int_x^1\rho_0 u_0\,dx \biggr| \leq C \biggl(\int u_0^2\,dx \biggr)^{1/2}.
		\end{aligned}
	\end{displaymath}
	Notice, since $ r_x(x,0) = 1 $, 
	\begin{displaymath}
		\begin{aligned}
			& \int u_0^2 \,dx \leq \sup_{t \geq 0} \int r_x v^2 \,dx \\
			& ~~~~~~ \leq \int\,dt \int r_x v^2 \, dx + \int \,dt \int (r_x v^2)_t \,dx \leq C ( (1+\beta ) E_0 + E_1).
		\end{aligned}
	\end{displaymath}
%	Moreover, from \eqref{lm303}, $ \exists \sigma > 0 $ such that $ R(t) > \sigma $. Hence,
%	\begin{equation}
%		\begin{aligned}
%			&\int \dfrac{x^2}{r^2} \biggl( \dfrac{v}{r}\biggr)^2 \,dx \leq C \int r_x \biggl( \dfrac{v}{r}\biggr)^2 \,dx = \int_0^{R(t)} \biggl( \dfrac{u}{r} \biggr)^2\,dr\\
%			& = \int_0^\sigma \dfrac{u^2}{r^2} \,dr + \int_\sigma^{R(t)} \dfrac{u^2}{r^2}\,dr \leq C \int_0^\sigma u_r^2 \,dr + \int_\sigma^{R(t)}		\end{aligned}
%	\end{equation}	
%	
	Also, \eqref{PrAssum}, \eqref{PrAssum2}, \eqref{lm503}, \eqref{lm1}
	\begin{displaymath}
		\begin{aligned}
			& \int \int \dfrac{x^2}{r^2} \biggl(\dfrac{v}{r}\biggr)^2\,dx\,dt \leq C \alpha^2 \int \biggl\Arrowvert \dfrac{v}{r} \biggr\Arrowvert_{L_x^\infty}^2\,dt\leq C \alpha^2 (\alpha^3 E_0 + \alpha^5 E_1 + 1) ,\\
			& \int \int \dfrac{x^2}{r^2} v^2 \,dx \,dt \leq C \alpha^3 \int \int r_x v^2 \,dx \leq C \alpha^3 E_0.
		\end{aligned}
	\end{displaymath}
	Therefore, from  \eqref{PrAssum2}, \eqref{uee0101} and the fact $ \alpha> 1 $, 
	\begin{equation}
		h \leq C \alpha^{5/2} ( (1 + \beta )E_0 + E_1 )^{1/2} + C \alpha^{5/2} ( \alpha^3 E_0 + \alpha^5 E_1 + 1)^{1/2} E_0^{1/2}.
	\end{equation}
	Consequently, % since $ \alpha > 1 $, 
	\begin{equation}
		\mathcal X \leq \exp \biggl\lbrace C\alpha^{5/2} ( (1 + \beta )E_0 + E_1 )^{1/2} + C \alpha^{5/2} ( \alpha^3 E_0 + \alpha^5 E_1 + 1)^{1/2} E_1^{1/2}\biggr\rbrace.
	\end{equation}
\end{pf}
	
With Lemma \ref{lm:pointwisebound}, it is possible to design the conditions on $ \mathcal E_0, \mathcal E_1, \mathcal E_2 $ such that the a priori assumption \eqref{PrAssum} can be verified for a smooth solution. 
	
\begin{lm}\label{lm:extendinglemma}
For a fixed $ \alpha $, there is a constant $ \epsilon_0 > 0 $, depending on $ \alpha, \bar\rho, \mu, \lambda, M $, such that for $ \mathcal E_0, \mathcal E_1, \mathcal E_2 \leq \epsilon_0 $, $ \exists  \beta_0 =  \beta_0 (\alpha, \mathcal E_0, \mathcal E_1, \mathcal E_2,M) $, $ \beta_1 =  \beta_1 (\alpha, \mathcal E_0, \mathcal E_1, \mathcal E_2,M) $, $ \beta_1 \geq \beta_0 > M $, satisfying the following. If $ \beta_0 \leq \beta \leq  \beta_1 $ in \eqref{PrAssum}, then the following inequality holds,
	\begin{equation}\label{lm903}
		\biggl\Arrowvert \dfrac{v_x}{r_x} \biggr\Arrowvert_{L_x^\infty}^2, \biggl\Arrowvert \dfrac{v}{r}\biggr\Arrowvert_{L^{\infty}_x}^2 \leq \dfrac{1}{2} \beta^2 < \beta^2.
	\end{equation}
	In addition, there exists a constant $ \epsilon_1 \leq \epsilon_0 $, depending on $ \alpha, \bar\rho, \mu, \lambda, M $, such that for $ \mathcal E_0, \mathcal E_1 \leq \epsilon_1, \mathcal E_2 \leq \epsilon_0 $,  if $\beta\leq \beta_1 (\alpha,\mathcal E_0, \mathcal E_1, \mathcal E_2, M) $, it holds, 
	\begin{equation}\label{lm904}
		\biggl\Arrowvert\dfrac{x^2}{r^2r_x}\biggr\Arrowvert_{L^{\infty}_x} \leq \alpha^{3/2} < \alpha^3.
	\end{equation}
	Moreover, $ \epsilon_0, \epsilon_1 $ can be chosen such that,
	\begin{equation}\label{lm905}
		\epsilon_0 = \epsilon_1 = \epsilon(\alpha),
	\end{equation}
	with $ \epsilon(\alpha) $ being a bounded continuous function of $ \alpha \in (1,  +\infty) $. We denote the following quantity which is the constant in \eqref{InitialAssum4},
	\begin{equation}\label{lm906}
		\bar\epsilon := \sup_{1 < \alpha < +\infty} \epsilon(\alpha) < \infty.
	\end{equation}
\end{lm}
	
\begin{pf}
	From \eqref{lm901}, \eqref{lm902}, there exists $ C_1 > 0 $ such that
	\begin{equation}\label{uee0091}
		\begin{aligned}
			& \biggl\Arrowvert \dfrac{v_x}{r_x} \biggr\Arrowvert_{L_x^\infty}^2, \biggl\Arrowvert \dfrac{v}{r}\biggr\Arrowvert_{L^{\infty}_x}^2 \leq C_1 \alpha^2 ( ( \alpha^2 E_0 + (\alpha^2 +\Gamma)E_1 + \alpha^{2}\Gamma^{2\gamma-1}) \\
			& ~~~~~~~\times ( \alpha^3 E_0 + \alpha^5 E_1 + 1) + E_2 ) + C_1 \Gamma^{2\gamma}.
		\end{aligned}
	\end{equation}
	From definition of $ E_0, E_1, E_2 $ in \eqref{TEnergy0} and \eqref{lm602}, \eqref{uee0091} can be written as,
	\begin{equation}
		\biggl\Arrowvert \dfrac{v_x}{r_x} \biggr\Arrowvert_{L_x^\infty}^2, \biggl\Arrowvert \dfrac{v}{r}\biggr\Arrowvert_{L^{\infty}_x}^2 \leq V_{\alpha}(\beta, \mathcal E_0, \mathcal E_1, \mathcal E_2),
	\end{equation}
	for some positive, increasing, continuous function $ V_\alpha $ satisfying,
	\begin{displaymath}
		V_\alpha(0,0,0,0) = C_1 \alpha^4 \Gamma_0^{2\gamma-1}  + C_1\Gamma_0^{2\gamma}, ~~ \Gamma_0 := \min  \lbrace \bar\rho \alpha^3, \bar\rho \rbrace.
	\end{displaymath}
	In addition, there exist $ C > 0 $ and positive integers $ l, n > 2 $ such that, % $ V_\alpha $ satisfies
	\begin{equation}\label{uee0092}
		  V_\alpha \leq C \alpha^l + A_{1} \beta + A_{2} \beta^n,
	\end{equation}
	for some non-negative, increasing, continuous functions $ A_{1} =  A_{1}(\alpha,\mathcal E_0, \mathcal E_1, \mathcal E_2)$, $A_{2} = A_{2} (\alpha,\mathcal E_0, \mathcal E_1, \mathcal E_2) $ satisfying
	\begin{displaymath} \begin{aligned}
		& A_1 (\alpha, 0,0,0)=0, ~ A_2(\alpha, 0,0,0) = 0.
	\end{aligned}
	\end{displaymath}
	Moreover, $ A_1, A_2 > 0 $ for $ (\mathcal E_0, \mathcal E_1, \mathcal E_2 ) \neq (0,0,0) $.
	Therefore, for some $ \epsilon_0  > 0 $ depending on $ \alpha, \bar\rho, \mu, \lambda, M $, it holds for $ 0 < \mathcal E_0, \mathcal E_1, \mathcal E_2 \leq \epsilon_0 $, the set % there exists a constant $ \beta_0 > 0 $ depending on $ \alpha, \bar\rho, \mu, \lambda, \eta, \mathcal E_0, \mathcal E_1, \mathcal E_2 $, such that 
	\begin{equation}
		S_\alpha : = \biggl\lbrace \beta \geq M | V_\alpha(\beta, \mathcal E_0, \mathcal E_1, \mathcal E_2) \leq \dfrac {1}{2} \beta^2 \biggr\rbrace% \subset \biggl\lbrace \beta > 0: \beta \leq  \beta_0 \biggr\rbrace.
	\end{equation}
	is nonempty and bounded. %Since $ V_\alpha $ is increasing in $ \alpha $, it is easy to show
%	\begin{equation}\label{uee0092}
%		S_{\alpha'} \subset S_\alpha, ~~ \text{for} ~ \alpha' \geq \alpha.
%	\end{equation}
	Define
	\begin{equation}
		\beta_0 = \beta_0(\alpha, \mathcal E_0, \mathcal E_1, \mathcal E_2,M) := \inf S_\alpha, ~ \beta_1 =\beta_1(\alpha, \mathcal E_0, \mathcal E_1, \mathcal E_2,M) := \sup S_\alpha . %\biggl\lbrace \beta > 0 | V_\alpha(\beta, \mathcal E_0, \mathcal E_1, \mathcal E_2) \leq \beta^2 \biggr\rbrace.
	\end{equation}
	Using this notation, $ S_\alpha = [\beta_0, \beta_1] $. 
	Then if $ \mathcal E_0, \mathcal E_1, \mathcal E_2 \leq \epsilon_0 $, for $ \beta \in S_\alpha $,  
	\begin{displaymath}
		\biggl\Arrowvert \dfrac{v_x}{r_x} \biggr\Arrowvert_{L_x^\infty}^2, \biggl\Arrowvert \dfrac{v}{r}\biggr\Arrowvert_{L^{\infty}_x}^2 \leq V_{\alpha}(\beta, \mathcal E_0, \mathcal E_1, \mathcal E_2) \leq \dfrac{1}{2} \beta^2  < \beta^2.
	\end{displaymath}
% ==========
	To show \eqref{lm904}, from \eqref{lm9021} we have,
	\begin{equation}
		\left\Arrowvert \dfrac{x^2}{r^2r_x} \right\Arrowvert_{L_x^\infty} \leq \exp \left\lbrace W_\alpha(\beta, \mathcal E_0, \mathcal E_1)\right\rbrace.
	\end{equation}
	$ W_\alpha $ is a non-negative, increasing, continuous function such that
	\begin{equation}\label{uee0093}
		W_\alpha \leq A_{3}\alpha^l,
	\end{equation}
	for some $ l > 0 $, where $ A_{3} = A_{3}(\beta,\mathcal E_0, \mathcal E_1) $ is a non-negative, increasing, continuous function and
	\begin{displaymath}
		A_{3}(\beta, 0, 0 ) = 0.
	\end{displaymath}
	Therefore, there exists $ \epsilon_1 \leq \epsilon_0 $ depending on $ \alpha $ such that for $ \mathcal E_0, \mathcal E_1 \leq \epsilon_1, \mathcal E_2 \leq \epsilon_0, \beta\leq \beta_1 $, 
	\begin{displaymath}
		W_\alpha(\beta,\mathcal E_0, \mathcal E_1 ) \leq W_\alpha(\beta_{1}, \mathcal E_0, \mathcal E_1 ) \leq \dfrac{3}{2} \ln {\alpha},
	\end{displaymath}
	and therefore, since $ \alpha > 1 $,
	\begin{equation}
		\left\Arrowvert \dfrac{x^2}{r^2r_x} \right\Arrowvert_{L_x^\infty} \leq \alpha^{3/2} < \alpha^3.
	\end{equation}
	\eqref{lm905}, \eqref{lm906} follow from \eqref{uee0092}, \eqref{uee0093}  and the continuity of $V_\alpha, W_\alpha$ in $ \alpha \in (1, +\infty)$ and similar arguments. 

\end{pf}

\subsection{Regularity}\label{sec:regularity}

In this section, the goal is to study the regularity propagated by the strong solution to \eqref{LgNS} with the a prior bound \eqref{PrAssum}. %It is further assumed 
%\begin{equation}
%	\left\Arrowvert (\rho_0^\gamma)_x \right\Arrowvert_{L_x^2} < \infty.
%\end{equation}

\begin{lm} For any $ T >0 $, there are constants $ C > 0 $ depending on $ \alpha, \beta$, and $ C_T > 0 $ depending on $ T $ such that for $ 0 \leq t \leq T $, 
	\begin{gather}
	C \leq \dfrac{r}{x} \leq C e^{Ct}, ~~~ C e^{- C t} \leq r_x \leq C e^{Ct},\label{lm1001}	\\
	\int_0^T \, dt \int \left( x^2 v_x^2 + x^2 v_{xt}^2 + v^2 + v_t^2 + v_{xt}^2 + \dfrac{v_t^2}{x^2} \right) \,dx \leq C_T, \label{lm1002}\\
	\int \left( v_x^2 + \dfrac{v^2}{x^2} + \rho_0 v_t^2 \right) \,dx \leq C_T, \label{lm1003}
	\end{gather}
and therefore,
\begin{equation}\label{lm1004}
	\int v^2 \,dx \leq C_T.
\end{equation}
\end{lm}

\begin{pf}
%	From \eqref{lm903}, 
%	\begin{equation*}
%		c\leq \dfrac{\dt r}{r} \leq C ,~~ c \leq \dfrac{\dt r_x}{r_x} \leq C,
%	\end{equation*}
%	integration in temporal variable yileds
%	\begin{displaymath}
%		c t \leq \ln{\dfrac{r}{x}} \leq C t, ~~ ct \leq \ln{r_x} \leq C t.
%	\end{displaymath}
%	with, together with \eqref{lm904} shows \eqref{lm1001}.

\eqref{lm1001} is a direct consequence of \eqref{PrAssum}, \eqref{PrAssum2}.
\eqref{lm1002} and \eqref{lm1003} follow from \eqref{lm1}, \eqref{lm2}, \eqref{lm6}, \eqref{lm801}. Embedding theory yields \eqref{lm1004}.
\end{pf}

Next lemma is concerning the regularity at the centre inspired by \cite{LuoXinZeng2016}.
\begin{lm}\label{lm:D2regu} There is a constant $ C_T>0 $ % depending on $ \mu, \lambda, \eta, \bar\rho, \alpha, \beta, T $
 such that,
	\begin{align}
	& \int r_{xx}^2 \,dx + \int\left\lbrack\left(\dfrac{r}{x}\right)_x\right\rbrack^2\,dx \leq C_T,      \label{lm1102} \\
	& \int v_{xx}^2\,dx + \int \left\lbrack\left(\dfrac{v}{x}\right)_x\right\rbrack^2\,dx \leq C_T.   \label{lm1101}
	\end{align}

\end{lm}

\begin{pf}
Here, we use the structure from \cite{LuoXinZeng2016} to obtain high-order estimates. 
	Define \begin{equation}
		\mathcal{G} = \ln \left( \dfrac{r^2 r_x}{x^2} \right).
	\end{equation}
	\eqref{LgNS} can be written as
	\begin{equation}\label{rg1001}
		(2\mu  + \lambda ) \mathcal{G}_{xt} + \gamma \left( \dfrac{x^2 \rho_0}{r^2r_x} \right)^\gamma \mathcal{G}_x = \left(\dfrac{x}{r}\right)^2 \rho_0 v_t + \left(\dfrac{x^2}{r^2r_x}\right)^\gamma (\rho_0^\gamma)_x.
	\end{equation}
	Multiple \eqref{rg1001} with $ \mathcal{G}_x $ and integrate the resulting equation in the spatial variable,
	\begin{equation*}
	\begin{aligned}
		& \dfrac{d}{dt} \left\lbrace \dfrac{2\mu+\lambda}{2} \int \mathcal{G}_x^2 \,dx \right\rbrace + \gamma \int \left( \dfrac{x^2 \rho_0}{r^2r_x}\right)^\gamma \mathcal{G}_x^2\,dx = \int \left(\dfrac{x}{r}\right)^2\rho_0v_t \mathcal{G}_x \,dx \\
		& ~~~~  + \int \left(\dfrac{x^2}{r^2r_x}\right)^\gamma\left(\rho_0^\gamma\right)_x \mathcal{G}_x\,dx \leq \int \mathcal{G}_x^2 \,dx + C \int v_t^2 \,dx  + C \int (\left( \rho_0^\gamma \right)_x )^2 \,dx. \\
	\end{aligned}
	\end{equation*}
	From \eqref{lm1002}, Gr\"onwall's inequality then yields
	\begin{equation}\label{rg1002}
		\int \mathcal{G}_x^2\,dx \leq C_T,
	\end{equation} 
	and  as a consequence of \eqref{lm1003}, \eqref{rg1001}, \eqref{rg1002},
	\begin{equation}\label{rg1003}
		\int \mathcal{G}_{xt}^2 \,dx \leq C_T.
	\end{equation}
	Notice
	\begin{displaymath}
		\begin{aligned}
		& \mathcal{G}_x = \dfrac{x}{rr_x} \left( 2 r_x \left( \dfrac{r}{x} \right)_x + r_{xx} \left( \dfrac{r}{x}\right) \right), \\
		& \mathcal{G}_{xt} = \dfrac{x}{rr_x} \left( 2 r_x \left( \dfrac{v}{x} \right)_x + v_{xx} \left( \dfrac{r}{x}\right) + 2 v_x \left(\dfrac{r}{x}\right)_x + r_{xx}\left(\dfrac{v}{x}\right) \right) \\
		& ~~~~~~ - \dfrac{x}{rr_x} \left(\dfrac{v_x}{r_x} + \dfrac{v}{r}  \right) \left( 2 r_x \left( \dfrac{r}{x} \right)_x + r_{xx} \left( \dfrac{r}{x}\right) \right).
		\end{aligned}
	\end{displaymath}
	Therefore \eqref{lm1001}, \eqref{rg1002} imply
	\begin{equation}\label{rg10031}
		\begin{aligned}
			& C_T \geq \int \left( 2 r_x \left( \dfrac{r}{x} \right)_x + r_{xx} \left( \dfrac{r}{x}\right) \right)^2 \,dx = \underbrace{\int r_{xx}^2 \left( \dfrac{r}{x} \right)^2\,dx }\\
			& ~~~~~~ + \underbrace{\int 4 r_x r_{xx} \left( \dfrac{r}{x} \right) \left( \dfrac{r}{x}\right)_x + 4 r_x^2 \left\lbrack\left( \dfrac{r}{x} \right)_x\right\rbrack^2 \,dx} : = A+B.
		\end{aligned}
	\end{equation}
	At the same time, we shall use the following identity to manipulate $ B $,
	\begin{displaymath}
		r_x = x \left( \dfrac{r}{x} \right)_x + \dfrac{r}{x} , ~ r_{xx} = x \left( \dfrac{r}{x}\right)_{xx} + 2 \left( \dfrac{r}{x}\right)_x.
	\end{displaymath}
	Consequently, it follows from a complicated and direct calculation and integration by parts,  
	\begin{equation}
		\begin{aligned}
			& B = 4 \int \biggl\lbrace \left( x \left(\dfrac{r}{x}\right)_x + \dfrac{r}{x}\right) \left( x \left( \dfrac{r}{x} \right)_{xx} + 2 \left( \dfrac{r}{x} \right)_x \right) \left( \dfrac{r}{x}\right)\left(\dfrac{r}{x}\right)_x  \\ 
			& ~~~~~~ + \left( x\left(\dfrac{r}{x}\right)_x + \dfrac{r}{x} \right)^2 \left\lbrack\left(\dfrac{r}{x}\right)_x\right\rbrack^2 \biggr\rbrace  \,dx \\
%			& = 4\int x \left(\dfrac{r}{x}\right)^2 \left(\dfrac{r}{x}\right)_x\left(\dfrac{r}{x}\right)_{xx} \,dx + 4 \int x^2 \left( \dfrac{r}{x}\right) \left(\dfrac{r}{x}\right)^2_x \left(\dfrac{r}{x}\right)_{xx} \,dx \\
%			& ~~~~~~ + 12 \int \left(\dfrac{r}{x}\right)^2 \left(\dfrac{r}{x}\right)_{x}^2 \,dx + 16 \int x \left(\dfrac{r}{x}\right) \left(\dfrac{r}{x}\right)_x^3\,dx + 4 \int x^2 \left(\dfrac{r}{x}\right)^4_x\,dx\\
			& = \left. \left( 2 x \left(\dfrac{r}{x}\right)^2 \left\lbrack\left(\dfrac{r}{x}\right)_x\right\rbrack^2 + \dfrac{4}{3}x^2\left(\dfrac{r}{x}\right)\left\lbrack\left(\dfrac{r}{x}\right)_x\right\rbrack^3 \right) \right|_{x=1} + 10 \int \left(\dfrac{r}{x}\right)^2 \left\lbrack\left(\dfrac{r}{x}\right)_x\right\rbrack^2 \,dx \\
			& ~~~~~~ + \dfrac{28}{3} \int x \left(\dfrac{r}{x}\right) \left\lbrack\left(\dfrac{r}{x}\right)_x\right\rbrack^3\,dx + \dfrac{8}{3} \int x^2 \left\lbrack\left(\dfrac{r}{x}\right)_x\right\rbrack^4\,dx\\
			& = \left. \left( 2 x \left(\dfrac{r}{x}\right)^2 \left\lbrack\left(\dfrac{r}{x}\right)_x\right\rbrack^2 + \dfrac{4}{3}x^2\left(\dfrac{r}{x}\right)\left\lbrack\left(\dfrac{r}{x}\right)_x\right\rbrack^3 \right) \right|_{x=1}\\
			& ~~~~~~ + \dfrac{2}{3} \int \left( 2 x \left\lbrack\left(\dfrac{r}{x}\right)_x \right\rbrack^2 + \dfrac{7}{2} \dfrac{r}{x}\left(\dfrac{r}{x}\right)_x \right)^2 \,dx + \dfrac{11}{6} \int \left(\dfrac{r}{x}\right)^2 \left\lbrack\left(\dfrac{r}{x}\right)_x\right\rbrack^2\,dx. \\
			%& \simeq \left. \left( C \dfrac{r^4}{x^4} + C \dfrac{r^3}{x^3}r_x + C \dfrac{r^2}{x^2}r_x^2 + C \dfrac{r}{x}r_x^3 \right)\right|_{x=1} 
%			& \geq -C_T + C \int x^2 \left\lbrack\left(\dfrac{r}{x}\right)_x\right\rbrack^4\,dx + C \int \left(\dfrac{r}{x}\right)^2 \left\lbrack\left(\dfrac{r}{x}\right)_x\right\rbrack^2\,dx
		\end{aligned}
	\end{equation}
	Meanwhile, from \eqref{lm1001}
	\begin{displaymath}
		\left| \left. \left(\dfrac{r}{x}\right)_x \right|_{x=1}\right| = \left|\left. ( r_x - r )\right|_{x=1} \right| \leq C_T.
	\end{displaymath}
	Consequently, from \eqref{rg10031}, it holds,
	\begin{equation}
		\int \left(\dfrac{r}{x}\right)^2 r_{xx}^2 \,dx + \int \left(\dfrac{r}{x}\right)^2  \left\lbrack\left(\dfrac{r}{x}\right)_x\right\rbrack^2\,dx+ \int x^2 \left\lbrack\left(\dfrac{r}{x}\right)_x\right\rbrack^4\,dx \leq A+B+C_T \leq C_T,
	\end{equation}
	and together with \eqref{lm1001}, \eqref{lm1102} follows.
	Similarly, from \eqref{PrAssum}, \eqref{lm1001}, \eqref{lm1102}, \eqref{rg1003}, 
	\begin{equation}
		C_T \geq C_T \int \mathcal{G}_{xt}^2 \,dx \geq \int \left(2r_x \left(\dfrac{v}{x} \right)_x + v_{xx} \left(\dfrac{r}{x}\right)\right)^2\,dx - C_T,
	\end{equation}
	and hence,
	\begin{equation}\label{rg10032}
		\begin{aligned}
			& C_T \geq \int \left(2r_x \left(\dfrac{v}{x} \right)_x + v_{xx} \left(\dfrac{r}{x}\right)\right)^2\,dx = \underbrace{ 4 \int r_x^2 \left\lbrack\left(\dfrac{v}{x}\right)_x\right\rbrack^2 \,dx + \int \left(\dfrac{r}{x}\right)^2 v_{xx}^2 \,dx}\\
			& ~~~~~~ + \underbrace{4 \int r_x \left(\dfrac{r}{x}\right) \left( \dfrac{v}{x}\right)_x v_{xx}\,dx} := \tilde A + \tilde B.
		\end{aligned}
	\end{equation}
	Making use of the identity,
	$$ v_{xx} = x \left( \dfrac v x \right)_{xx} + 2 \left( \dfrac v x \right)_x, $$
	integration by parts yields
	\begin{equation}
		\begin{aligned}
			& \tilde B = 4 \int r_x \left(\dfrac{r}{x}\right) \left( \dfrac{v}{x}\right)_x \left( x \left(\dfrac{v}{x}\right)_{xx}+ 2 \left(\dfrac{v}{x}\right)_x \right)\,dx \\
			& ~~ = \left.\left( 2 x r_x \left(\dfrac{r}{x}\right)\left\lbrack\left(\dfrac{v}{x}\right)_x\right\rbrack^2 \right)\right|_{x=1} + 6 \int r_x \left(\dfrac{r}{x}\right) \left\lbrack\left(\dfrac{v}{x}\right)_x\right\rbrack^2\,dx\\
			& ~~~~~~ - 2 \int x r_{xx}\left(\dfrac{r}{x}\right) \left\lbrack\left(\dfrac{v}{x}\right)_x\right\rbrack^2\,dx - 2 \int x r_{x}\left(\dfrac{r}{x}\right)_x \left\lbrack\left(\dfrac{v}{x}\right)_x\right\rbrack^2\,dx,\\
%			& =  \left.\left( 2 x r_x \left(\dfrac{r}{x}\right)\left\lbrack\left(\dfrac{v}{x}\right)_x\right\rbrack^2 \right)\right|_{x=1} + 6 \int \left(\dfrac{r}{x}\right)^2 \left\lbrack \left(\dfrac{v}{x}\right)_x\right\rbrack^2\,dx\\
%			& ~~~~~~ - 2 \int x^2 \left(\dfrac{r}{x}\right)_{xx} \left(\dfrac{r}{x}\right) \left\lbrack\left(\dfrac{v}{x}\right)_x\right\rbrack^2\,dx - 2\int x^2 \left\lbrack\left(\dfrac{r}{x}\right)_x\right\rbrack^2 \left\lbrack\left(\dfrac{v}{x}\right)_x\right\rbrack^2\,dx  
		\end{aligned}
		\end{equation}
%		where the last equality follows the fact
%		\begin{displaymath}
%			\begin{aligned}
%				& r_{xx} = x \left(\dfrac{r}{x}\right)_{xx} + 2 \left(\dfrac{r}{x}\right)_x, r_x = x \left(\dfrac{r}{x}\right)_x + \dfrac{r}{x}.
%			\end{aligned}
%		\end{displaymath}
%		Meanwhile, 
%		\begin{displaymath}
%			\begin{aligned}
%				& \int x^2 \left(\dfrac{r}{x}\right)_{xx} \left(\dfrac{r}{x}\right) \left\lbrack\left(\dfrac{v}{x}\right)_x\right\rbrack^2\,dx \leq 3 \int \left(\dfrac{r}{x}\right)^2 \left\lbrack\left(\dfrac{v}{x}\right)_x \right\rbrack^2\,dx\\
%				& ~~~~~~ + C \int x^2 \left\lbrack\left(\dfrac{r}{x}\right)_{xx}\right\rbrack^2 \left| x \left(\dfrac{v}{x}\right)_x\right|^2\,dx \leq  3 \int \left(\dfrac{r}{x}\right)^2 \left\lbrack\left(\dfrac{v}{x}\right)_x \right\rbrack^2\,dx\\
%				& ~~~~~~ + C \int x^2 \left( \dfrac{r_{xx}}{x} - 2 \dfrac{r_x}{x^2} + 2 \dfrac{r}{x^3}\right)^2\,dx \leq 3 \int \left(\dfrac{r}{x}\right)^2 \left\lbrack\left(\dfrac{v}{x}\right)_x \right\rbrack^2\,dx\\
%				& ~~~~~~ + C \int \left(r_{xx}^2 + r_{x}^2 + r^2 \right) \,dx 
%			\end{aligned}
%		\end{displaymath}
		where
		\begin{displaymath}
%			\begin{aligned}
				\left| \left. \left( \dfrac{v}{x}\right)_x \right|_{x=1} \right| =  \left| \left. ( v_x - v )   \right|_{x=1} \right| \leq C_T.
%			\end{aligned}
		\end{displaymath}
		Therefore, from \eqref{rg10032}, we shall have
		\begin{equation} \label{rg1004}
%			\begin{split}
			\begin{aligned}
				& \int \left(\dfrac{r}{x}\right)^2 v_{xx}^2 \,dx + \int \left(4 r_x^2 + 6 r_x \left(\dfrac{r}{x}\right) - 2 x r_x \left(\dfrac{r}{x}\right)_x \right) \left\lbrack \left(\dfrac{v}{x}\right)_x \right\rbrack^2\,dx\\
				& ~~~~ \leq C_T + 2 \int x r_{xx}\left(\dfrac{r}{x}\right) \left\lbrack\left(\dfrac{v}{x}\right)_x\right\rbrack^2\,dx\\
				& ~~~~ \leq C_T + \int \left(\dfrac r x \right)^2 \left\lbrack\left(\dfrac v x \right)_x\right\rbrack^2\,dx + C \int r_{xx}^2 \left|x \left(\dfrac{v}{x}\right)_x\right|^2\,dx.
			\end{aligned}
%			\end{split}
		\end{equation}
		Then, since $ r_x = \frac{r}{x} + x \left(\frac{r}{x}\right)_x $,
		\begin{equation}\label{rg1101}
			\begin{aligned}
				& 4r_x^2 + 6r_x \left(\dfrac{r}{x}\right) - 2 x r_x \left( \dfrac r x \right)_x = 2 x^2 \left(\dfrac r x\right)_x^2 + 12 x \dfrac r x \left( \dfrac r x \right)_x + 10 \left(\dfrac r x\right)^2 \\
				& ~~~~~~ \geq 9\left( \dfrac{r}{x}\right)^2 - Cx^2 \left( \dfrac r x\right)_x^2.
			\end{aligned}
		\end{equation} Consequently, from \eqref{rg1004}, it follows
		\begin{displaymath}
			\begin{aligned}
				& \int \left(\dfrac{r}{x}\right)^2 v_{xx}^2 \,dx + 8 \int \left( \dfrac r x \right)^2 \left\lbrack \left(\dfrac v x \right)_x \right\rbrack^2 \,dx \leq C_T + C_T \int r_{xx}^2 \,dx \\
				& ~~~~~~ + C \int x^2 \left( \dfrac r x \right)_x^2 \left\lbrack \left(\dfrac v x \right)_x \right\rbrack^2\,dx \leq C_T + C_T \left\lbrace \int r_{xx}^2 \,dx + \int \left(\dfrac r x \right)_x^2 \,dx \right\rbrace,
			\end{aligned}
		\end{displaymath}
		where it has been applied \eqref{PrAssum}, \eqref{lm1001} and
		\begin{displaymath}
		\begin{aligned}
			& \left|x \left(\dfrac{v}{x}\right)_x\right| = \left|v_x - \dfrac{v}{x}\right| \leq C_T. \\
			\end{aligned}
		\end{displaymath}
		Therefore 
		\begin{equation}
			\int v_{xx}^2\,dx + \int \left\lbrack\left(\dfrac{v}{x}\right)_x\right\rbrack^2\,dx \leq C_T.
		\end{equation}
\end{pf}

%\begin{rmk}
%	Here we put a note about the vanishing of the boundary values at $ x=0 $, which is used in the calculation of $ B , \tilde B $ and in the following. We shall use the case when calculating $ B $ to illustrate the idea.  
%	Indeed, we are looking for the solution with $ \int \left\lbrack\left(\frac{r}{x} \right)_{x}\right\rbrack ^2 \,dx $ bounded. It suffices to show $$ \lim_{x\rightarrow 0^+} \left. x \left\lbrack\left(\dfrac{r}{x}\right)_x\right\rbrack^2 \right. = 0. $$ 
%	By the mean value theorem, $ \exists x/2 \leq \tilde x \leq x $ such that 
%	$$ \left. \left\lbrack\left(\frac{r}{x}\right)_x\right\rbrack^2 \right|_{x=\tilde x} = \dfrac{2}{x} \int_{x/2}^x \left\lbrack\left(\frac{r}{x}\right)_x\right\rbrack^2,$$
%	where $ 0 < x < 1 $.
%	Therefore,
%	$$  \left. x  \left\lbrack\left(\dfrac{r}{x}\right)_x\right\rbrack^2 \right|_{x=\tilde x} = \tilde x \cdot \dfrac{2}{x} \int_{x/2}^x \left\lbrack\left(\frac{r}{x}\right)_x\right\rbrack^2 \leq 2 \int_{x/2}^x \left\lbrack\left(\frac{r}{x}\right)_x\right\rbrack^2 \rightarrow 0 ~~~ \text{as} ~ x \rightarrow 0. $$
%	This finishes the proof in the continuous case. In general, it would follow by a density argument. 	
%\end{rmk}

% ============

\subsection{Higher Regularity}\label{sec:highregularity}

In this section, we discuss the regularity propagated by the classical solution to \eqref{LgNS}. It is assumed \eqref{PrAssum} and 
\begin{equation}
	\left\Arrowvert (\rho_0)_x \right\Arrowvert_{L_x^2} < \infty, ~\left\Arrowvert (\rho_0^\gamma)_{xx}\right\Arrowvert_{L_x^2} < \infty.
\end{equation}
%Also, denote
%\begin{align}
%	& \mathcal E_3 = \dfrac{1}{2}\int x^2 \rho_0 u_{2}^2\,dx,\\
%	& \mathcal E_4 = \dfrac{1}{2}\int \rho_0 u_{2}^2\,dx.
%\end{align}
%where 
%\begin{equation}\label{intut2}\begin{aligned}
%	& u_{2} = \dfrac{1}{\rho_0}\left\lbrace (2\mu + \lambda) \left\lbrack\dfrac{(x^2u_1)_x}{x^2} - u_{0,x}^2 - 2 \dfrac{u_0^2}{x^2} \right\rbrack_x + \gamma \left\lbrack\rho_0^\gamma \dfrac{(x^2 u_0)_x}{x^2} \right\rbrack_x  \right\rbrace  + 2 \dfrac{u u_1}{x}
%	\end{aligned}
%\end{equation}
%where $ u_1 $ is defined in \eqref{intut1}.\\
After taking the time derivative of \eqref{LgNS1}, one has
\begin{equation}\label{LgNS2}
	\begin{aligned}
		& x^2 \rho_0 v_{ttt} + r^2 \left\lbrack \left(\dfrac{x^2\rho_0}{r^2r_x}\right)^\gamma\right\rbrack_{xtt} = r^2 \mathfrak{B}_{xtt} +4\mu r^2 \left(\dfrac{v}{r}\right)_{xtt} \\
		& ~~~~~~ + 4 r v \left( \mathfrak{B}_{xt} + 4\mu \left(\dfrac{v}{r}\right)_{xt} \right) - 4 r v \left\lbrack \left(\dfrac{x^2\rho_0}{r^2r_x}\right)^\gamma \right\rbrack_{xt}\\
		& ~~~~~~ + 2 (rv)_t \left(\mathfrak{B}_{x} + 4\mu \left(\dfrac{v}{r}\right)_{x}\right) - 2 \left(rv\right)_t \left\lbrack\left(\dfrac{x^2\rho_0}{r^2r_x}\right)^\gamma\right\rbrack_x,
	\end{aligned}
\end{equation}
with the additional boundary conditions
\begin{equation}
	\left.\left\lbrack\left(\dfrac{x^2\rho_0}{r^2r_x}\right)^\gamma\right\rbrack_{tt} - \mathfrak{B}_{tt}\right|_{x=1} = 0,~ v_{tt}(0,t) = 0.
\end{equation}
In the following, unless it is stated otherwise, $ C > 0 $ is a generic constant depending on $ \mu, \lambda, \alpha, \beta, \bar\rho, \gamma $.
\begin{lm} There is a constant $ C > 0 $ % depending on $ \mu,\lambda, \bar\rho,\alpha,\beta $ 
such that,
	\begin{equation}\label{lm12}
	\begin{aligned}
		& \dfrac{1}{2}\int x^2 \rho_0 v_{tt}^2\,dx + \int \int r^2 r_x \left( \dfrac{v_{xtt}^2}{r_x^2} + 2 \dfrac{v_{tt}^2}{r^2} \right) \,dx  \leq C(\mathcal E_0 + \mathcal E_1 + \mathcal E_3).
	\end{aligned}
	\end{equation}
\end{lm}

\begin{pf}

Multiply \eqref{LgNS2} with $ v_{tt} $ and integrate the resulting in spatial variable.
%\begin{equation}
%	\begin{aligned}
%	&	\dfrac{d}{dt}\dfrac{1}{2} \int x^2 \rho_0 v_{tt}^2\,dx - \int (r^2 v_{tt})_x \left\lbrack\left(\dfrac{x^2\rho_0}{r^2r_x}\right)^\gamma\right\rbrack_{tt} \,dx \\
%	& ~~~~ = - \int (r^2v_{tt})_x\mathfrak{B}_{tt} + 4\mu r^2 \left(\dfrac{v}{r}\right)_{xtt}v_{tt}\,dx \\
%	& ~~~~~~ - \int 4 (rvv_{tt})_x\left( \mathfrak{B}_t - \left\lbrack\left(\dfrac{x^2\rho_0}{r^2r_x}\right)^\gamma\right\rbrack_{t}\right)\,dx + 16 \mu \int r v \left(\dfrac{v}{r}\right)_{xt}v_{tt}\,dx \\
%	& ~~~~~~ -\int 2 ((rv)_tv_{tt})_x \left( \mathfrak{B} - \left(\dfrac{x^2\rho_0}{r^2r_x}\right)^\gamma\right)\,dx + 8 \mu \int (rv)_t\left(\dfrac{v}{r}\right)_x v_{tt}\,dx
%	\end{aligned}
%\end{equation}
Similar as before, integration by parts then yields
\begin{equation}\label{rg1201}
	\begin{aligned}
		& \dfrac{d}{dt}\dfrac{1}{2}\int x^2 \rho_0 v_{tt}^2\,dx + 2 \mu \int r^2 r_x \left( \dfrac{v_{xtt}^2}{r_x^2} + 2 \dfrac{v_{tt}^2}{r^2} \right) \,dx \\
		& ~~~~~~ + \lambda \int r^2 r_x \left( \dfrac{v_{xtt}}{r_x} + 2 \dfrac{v_{tt}}{r}\right)^2\,dx = L_1+ L_2+ L_3,
	\end{aligned}
\end{equation}
with
%\begin{displaymath}
	\begin{align*}
		& L_1 = \int \biggl( 6 (2\mu +\lambda) \dfrac{r v_x v_{xt} v_{tt}}{r_x} + 12 (2 \mu + \lambda)\dfrac{r_x v v_t v_{tt}}{r} \\
		& ~~~~ + 3 ( 2\mu + \lambda) \dfrac{r^2 v_x v_{xt} v_{xtt}}{r_x^2} + 6\lambda v v_t v_{xtt} - 4 (2\mu +\lambda) \dfrac{r v_x^3 v_{tt}}{r_x^2}\\
		& ~~~~ - 8 (  \lambda + 3\mu ) \dfrac{r_x v^3 v_{tt}}{r^2} - 2 (2\mu + \lambda) \dfrac{r^2 v_x^3 v_{xtt}}{r_x^3} - 4\lambda \dfrac{v^3 v_{xtt}}{r} \\
		& ~~~~ - 12 \mu v v_{xt} v_{tt} - 12 \mu v_x v_t v_{tt} + 24 \mu \dfrac{v^2v_xv_{tt}}{r}\biggr)\,dx,\\
		& L_2 = \int (r^2 v_{tt})_x \left\lbrack\left(\dfrac{x^2\rho_0}{r^2r_x}\right)^\gamma\right\rbrack_{tt} \,dx,\\
		& L_3 = - \int 4 (rvv_{tt})_x\left( \mathfrak{B}_t - \left\lbrack\left(\dfrac{x^2\rho_0}{r^2r_x}\right)^\gamma\right\rbrack_{t}\right)\,dx + 16 \mu \int r v \left(\dfrac{v}{r}\right)_{xt}v_{tt}\,dx \\
		& ~~~~ -\int 2 ((rv)_tv_{tt})_x \left( \mathfrak{B} - \left(\dfrac{x^2\rho_0}{r^2r_x}\right)^\gamma\right)\,dx + 8 \mu \int (rv)_t\left(\dfrac{v}{r}\right)_x v_{tt}\,dx.
	\end{align*}
%\end{displaymath}
H\"older inequaltiy together with \eqref{PrAssum} then yields 
\begin{equation}
	\begin{aligned}
		& L_1,L_2,L_3 \leq \delta \int r^2 r_x\left(\dfrac{v_{xtt}^2}{r_x^2} + \dfrac{v_{tt}^2}{r^2}\right)\,dx + C_\delta \int r^2 r_x\left(\dfrac{v_{xt}^2}{r_x^2} + \dfrac{v_{t}^2}{r^2}\right)\,dx\\
		& ~~~~~~ + C_\delta \int r^2 r_x\left(\dfrac{v_{x}^2}{r_x^2} + \dfrac{v^2}{r^2}\right)\,dx.
	\end{aligned}
\end{equation}
From \eqref{lm1}, \eqref{lm2}, integration in the temporal variable of \eqref{rg1201} yields \eqref{lm12}.

\end{pf}

\begin{lm} There is a constant $ C>0 $ %depending on $ \eta, \alpha, \beta $ 
such that
	\begin{align}
%	& \left. v_t \right|_{x=1} \leq C  \label{lm1301} \\
%	& \int \left. v_t^2 \right|_{x=1} \,dt \leq C  \label{lm1302}\\
	& \int \left. v_{tt}^2 \right|_{x=1}\,dt \leq C (\mathcal E_0 + \mathcal E_1 + \mathcal E_3 ). \label{lm1303}
	\end{align}
\end{lm}

\begin{pf}
It follows from similar arguments as in Lemma \ref{lm:boundary}. Proof is omitted here.	
\end{pf}

\begin{lm} There is a polynomial $ P = P(y_1, y_2, y_3, y_4, y_5) $ such that,
\begin{equation}\label{lm1401}
	\dfrac{1}{2} \int \left(\dfrac{x}{r}\right)^2 \rho_0 v_{tt}^2 \,dx + \int \int \left(\dfrac{v_{xtt}^2}{r_x} + r_x\dfrac{v_{tt}^2}{r^2}\right) \,dx\,dt \leq P(\mathcal E_0, \mathcal E_1, \mathcal E_2, \mathcal E_3, \mathcal E_4 ).
\end{equation}
In particular, for any $ T> 0 $, there is a constant $ C_T > 0 $ such that % depending on $ \mu, \lambda, \eta, \bar\rho, \alpha, \beta, T, \mathcal E_0, \mathcal E_1, \mathcal E_2, \mathcal E_3, \mathcal E_4  $,
	\begin{gather}
%		&  \dfrac{1}{2} \int \left(\dfrac{x}{r}\right)^2 \rho_0 v_{tt}^2 \,dx + \int \,dt \int \left(\dfrac{v_{xtt}^2}{r_x} + r_x\dfrac{v_{tt}^2}{r^2}\right) \,dx \leq C \label{lm1401}\\
 		\int \rho_0 v_{tt}^2\, dx + \int_0^T \int \left( v_{xtt}^2 + \dfrac{v_{tt}^2}{x^2} \right)\,dx  \,dt \leq C_T, \label{lm1403}	\\
 		\int v_{xt}^2 + \left(\dfrac{v_t}{x}\right)^2 \,dx \leq C_T.\label{lm1402}
	 \end{gather}
\end{lm}

\begin{pf}
	Taking the first derivative in the temporal variable of \eqref{LgNS01} yields
	\begin{equation}\label{LgNS02}
		\begin{aligned}
			& \left( \dfrac{x}{r}\right)^2 \rho_0 v_{ttt} + \left\lbrack\left(\dfrac{x^2\rho_0}{r^2r_x}\right)^\gamma\right\rbrack_{xtt} = (2\mu +\lambda) \left\lbrack\dfrac{(r^2v)_x}{r^2r_x} \right\rbrack_{xtt}\\
			& ~~~~~~ + 2 \left(\dfrac{x^2}{r^2}\dfrac{v}{r}\rho_0v_t\right)_t + 2 \left(\dfrac{x}{r}\right)^2\dfrac{v}{r} \rho_0 v_{tt}.
		\end{aligned}
	\end{equation}
	Multiply \eqref{LgNS02} with $ v_{tt} $ and integrate the resulting in the spatial variable. It follows after integration by parts that
	\begin{equation}\label{rg1407}
		\dfrac{d}{dt} \dfrac{1}{2} \int \left(\dfrac{x}{r}\right)^2 \rho_0 v_{tt}^2 \,dx + (2\mu+\lambda)\int \left( \dfrac{v_{xtt}^2}{r_x} + r_x \dfrac{v_{tt}^2}{r^2}\right)\,dx = L_1 + L_2 + L_3 + L_4 + L_5,
	\end{equation}
	with
	\begin{displaymath}
		\begin{aligned}
			& L_1 = \int \left(\dfrac{x}{r}\right)^2 \dfrac{v}{r}\rho_0 v_{tt}^2\,dx,\\
			& L_2 = (2\mu+\lambda) \int \left( - 2 \left(2 \dfrac{v^3}{r^3} v_{xtt} + \dfrac{v_x^3}{r_x^3}v_{xtt}\right) + 3 \left( 2 \dfrac{v_t v v_{xtt}}{r^2} + \dfrac{v_x v_{xt} v_{xtt}}{r_x^2}\right)\right)\,dx,\\
			% & ~~~~~~ + (2\mu+\lambda) \int \chi' \dfrac{v_{tt}^2}{r}\,dx - \int \chi' \left( \dfrac{(r^2v)_x}{r^2r_x}\right)_{tt} v_{tt}\,dx \\\
			& L_3 = \int  v_{xtt}\left\lbrack\left(\dfrac{x^2\rho_0}{r^2r_x}\right)^\gamma\right\rbrack_{tt}\,dx, ~~~~ L_4 = 2 \int \left(\dfrac{x^2}{r^2}\dfrac{v}{r}\rho_0v_t\right)_t v_{tt}\,dx, \\
			& L_5 = - (2\mu +\lambda) \left. \dfrac {v_{tt}^2} {r} \right|_{x=1} + 4\mu \left. \left( \dfrac v r \right)_{tt} v_{tt} \right|_{x=1}. 
		\end{aligned}
	\end{displaymath}
%	By using \eqref{lm903}, \eqref{lm904}, $ L_1, L_2, L_3 $ can be bounded by the following
%	\begin{align}
%		& \int \chi \left( r_x v_{tt}^2 + \dfrac{vv_{xtt}}{r} + \dfrac{v_tv_{xtt}}{r} + \dfrac{v_xv_{xtt}}{r_x}  + \dfrac{v_{xt}v_{xtt}}{r_x}  \right) \,dx  \label{rg1401} \\
%		& \int \chi' \left(  \dfrac{v v_{tt}}{r} + \dfrac{v_t v_{tt}}{r} + \dfrac{v_{tt}^2}{r}   + \dfrac{v_xv_{tt}}{r_x}  + \dfrac{v_{xt}v_{tt}}{r_x} + \dfrac{v_{xtt}v_{tt}}{r_x} \right) \,dx. \label{rg1402}
%	\end{align}
	Meanwhile,
	\begin{equation}
		L_4 = 2 \int  \dfrac{x^2}{r^2}\dfrac{v}{r} \rho_0 v_{tt}^2\,dx - 6\int  \dfrac{x^2}{r^2}\dfrac{v^2}{r^2} \rho_0 v_t v_{tt}\,dx + 2 \int  \dfrac{x^2}{r^2}\rho_0 \dfrac{v_t^2}{r} v_{tt}\,dx.
	\end{equation}
%	where the first two terms can be bounded similarly by using \eqref{lm904}. 
	Additionally, from \eqref{LgNS},
	\begin{equation}
		\begin{aligned}
			& \int \dfrac{x^2}{r^2}\rho_0 \dfrac{v_t^2}{r} v_{tt}\,dx = \int  \dfrac{v_t v_{tt}}{r} \left\lbrack (2\mu+\lambda)\left(\dfrac{(r^2v)_x}{r^2r_x}\right)_x - \left( \left( \dfrac{x^2\rho_0}{r^2r_x}\right)^\gamma \right)_x \right\rbrack\,dx\\
			& ~~~~ = \int \left( \dfrac{v_t v_{tt}}{r} \right)_x \cdot A\,dx + \left.4\mu \dfrac{vv_t v_{tt}}{r^2} \right|_{x=1} \\
			& ~~~~ = \int \left( \dfrac{v_{xt}v_{tt}}{r} + \dfrac{v_tv_{xtt}}{r} - r_x \dfrac{v_tv_{tt}}{r^2}  \right)\cdot A \,dx + \left.4\mu \dfrac{vv_t v_{tt}}{r^2} \right|_{x=1},
		\end{aligned}
	\end{equation}
	where, from \eqref{PrAssum},
	\begin{displaymath}
	\begin{aligned}
		& A = - (2\mu + \lambda) \left(\dfrac{v_x}{r_x} + 2 \dfrac{v}{r} \right) + \left( \dfrac{x^2\rho_0}{r^2r_x}\right)^\gamma \leq C, \\
		& \left.4\mu \dfrac{vv_t v_{tt}}{r^2} \right|_{x=1} \leq C \left. ( v_{t}^2 + v_{tt}^2 ) \right|_{x=1}.
	\end{aligned}
	\end{displaymath}
%	Therefore, $ L_4 $ can be bounded by
%	\begin{equation}\label{rg1403}
%		\int \chi \left( r_x v_{tt} + r_x v_t v_{tt} + r_x \dfrac{v_t v_{tt}}{r^2} + \dfrac{v_{xt}v_{tt}}{r} + \dfrac{v_t v_{xtt}}{r} \right)\,dx + \int \chi' \dfrac{v_t v_{tt}}r \,dx.
%	\end{equation}
%	Collecting \eqref{rg1401}, \eqref{rg1402}, \eqref{rg1403} and dividing them into three classes, i.e.
Therefore, $ L_1, L_2, L_3, L_4 $ can be bounded by the following, 
%	\begin{align}
%	& \int \left( r_x v_{tt}^2 + r_x v_{t} v_{tt} \right) \,dx \label{rg1404} \\
%	& \int \left( \dfrac{v^2 v_{xtt}}{r^2} + \dfrac{v_x^2 v_{xtt}}{r_x^2}  + \dfrac{v_t v_{xtt}}{r} + \dfrac{v_{xt}v_{xtt}}{r_x} + r_x \dfrac{v_tv_{tt}}{r^2} + \dfrac{v_{xt}v_{tt}}{r} \right) \,dx  \label{rg1405}\\
%	& \int \left( \dfrac{v v_{tt}}{r} + \dfrac{v_t v_{tt}}{r} + \dfrac{v_{tt}^2}{r}   + \dfrac{v_xv_{tt}}{r_x}  + \dfrac{v_{xt}v_{tt}}{r_x} + \dfrac{v_{xtt}v_{tt}}{r_x}  \right) \,dx. \label{rg1406}
%	\end{align}
%	Similar as before, by using \eqref{lm702} in \eqref{rg1406}, H\"older inequality yields \eqref{rg1404}, \eqref{rg1405}, and \eqref{rg1406} can be bounded by
	\begin{displaymath}
	\begin{aligned}
		& \delta \biggl\lbrace \int r_x \dfrac{v_{tt}^2}{r^2}\,dx + \int \dfrac{v_{xtt}^2}{r_x}\,dx \biggr\rbrace + C_\delta \biggl\lbrace \int r_x \dfrac{v_t^2}{r^2}\,dx + \int \dfrac{v_{xt}^2}{r_x}\,dx + \int r_x v_{tt}^2\,dx \\& ~~~~ + \int r_x v_t^2 \,dx  \biggr\rbrace + \biggl\lbrace \int r_x \dfrac{v^2}{r^2}\,dx   + \int \dfrac{v_{x}^2}{r_x}\,dx \biggr\rbrace\biggl\lbrace \left\Arrowvert \dfrac{v}{r}\right\Arrowvert_{L_x^\infty}^2 + \left\Arrowvert \dfrac{v_x}{r_x}\right\Arrowvert_{L_x^\infty}^2 \\&  ~~~~ + \left\Arrowvert \left(\dfrac{x^2\rho_0}{r^2r_x}\right)^{2\gamma} \right\Arrowvert_{L_x^\infty} \biggr\rbrace  + C \left. ( v_{t}^2 + v_{tt}^2 ) \right|_{x=1}.
	\end{aligned}
	\end{displaymath}
	Meanwhile, from \eqref{lm402},
	\begin{displaymath}
		L_5 \leq % C \left. (v_{tt}^2 + v^2v_t^2 + v^6) \right|_{x=1} \leq C \left. (v_{tt}^2 + (E_0+E_1) v_t^2 + (E_0+E_1)^2 v^2 ) \right|_{x=1}.
		C  ( v^2 + v_{t}^2 + v_{tt}^2 ) \bigr|_{x=1}.
	\end{displaymath}
	Hence, integration in temporal variable of \eqref{rg1407}, together with \eqref{lm1}, \eqref{lm2}, \eqref{lm501}, \eqref{lm503}, \eqref{lm6}, \eqref{lm801}, \eqref{lm12},  \eqref{lm401}, \eqref{lm403}, and \eqref{lm1303} yields \eqref{lm1401}. \eqref{lm1402} and \eqref{lm1403} are consequence of \eqref{lm1401}, \eqref{lm1001}, \eqref{lm1002}.
\end{pf}

\begin{lm} There is a constant $ C_T >0 $ such that,  %depending on $ \mu, \lambda, \bar\rho, \alpha, \beta, T $ such that
	\begin{align}
		\int v_{xxt}^2\,dx + \int \left\lbrack\left(\dfrac{v_t}{x}\right)_x\right\rbrack^2 \,dx \leq C_T. \label{lm1404}
	\end{align}
\end{lm}

\begin{pf}
	Taking one temporal derivative of \eqref{rg1001} yields
	\begin{equation}\label{rg1501}
		\begin{aligned}
		& (2\mu + \lambda)\mathcal{G}_{xtt} + \gamma\left(\dfrac{x^2\rho_0}{r^2r_x}\right)^\gamma \mathcal{G}_{xt} = \gamma^2 \left( \dfrac{x^2 \rho_0}{r^2r_x}\right)^\gamma \left( \dfrac{v_x}{r_x} + 2\dfrac{v}{r} \right) \mathcal{G}_{x} \\
		& ~~~~ + \left( \dfrac{x}{r} \right)^2 \rho_0 v_{tt} - 2 \dfrac{x^2}{r^2}\dfrac{v}{r} \rho_0 v_t - \gamma (\rho_0^{\gamma})_x \left(\dfrac{x^2}{r^2r_x}\right)^\gamma \left( \dfrac{v_x}{r_x} + 2 \dfrac{v}{r} \right).
		\end{aligned}
	\end{equation}
	By \eqref{PrAssum}, \eqref{lm1001}, \eqref{rg1002}, \eqref{rg1003}, \eqref{lm1401}, \eqref{lm801}, it holds,
	\begin{equation}\label{rg1502}
		\int \mathcal{G}_{xtt}^2\,dx \leq C_T.
	\end{equation}
	Meanwhile
	\begin{displaymath}
		\begin{aligned}
		& \mathcal{G}_{xtt} = \dfrac{x}{rr_x} \left( 2 r_x \left(\dfrac{v_t}{x}\right)_x + v_{xxt}\left(\dfrac{r}{x}\right)\right) + l_1 + l_2 + l_3 + l_4 + l_5,
		\end{aligned}
	\end{displaymath}
	with
	\begin{displaymath}
		\begin{aligned}
			& l_1 = \dfrac{x}{rr_x} \left\lbrack 2\left(2v_x\left(\dfrac{v}{x}\right)_x + v_{xx}\left(\dfrac{v}{x}\right)\right) + 2 v_{xt}\left(\dfrac{r}{x}\right)_x + r_{xx}\left(\dfrac{v_t}{x}\right) \right\rbrack,\\
			& l_2 = - \dfrac{x}{rr_x} \left(\dfrac{v_x}{r_x} + \dfrac{v}{r}\right)\left\lbrack 2r_x \left(\dfrac{v}{x}\right)_x + v_{xx}\left(\dfrac{r}{x}\right) + 2 v_x \left(\dfrac{r}{x}\right)_x + r_{xx}\left(\dfrac{v}{x}\right) \right\rbrack,\\
			& l_3 = \dfrac{x}{rr_x} \left( \dfrac{v_x}{r_x} + \dfrac{v}{r}\right)^2 \left\lbrack 2 r_x \left(\dfrac{r}{x}\right)_x + r_{xx}\left(\dfrac{r}{x}\right)\right\rbrack,\\
			& l_4 = -\dfrac{x}{rr_x} \left\lbrack \dfrac{v_{xt}}{r_x} + \dfrac{v_t}{r} - \left( \left(\dfrac{v_x}{r_x}\right)^2 + \left(\dfrac{v}{r}\right)^2 \right) \right\rbrack \left\lbrack2r_x\left(\dfrac{r}{x}\right)_x +r_{xx}\left(\dfrac{r}{x}\right)\right\rbrack ,\\
			& l_5 = - \dfrac{x}{rr_x}\left(\dfrac{v_x}{r_x}+\dfrac{v}{r}\right)\left\lbrack 2 v_x \left(\dfrac{r}{x}\right)_x + v_{xx}\left(\dfrac{r}{x}\right) +2 r_x \left(\dfrac{v}{x}\right)_x +r_{xx}\left(\dfrac{v}{x}\right)\right\rbrack,
		\end{aligned}
	\end{displaymath}
		satisfying 
		\begin{equation}
			\begin{aligned}
			& |l_1|, |l_2|, |l_3|, |l_4|, |l_5| \leq C_T \left( \left(\dfrac{v}{x}\right)_x+ v_{xx} + \left(\dfrac{r}{x}\right)_x + r_{xx} \right.\\
			 & ~~~~~~~~ \left.+ \left( v_{xt} + \dfrac{v_t}{x}\right)\left(\left(\dfrac{r}{x}\right)_x + r_{xx} \right)  \right).
			 \end{aligned}
		\end{equation}
		Therefore, from \eqref{rg1502}, \eqref{lm1102} and \eqref{lm1101},
		\begin{equation}\label{rg1503}
%		\begin{aligned}
			 \int \left(  2 r_x \left(\dfrac{v_t}{x}\right)_x + v_{xxt}\left(\dfrac{r}{x}\right) \right)^2\,dx \leq C_T + C_T \left( \left\Arrowvert v_{xt}\right\Arrowvert_{L_x^\infty}^2 + \left\Arrowvert \dfrac{v_t}{x}\right\Arrowvert_{L_x^\infty}^2 \right) .
%			& ~~~~ \leq C_T + \delta \left( \int v_{xxt}^2 \,dx + \int \left(\dfrac{v_t}{x}\right)_x^2 \,dx \right)
%		\end{aligned}
		\end{equation}
%		where the last inequality is due to the fact
%		\begin{displaymath}
%		f^2(x) \leq \int_0^1 f^2\,dx + \int 2 f f_x \,dx \leq \delta \int f_x^2\,dx + (1+C_\delta)\int f^2	 \,dx
%%			\left\Arrowvert v_{xt}\right\Arrowvert_{L_x^\infty}^2 + \left\Arrowvert \dfrac{v_t}{x}\right\Arrowvert_{L_x^\infty}^2 \leq  \int \left( C_\delta v_{xt}^2 + C_\delta \left(\dfrac{v_t}{x}\right)^2 + \delta v_{xxt}^2 + \delta \left(\dfrac{v_t}{x}\right)_x^2  \right)\,dx
%		\end{displaymath}
		Meanwhile,
		\begin{displaymath}\begin{aligned}
			& \int \left(  2 r_x \left(\dfrac{v_t}{x}\right)_x + v_{xxt}\left(\dfrac{r}{x}\right) \right)^2\,dx = \int \left(\dfrac{r}{x}\right)^2 v_{xxt}^2\,dx + 4 \int r_x^2 \left(\dfrac{v_t}{x}\right)_x^2\,dx \,dx\\
			& ~~~~ + 4 \int r_x \left(\dfrac{r}{x}\right) \left(\dfrac{v_t}{x}\right)_x v_{xxt}\,dx.
			\end{aligned}
		\end{displaymath}
		Similarly as before, since $ v_{xxt} = x \left(\dfrac{v_t}{x}\right)_{xx} + 2 \left(\dfrac{v_t}{x}\right)_x $, 
%		 $r_x = \frac{r}{x}  + x \left(\dfrac{r}{x}\right)_x $ and 
		\begin{displaymath} \begin{aligned}
			& \int r_x \left(\dfrac{r}{x}\right) \left(\dfrac{v_t}{x}\right)_x v_{xxt}\,dx = \dfrac{1}{2}\int x r_x \left(\dfrac{r}{x}\right) \left\lbrack\left(\dfrac{v_t}{x}\right)_{x}^2\right\rbrack_{x} \,dx + 2 \int r_x \left(\dfrac{r}{x}\right) \left(\dfrac{v_t}{x}\right)_x^2\,dx\\
			& = \left. \left(  \dfrac{1}{2}x r_x \left(\dfrac{r}{x}\right) \left(\dfrac{v_t}{x}\right)_x^2  \right)\right|_{x=1} + \dfrac{3}{2}\int r_x \left(\dfrac{r}{x}\right)\left(\dfrac{v_t}{x}\right)_x^2\,dx - \dfrac{1}{2} \int x r_{xx}\left(\dfrac{r}{x}\right) \left(\dfrac{v_t}{x}\right)_x^2\,dx\\
			& ~~~~ - \dfrac{1}{2} \int x r_x \left(\dfrac{r}{x}\right)_x\left(\dfrac{v_t}{x}\right)_x^2\,dx.
			\end{aligned}
		\end{displaymath}
		\eqref{rg1503} then implies
		\begin{equation}\label{rg1504}
			\begin{aligned}
				& \int \left(\dfrac{r}{x}\right)^2 v_{xxt}^2\,dx  + \int \left( 4 r_x^2 + 6 r_x \left(\dfrac{r}{x}\right) - 2 x r_x \left(\dfrac{r}{x}\right)_x\right) \left(\dfrac{v_t}{x}\right)_x^2\,dx \\
				& ~~~~ \leq C_T + 2 \int x r_{xx} \left(\dfrac{r}{x}\right)\left(\dfrac{v_t}{x}\right)_x^2\,dx + C_T \left( \left\Arrowvert v_{xt}\right\Arrowvert_{L_x^\infty}^2 + \left\Arrowvert \dfrac{v_t}{x}\right\Arrowvert_{L_x^\infty}^2 \right).
			\end{aligned}
		\end{equation}
		Then from \eqref{rg1101} again, it follows
		\begin{equation*}
			\begin{aligned}
				& \int \left(\dfrac{r}{x}\right)^2v_{xxt}^2\,dx + 9 \int \left(\dfrac{r}{x}\right)^2 \left(\dfrac{v_t}{x}\right)_x^2\,dx \leq C_T + C \int x^2 \left(\dfrac{r}{x}\right)_x^2\left(\dfrac{v_t}{x}\right)_x^2\,dx \\
				& ~~~~ +  2 \int x r_{xx} \left(\dfrac{r}{x}\right)\left(\dfrac{v_t}{x}\right)_x^2\,dx + C_T \left( \left\Arrowvert v_{xt}\right\Arrowvert_{L_x^\infty}^2 + \left\Arrowvert \dfrac{v_t}{x}\right\Arrowvert_{L_x^\infty}^2 \right)\\
				& \leq C_T  + \delta\int \left(\dfrac{r}{x}\right)^2 \left(\dfrac{v_t}{x}\right)_x^2\,dx + C_\delta \left\Arrowvert x \left(\dfrac{v_t}{x}\right)_x\right\Arrowvert_{L_x^\infty}^2 \left( \int r_{xx}^2 \,dx + \int \left(\dfrac{r}{x}\right)_x^2\,dx \right) \\ & ~~~~ + C_T \left( \left\Arrowvert v_{xt}\right\Arrowvert_{L_x^\infty}^2 + \left\Arrowvert \dfrac{v_t}{x}\right\Arrowvert_{L_x^\infty}^2 \right)\\
				& \leq C_T + \delta\int \left(\dfrac{r}{x}\right)^2 \left(\dfrac{v_t}{x}\right)_x^2\,dx + C_T\left( \left\Arrowvert v_{xt}\right\Arrowvert_{L_x^\infty}^2 + \left\Arrowvert \dfrac{v_t}{x}\right\Arrowvert_{L_x^\infty}^2 \right),
			\end{aligned}
		\end{equation*}
	where the last inequality is due to \eqref{lm1102} and
	\begin{displaymath}
		x \left(\dfrac{v_t}{x}\right)_x = v_{xt} - \dfrac{v_t}{x} .	\end{displaymath}
	Therefore, by suitably small $ \delta > 0 $, 
	\begin{equation*}
		\int \left(\dfrac{r}{x}\right)^2v_{xxt}^2\,dx + \int \left(\dfrac{r}{x}\right)^2 \left(\dfrac{v_t}{x}\right)_x^2\,dx \leq C_T + C_T\left( \left\Arrowvert v_{xt}\right\Arrowvert_{L_x^\infty}^2 + \left\Arrowvert \dfrac{v_t}{x}\right\Arrowvert_{L_x^\infty}^2 \right),
	\end{equation*}
	or, together with \eqref{lm1001},
	\begin{equation} \begin{aligned}
		& \int v_{xxt}^2\,dx + \int \left(\dfrac{v_t}{x}\right)_x^2 \,dx \leq C_T + C_T\left( \left\Arrowvert v_{xt}\right\Arrowvert_{L_x^\infty}^2 + \left\Arrowvert \dfrac{v_t}{x}\right\Arrowvert_{L_x^\infty}^2 \right)\\
		& \leq (1+C_\delta )C_T + \delta \left(\int v_{xxt}^2\,dx + \int \left(\dfrac{v_t}{x}\right)_x^2 \,dx\right),
		\end{aligned}
	\end{equation}
	where the last inequality is from \eqref{lm1402} and 
	\begin{equation}\label{rg1505}
		f^2(x) \leq \int_0^1 f^2\,dx + \int 2 f f_x \,dx \leq \delta \int f_x^2\,dx + (1+C_\delta)\int f^2	 \,dx.
	\end{equation}
	\eqref{lm1404} is then obtained after choosing $ \delta $ small enough. 
\end{pf}

\begin{lm} There is a constant $ C_T > 0 $,
	\begin{equation}\label{rg1506}
		\left\Arrowvert v_t , v_{xt},  \dfrac{v_t}{x}  \right\Arrowvert_{L_x^\infty}^2 \leq C_T.
	\end{equation}
\end{lm}

\begin{pf}
	This is a direct consequence of \eqref{lm1402}, \eqref{lm1404}.
\end{pf}

\begin{lm}\label{lm:D3regu} There is a constant $ C_T > 0 $ such that, 
\begin{align}
& \int r_{xxx}^2 \,dx + \int \left(\dfrac{r}{x}\right)_{xx}^2\,dx \leq C_T, \label{lm1601}\\
& \int v_{xxx}^2 \,dx + \int \left(\dfrac{v}{x}\right)_{xx}^2\,dx \leq C_T .\label{lm1602}	
\end{align}

\end{lm}

\begin{pf}
Take one spatial derivative of \eqref{rg1001}. The resulting equation is
\begin{equation}\label{rg1600}
	\begin{aligned}
		& (2\mu+\lambda)\mathcal{G}_{xxt} + \gamma \left(\dfrac{x^2 \rho_0}{r^2r_x}\right)^\gamma \mathcal{G}_{xx} = -2 \dfrac{x^3}{r^3} \left(\dfrac{r}{x}\right)_x \rho_0 v_t + \dfrac{x^2}{r^2} \rho_0 v_{xt} + \dfrac{x^2}{r^2} v_{t} (\rho_0)_x \\
		& ~~~~~~ + \left\lbrack \left(\dfrac{x^2\rho_0}{r^2r_x}\right)^\gamma \right\rbrack_x \left( - \gamma \mathcal{G}_x + (\rho_0^\gamma)_{x} \right) + \left(\dfrac{x^2\rho_0}{r^2r_x}\right)^\gamma (\rho_0^\gamma)_{xx}.
	\end{aligned}
\end{equation}
Then it follows
\begin{equation}\label{rg1601}
	\begin{aligned}
		& \dfrac{d}{dt} \int \mathcal{G}_{xx}^2 \,dx + \int \left(\dfrac{x^2\rho_0}{r^2r_x}\right)^\gamma\mathcal{G}_{xx}^2\,dx \leq \int \mathcal{G}_{xx}^2\,dx + C_T \int \left\lbrack \left( \dfrac{r}{x}\right)_x \right\rbrack^2\,dx \\
		& ~~~~~~ + C_T \int v_{xt}^2\,dx + C_T \int ((\rho_0)_x)^2\,dx + C_T \int ((\rho_0^\gamma)_{xx})^2 \,dx\\
		& ~~~~~~ + C_T \int \left\lbrack ((\rho_0^\gamma)_x)^2 + r_{xx}^2 + \left(\dfrac{r}{x}\right)_x^2 \right\rbrack \cdot \left( \mathcal{G}_x^2 + (\rho_0^\gamma)_x^2 \right)\,dx, \\
	\end{aligned}
\end{equation}
where \eqref{rg1506} ,  \eqref{PrAssum}, \eqref{lm1001} are applied above.
Meanwhile, notice
\begin{align*}
	& ((\rho_0^\gamma)_x)^2 \leq C \left( \int ((\rho_0)_x)^2\,dx + \int ((\rho_0^\gamma)_{xx})^2\,dx \right), \\
	& \mathcal{G}_x^2 \leq C \left( \int \mathcal{G}_x^2\,dx + \int \mathcal{G}_{xx}^2\,dx \right).
\end{align*}
\eqref{rg1601} together with \eqref{lm1102} and \eqref{lm1402} yields
\begin{equation}\begin{aligned}
	&	\dfrac{d}{dt} \int \mathcal{G}_{xx}^2\,dx + \int \left(\dfrac{x^2\rho_0}{r^2r_x}\right)^\gamma\mathcal{G}_{xx}^2\,dx \leq C_T + C_T\int\mathcal{G}_{xx}^2\,dx \\
	& ~~~~~~ + C_T \left\lbrack \int \left((\rho_0)_x\right)^2\,dx + \int (( \rho_0^\gamma)_{xx})^2 \,dx  \right\rbrack.
	\end{aligned}
\end{equation}
Gr\"onwall's inequality then shows
\begin{equation}\label{rg1602}
	\int \mathcal{G}_{xx}^2 \,dx \leq C_T.
\end{equation}
Consequently, from \eqref{rg1600} it holds,
\begin{equation}\label{rg1603}
	\int \mathcal{G}_{xxt}^2 \,dx \leq C_T.
\end{equation}
Notice,
\begin{align}
	& \mathcal{G}_{xx} = \dfrac{x}{rr_x}\left( 2r_x \left(\dfrac{r}{x}\right)_{xx} + r_{xxx}\left(\dfrac{r}{x}\right) \right) \nonumber \\
	& ~~~~ + \left(\dfrac{x}{rr_x}\right) \cdot 3r_{xx}\left(\dfrac{r}{x}\right)_x + \left( \dfrac{x}{rr_x}\right)_x \left( 2r_x\left(\dfrac{r}{x}\right)_x + r_{xx}\left(\dfrac{r}{x}\right)\right), \label{rg1604}  \\
	& \mathcal{G}_{xxt} = \dfrac{x}{rr_x}\left( 2r_x \left(\dfrac{v}{x}\right)_{xx} + v_{xxx}\left(\dfrac{r}{x}\right) \right) + l_1 + l_2 + l_3 + l_4, \label{rg1605} 
\end{align}
where
\begin{equation}\label{rg1606}
	\begin{aligned}
	& l_1 = \dfrac{x}{rr_x}\left( 2v_x \left(\dfrac{r}{x}\right)_{xx} + r_{xxx}\left(\dfrac{v}{x}\right) \right) - \dfrac{x}{rr_x}\left(\dfrac{v}{r} + \dfrac{v_x}{r_x} \right)\left( 2r_x \left(\dfrac{r}{x}\right)_{xx} \right.\\& \left. ~~~~~~~~~~  + r_{xxx}\left(\dfrac{r}{x}\right) \right), \\
	& l_2 = 3 \left(\dfrac{x}{rr_x}\right) \left( v_{xx}\left(\dfrac{r}{x}\right)_x + r_{xx}\left(\dfrac{v}{x}\right)_x\right) - 3 \left(\dfrac{x}{rr_x}\right)\left(\dfrac{v}{r}+\dfrac{v_x}{r_x}\right) r_{xx}\left(\dfrac{r}{x}\right)_x,\\
	& l_3 = \left(\dfrac{x}{rr_x}\right)_x \left( 2 v_x \left(\dfrac{r}{x}\right)_x + 2 r_x \left(\dfrac{v}{x}\right)_x + v_{xx}\left(\dfrac{r}{x}\right) + r_{xx} \left(\dfrac{v}{x}\right)\right),\\
	& l_4 = - \left( \dfrac{x^2}{r^2r_x^2} \left(r_x\dfrac{v}{x} + \dfrac{r}{x} v_x\right) \right)_x\left( 2r_x\left(\dfrac{r}{x}\right)_x + r_{xx}\left(\dfrac{r}{x}\right)\right).
	\end{aligned}
\end{equation}
Then, \eqref{rg1602}, \eqref{rg1604}, \eqref{lm1001} and \eqref{lm1102} yield 
\begin{equation}\label{rg1607}
	\begin{aligned}
		& \int \left( 2 r_x \left(\dfrac{r}{x}\right)_{xx}+r_{xxx}\left(\dfrac{r}{x}\right)\right)^2\,dx \leq C_T + C_T \int \left( r_{xx}^2 + \left(\dfrac{r}{x}\right)_x ^2\right) \\
		& ~~ \times \left( r_{xx}^2 + \left(\dfrac{r}{x}\right)_x^2\right) \,dx  \leq C_T + C_T\left( \left\Arrowvert r_{xx}^2 \right\Arrowvert_{L_x^\infty} + \left\Arrowvert \left(\dfrac{r}{x}\right)_x^2\right\Arrowvert_{L_x^\infty} \right)\\
		& ~~~~ \leq (1+C_\delta )C_T + \delta\left(\int r_{xxx}^2 \,dx + \int \left(\dfrac{r}{x}\right)_{xx}^2\,dx \right),
	\end{aligned}
\end{equation}
where the last inequality is due to \eqref{rg1505}.
Meanwhile,
\begin{displaymath}
	\begin{aligned}
		& \int \left( 2 r_x \left(\dfrac{r}{x}\right)_{xx}+r_{xxx}\left(\dfrac{r}{x}\right)\right)^2 \,dx  =  4 \int r_{x}^2 \left(\dfrac{r}{x}\right)_{xx}^2 \,dx + \int \left(\dfrac{r}{x}\right)^2 r_{xxx}^2\,dx\\
		& ~~~~~~~~ + 4 \int r_x \left(\dfrac{r}{x}\right) \left(\dfrac{r}{x}\right)_{xx}r_{xxx}\,dx,
	\end{aligned}
\end{displaymath}
in which, 
\begin{displaymath}
	\begin{aligned}
		& 4 \int r_x \left(\dfrac{r}{x}\right)\left(\dfrac{r}{x}\right)_{xx}r_{xxx}\,dx = 4 \int r_x \left(\dfrac{r}{x}\right) \left(\dfrac{r}{x}\right)_{xx} \left( x \left(\dfrac{r}{x}\right)_{xxx} + 3 \left(\dfrac{r}{x}\right)_{xx}\right)\,dx \\
		& = \left. 2 x r_x \left(\dfrac{r}{x}\right) \left(\dfrac{r}{x}\right)_{xx}^2 \right|_{x=1} + 10 \int r_x \left(\dfrac{r}{x}\right) \left(\dfrac{r}{x}\right)_{xx}^2 \,dx  - 2 \int x r_{xx}\left(\dfrac{r}{x}\right) \left(\dfrac{r}{x}\right)_{xx}^2 \,dx\\
		& ~~~~~~ - 2 \int x r_x \left( \dfrac{r}{x}\right)_x \left(\dfrac{r}{x}\right)_{xx}^2\,dx.
	\end{aligned}
\end{displaymath}
In the meantime, applying \eqref{rg1505}, the boundary term above admits
\begin{displaymath}
	\begin{aligned}
		& \left. x r_x \left(\dfrac{r}{x}\right) \left(\dfrac{r}{x}\right)_{xx}^2 \right|_{x=1} = \left. r r_x \left(r_{xx} - 2 r_{x} + 2 r\right)^2 \right|_{x=1} \\
		& ~~~~~~ \leq C_T + C_T  \left\Arrowvert r_{xx}^2 \right\Arrowvert_{L_x^\infty} \leq (1+C_\delta) C_T + \delta\int r_{xxx}^2\,dx.
	\end{aligned}
\end{displaymath}
Therefore, from \eqref{lm1102} and \eqref{rg1607} it follows
\begin{equation}\label{rg1609}
	\begin{aligned}
		& \int \left(\dfrac{r}{x}\right)^2 r_{xxx}^2\,dx + \int \left(4r_x^2 + 10 r_x \left(\dfrac{r}{x}\right) - 2 x r_x \left(\dfrac{r}{x}\right)_x \right) \left(\dfrac{r}{x}\right)_{xx}^2\,dx \\
		&  \leq (1+C_\delta)C_T + 2 \int xr_{xx}\left(\dfrac{r}{x}\right) \left(\dfrac{r}{x}\right)_{xx}^2  +   \delta\left(\int r_{xxx}^2 \,dx + \int \left(\dfrac{r}{x}\right)_{xx}^2\,dx \right)\\
		& \leq (1+C_\delta)C_T + C_\delta C_T \left\Arrowvert x\left(\dfrac{r}{x}\right)_{xx} \right\Arrowvert_{L_x^\infty}^2 \int r_{xx}^2\,dx \\
		& ~~~~ + \delta\left(\int r_{xxx}^2 \,dx + \int \left(\dfrac{r}{x}\right)_{xx}^2\,dx \right)\\
		& \leq (1+C_\delta)C_T + \delta \left(\int r_{xxx}^2 \,dx + \int \left(\dfrac{r}{x}\right)_{xx}^2\,dx \right),
	\end{aligned}
\end{equation}
where in the last inequality it is applied
\begin{equation}\label{rg1610}
	\begin{aligned}
	& \left\Arrowvert x\left(\dfrac{r}{x}\right)_{xx} \right\Arrowvert_{L_x^\infty}^2 = \left\Arrowvert r_{xx} - 2 \left(\dfrac{r}{x}\right)_x \right\Arrowvert_{L_x^\infty}^2 \leq (1+C_\delta )C_T \\
	& ~~~~ + \delta \int \left( r_{xxx}^2 + \left(\dfrac{r}{x}\right)_{xx}^2 \right) \,dx .
	\end{aligned}
\end{equation}
Again, since
\begin{equation}\label{rg1608}
	\begin{aligned}
		& 4r_x^2 + 10 r_x \left(\dfrac{r}{x}\right) - 2 x r_x \left(\dfrac{r}{x}\right)_x = 2 x^2 \left(\dfrac{r}{x}\right)_x^2 + 16 x \left(\dfrac{r}{x}\right)\left(\dfrac{r}{x}\right)_x + 14 \left(\dfrac{r}{x}\right)^2\\
		& ~~~~ \geq 13 \left(\dfrac{r}{x}\right)^2 - C x^2 \left(\dfrac{r}{x}\right)_x^2,
	\end{aligned}
\end{equation}
it can be derived from \eqref{rg1609},
\begin{equation}
	\begin{aligned}
		& \int \left(\dfrac{r}{x}\right)^2 r_{xxx}^2 \,dx + 13 \int \left(\dfrac{r}{x}\right)^2 \left(\dfrac{r}{x}\right)_{xx}^2\,x \leq  (1+C_\delta ) C_T \\
		& ~~~~ + \int x^2 \left(\dfrac{r}{x}\right)_x^2 \left(\dfrac{r}{x}\right)_{xx}^2 \,dx + \delta \left(\int r_{xxx}^2 \,dx + \int \left(\dfrac{r}{x}\right)_{xx}^2\,dx \right) \\
		& ~~~~ \leq (1+C_\delta ) C_T + \delta \left(\int r_{xxx}^2 \,dx + \int \left(\dfrac{r}{x}\right)_{xx}^2\,dx \right),
	\end{aligned}
\end{equation}
where \eqref{rg1610}, \eqref{lm1102} is applied again.
Therefore, \eqref{lm1601} follows by choosing small $\delta$ and \eqref{lm1001}.
Similarly, from \eqref{rg1605}, \eqref{rg1606}, \eqref{lm1001}, \eqref{lm1102}, \eqref{lm1101}, \eqref{PrAssum}, and \eqref{lm1601}, the following estimate holds,
\begin{equation}
	\begin{aligned}
		& \int \left(2 r_x \left(\dfrac{v}{x}\right)_{xx} + v_{xxx} \left(\dfrac r x \right) \right)^2 \,dx \leq C_T + C_T \left( \int r_{xxx}^2 \,dx + \int \left(\dfrac r x \right)_{xx}^2 \,dx \right. \\ & ~~~~~~ \left.  +  \int\left(v_{xx}^2 + \left(\dfrac{v}{r}\right)_x^2 + r_{xx}^2 +  \left(\dfrac{r}{x}\right)_{x}^2 \right) \cdot \left( r_{xx}^2 +  \left(\dfrac{r}{x}\right)_{x}^2  \right) \,dx \right) \\
		& \leq C_T + C_T\left( \left\Arrowvert r_{xx} \right\Arrowvert_{L_x^\infty}^2 + \left\Arrowvert \left(\dfrac{r}{x}\right)_{x} \right\Arrowvert_{L_x^\infty}^2  \right) \\
		& \leq C_T +C_T \int \left( r_{xx}^2 + r_{xxx}^2 + \left(\dfrac{r}{x}\right)_{x}^2 + \left(\dfrac{r}{x}\right)_{xx}^2 \right)\leq C_T.
	\end{aligned}
\end{equation}
Again, 
\begin{displaymath}
	\begin{aligned}
		& \int \left(2 r_x \left(\dfrac{v}{x}\right)_{xx} + v_{xxx} \left(\dfrac r x \right) \right)^2 \,dx = \int \left(\dfrac r x \right)^2 v_{xxx}^2 \,dx + 4 \int r_x^2 \left(\dfrac v x \right)_{xx}^2 \,dx \\
		& ~~~~~~~~ + 4 \int r_x \left( \dfrac r x \right) \left(\dfrac{v}{x}\right)_{xx} v_{xxx}\,dx,
	\end{aligned}
\end{displaymath}
and we shall apply integration by parts as in the following,
\begin{displaymath}
	\begin{aligned}
		& \int r_x \left( \dfrac r x \right) \left(\dfrac{v}{x}\right)_{xx} v_{xxx}\,dx = \int r_x \left(\dfrac r x \right) \left(\dfrac v x \right)_{xx} \left( x \left( \dfrac v x \right)_{xxx} + 3 \left(\dfrac v x \right)_{xx} \right)\,dx \\
		& = \left. \left( \dfrac{1}{2} x r_{x} \left(\dfrac r x \right) \left( \dfrac v x \right)_{xx}^2 \right) \right|_{x=1} + \dfrac{5}{2} \int r_x \left(\dfrac r x \right) \left( \dfrac v x \right)_{xx}^2 \,dx \\
		& ~~~~~~ - \dfrac 1 2 \int \left( x r_{xx}\left(\dfrac r x \right) + x r_{x} \left( \dfrac r x \right)_x \right) \left(\dfrac v x \right)_{xx}^2\,dx.
	\end{aligned}
\end{displaymath}
And hence it admits, 
\begin{equation}
	\begin{aligned}
		& \int \left(\dfrac r x \right)^2 v_{xxx}^2 \,dx + \int \left( 4 r_{x}^2 + 10 r_x \dfrac r x - 2 x r_x \left(\dfrac r x \right)_x \right) \left(\dfrac v x\right)_{xx}^2\,dx \\
		& \leq 2 \int x r_{xx}\left( \dfrac r x \right) \left( \dfrac v x \right)_{xx}^2 \,dx - \left. \left( 2 x r_{x} \left(\dfrac r x \right) \left( \dfrac v x \right)_{xx}^2 \right) \right|_{x=1} + C_T,
	\end{aligned}
\end{equation}
or, by noticing \eqref{rg1608},
\begin{equation}\label{rg1611}
	\begin{aligned}
		& \int \left(\dfrac r x \right)^2 v_{xxx}^2 \,dx + 13 \int \left(\dfrac r x \right)^2\left( \dfrac v x \right)_{xx}^2\,dx \leq C_T + C \int x^2 \left(\dfrac r x \right)_{x}^2 \left(\dfrac v x \right)_{xx}^2 \,dx \\
		& ~~~~~~ + 2 \int x r_{xx}\left( \dfrac r x \right) \left( \dfrac v x \right)_{xx}^2 \,dx - \left. \left( 2 x r_{x} \left(\dfrac r x \right) \left( \dfrac v x \right)_{xx}^2 \right) \right|_{x=1} \\
		& \leq (1+C_\delta)C_T + \delta\int \left(\left(\dfrac{v}{x}\right)_{xx}^2 + v_{xxx}^2 \right)\,dx,
	\end{aligned}
\end{equation}
where the last inequality is due to \eqref{lm1102}, \eqref{lm1101}, \eqref{rg1505} and 
\begin{displaymath}
	\begin{aligned}
		& \left. x r_{x} \left(\dfrac r x \right) \left( \dfrac v x \right)_{xx}^2 \right|_{x=1} = \left. r_x r \left(v_{xx} -  2 v_{x} + 2 v\right)^2\right|_{x=1} \leq C_T \left\Arrowvert v_{xx} \right\Arrowvert_{L_x^\infty} + C_T,\\
		& \left\Arrowvert x \left(\dfrac{v}{x}\right)_{xx} \right\Arrowvert_{L_x^\infty}^2 = \left\Arrowvert v_{xx} - 2 \left(\dfrac{v}{x}\right)_x \right\Arrowvert_{L_x^\infty}^2 \leq C \left\Arrowvert v_{xx} \right\Arrowvert_{L_x^\infty}^2+ C \left\Arrowvert \left( \dfrac v x \right)_x \right\Arrowvert_{L_x^\infty}^2,\\
		& \left\Arrowvert v_{xx} \right\Arrowvert_{L_x^\infty}^2 + \left\Arrowvert\left( \dfrac v x \right)_x \right\Arrowvert_{L_x^\infty}^2 \leq C_\delta C_T + \delta \int \left(v_{xxx}^2 + \left(\dfrac v x \right)_{xx}^2\right) \,dx .
	\end{aligned}
\end{displaymath}
Therefore, \eqref{lm1602} is consequence of \eqref{rg1611} and \eqref{lm1001}. 
\end{pf}

\section{Global Existence}

\subsection{Functional Framework for the Local Well-posedness Theory}

In this section, we shall discuss the appropriate functional framework for the local well-posedness theory for the system consisting of the equation \eqref{LgNS} with the boundary condition \eqref{LgNSBC} and the following initial data
\begin{equation}\label{local001}
	r(x,0) = {\bar r}(x), ~ v(x,0) = {\bar u}_0(x).
\end{equation}
%where $ \bar r, \bar v $ satisfy the conditions 
%\begin{equation}\label{local002}
%	0 < \dfrac{x^2}{\bar r^2 \bar r_x} \leq \alpha^3, ~ \left| \dfrac{\bar v}{\bar r} \right|, \left| \dfrac{{\bar v}_x}{{\bar r}_x} \right| \leq \beta. 
%\end{equation}
Moreover, the coefficient $ \rho_0 $ satisfies the condition \eqref{InitialAssum}, \eqref{InitialAssum2}. 
%\todo{oct092016}
Now we define the energy functionals. Denote
\begin{equation}\label{local003}
\begin{aligned}
	& \mathfrak E_0(t) = \int x^2 \rho_0 v^2 \,dx + \int r^2 r_x \left(\dfrac{x^2 \rho_0}{r^2r_x}\right)^\gamma\,dx  + \int x^2 \rho_0 v_t^2 \,dx \\
	& ~~~~~~~~~~~  + \int \chi \left(\dfrac{x}{r}\right)^2 \rho_0 v_t^2 \,dx,\\
	& \mathfrak E_1(t) = \int x^2 \rho_0 v_{tt}^2 \,dx + \int \chi \left( \dfrac{x}{r} \right)^2 \rho_0 v_{tt}^2 \,dx,  \\
	& \mathfrak D_{0}(t) = \int \left( r^2 \dfrac{v_x^2}{r_x} + r_x v^2 + r^2 \dfrac{v_{xt}^2}{r_x} + r_x v_{t}^2 \right)\,dx + \int \chi \left( r_x \dfrac{v_t^2}{r^2} + \dfrac{v_{xt}^2}{r_x} \right) \, dx , \\
	& \mathfrak D_1(t) = \int \left( r^2 \dfrac{v_{xtt}^2}{r_x} + r_x v_{tt}^2 \right) \,dx + \int \chi \left( r_x \dfrac{v_{tt}^2}{r^2} + \dfrac{v_{xtt}^2}{r_x} \right) \,dx,
\end{aligned}	
\end{equation}
where $ \chi $ is a cut-off function satisfying 
\begin{equation*}
	\chi = \begin{cases}
		1, & x\in (0, 1/4), \\
		0, & x \in ( 1/2, R_0) ,
	\end{cases}
\end{equation*}
and $ - 8 \leq \chi' \leq 0 $.
%==============
The following local well-posedness theorem would be applied.
\begin{lm}[Local Well-posedness]\label{lm:local}
	For the equation \eqref{LgNS} with boundary condition \eqref{LgNSBC}, there is a strong solution $ (r, v) = (r,v)(x,t) $ in $ t \in (0,T^*), ~ T^* > 0 $ with given initial data $ (r,v)(x,0) = (\bar r, \bar u_0)(x) $ satisfying
%	\begin{equation}\label{local:initial}
%		\mathfrak E_0(0) + \int_0^1 \left\lbrack \dfrac{x^2}{\bar r_0^2\bar r_{0,x}} \left( \dfrac{\bar r_0^2\bar r_{0,x}}{x^2}\right)_x \right\rbrack^2\,dx < \infty
%	\end{equation}
	\begin{equation}\label{local:initial}
		\begin{aligned}
			& \int x^2 \rho_0 \bar u_0^2\,dx + \int \bar{r}^2 \bar{r}_{x} \left(\dfrac{x^2 \rho_0}{\bar{r}^2\bar{r}_{x}}\right)^\gamma\,dx + \int x^2 \rho_0 \bar u_1^2\,dx \\
			& ~~~~~~~ + \int \chi \dfrac{x^2}{\bar r^2} \rho_0 \bar u_1^2 \,dx  + \int \left\lbrack \dfrac{x^2}{\bar r^2\bar r_{x}} \left( \dfrac{\bar r^2\bar r_{x}}{x^2}\right)_x \right\rbrack^2\,dx < \infty,
		\end{aligned}
	\end{equation}
	where %$ v_t(x,0) $ is given by the compatible condition
	\begin{equation}\label{local:initial3}
		\bar u_1 (x) = \dfrac{\bar r^2}{x^2 \rho_0}\left\lbrace(2\mu+\lambda)\left\lbrack\dfrac{(\bar r^2 \bar u_0)_x}{\bar r^2 \bar r_x}\right\rbrack_x - \left\lbrack\left(\dfrac{x^2\rho_0}{\bar r^2 \bar r_x}\right)^\gamma\right\rbrack_x \right\rbrace.
	\end{equation}
	$ (r, v) $ satisfies
%	\begin{equation}\label{local:pointwisebound}
%		\dfrac{x^2}{r^2r_x}, \dfrac{v_x}{r_x}, \dfrac{v}{r} < \infty,
%	\end{equation}
%	and 
	\begin{equation}\label{local:regularity1}
		\begin{cases}
			 r_x, v_x, \dfrac{r}{x}, \dfrac{v}{x} \in L_t^\infty((0,T^*), H_x^1(0,1)),\\
			 x \sqrt{\rho_0} v, x \sqrt{\rho_0} v_t, \sqrt{\rho_0} v_t, \in L_t^\infty((0,T^*),L_x^2(0,1)),\\
			 v, x v_x, v_t, x v_{xt}, \dfrac{v_t}{x}, v_{xt} \in L_t^2((0,T^*),L_x^2(0,1)).
		\end{cases}
	\end{equation}
	If in addition % the following additional regularity are 
	\begin{equation}\label{local:initial2}
	\begin{aligned}
		& \norm{(\rho_0)_x}{L_x^2(0,1)}, \norm{(\rho_0^\gamma)_{xx}}{L_x^2(0,1)} < \infty,  \\
		& \int x^2 \rho_0 \bar u_2^2\,dx + \int \chi \dfrac{x^2}{\bar r^2} \rho_0 \bar u_2^2\,dx + \int \left\lbrace \left\lbrack \dfrac{x^2}{\bar r^2\bar r_{x}} \left( \dfrac{\bar r^2\bar r_{x}}{x^2}\right)_x \right\rbrack_x \right\rbrace^2\,dx < \infty,
	 \end{aligned}
	\end{equation}
	where
	\begin{equation}\label{local:initial4}
	\begin{aligned}
			& \bar u_2  = \dfrac{\bar r^2}{x^2 \rho_0}\left\lbrace(2\mu+\lambda)\left\lbrack\dfrac{(\bar r^2 \bar u_1)_x}{\bar r^2 \bar r_x}\right\rbrack_x +\gamma \left\lbrack\left(\dfrac{x^2\rho_0}{\bar r^2 \bar r_x}\right)^\gamma\left(\dfrac{(\bar r^2 \bar u_0)_x}{\bar r^2 \bar r_x}\right)\right\rbrack_x \right\rbrace \\
			& ~~~~ - (2\mu + \lambda) \dfrac{\bar r^2}{x^2 \rho_0} \left\lbrack \dfrac{\bar u_{0,x}^2}{\bar r_x^2} + 2 \dfrac{\bar u_0^2}{\bar r^2} \right\rbrack_x + 2 \dfrac{\bar u_0}{\bar r} \bar u_1, \\
		\end{aligned}
	\end{equation}
	there is a smooth solution in $ t \in (0, T^{**}) $ for some $ T^{**} > 0 $, and \eqref{local:regularity1} holds. Moreover, the following regularity would also hold
	\begin{equation}
		\begin{cases}
			r_x, v_x, \dfrac{r}{x}, \dfrac{v}{x} \in L_t^\infty((0,T^{**}],H_x^2(0,1)),\\
			x \sqrt{\rho_0} v_{tt}, \sqrt{\rho_0} v_{tt}, v_{xxt}, \left(\dfrac{v_t}{x}\right)_x \in L_t^\infty ((0,T^{**}],L_x^2(0,1)),\\
			v_{tt}, x v_{xtt}, \dfrac{v_{tt}}{x}, v_{xtt} \in L_t^2 ((0,T^{**}], L_x^2(0,1)).
		\end{cases}
	\end{equation}
	In particular, 
	\begin{equation}\label{local:continuity2}
		\dfrac{x^2}{r^2r_x} , \dfrac{v}{x}, v_x \in \mathcal C([0,T^{**}]\times (0,1)).
	\end{equation}
	
\end{lm}

\begin{spf} Here we present the local a prior estimate. The local well-posedness would follow from a similar iteration argument as in \cite{Jang2010} or difference arguments as in \cite[Appendix A]{LuoXinZeng2016}.
	Let $ P ( \cdot ) $ be the generic polynomial which might be different from line to line in the following. Also, for convenience, denote $$ M(t) = \sup_{0\leq s \leq t } \biggl\lbrace \norm{\dfrac{r}{x}}{\supnorm_x}, \norm{r_x}{\supnorm_x}, \norm{\dfrac{x}{r}}{\supnorm_x}, \norm{\dfrac{1}{r_x}}{\supnorm_x},  \norm{\dfrac{v}{x}}{\supnorm_x}, \norm{v_x}{\supnorm_x}  \biggr\rbrace. $$
	The initial date $ \dt v(x,0) = {\bar u}_1 (x) , \dt^2  v(x, 0) = {\bar u}_2(x) $ is given by \eqref{local:initial3}, \eqref{local:initial4}.
	We claim first and prove later, $ \exists \tilde T > 0 $, any $ 0 \leq t < \tilde T $, 
	\begin{equation}\label{local004}
		M(t) \leq C_{\tilde T} P(\sup_{0\leq s\leq t} \mathfrak E_0(s)). 
	\end{equation}
	Multiply \eqref{LgNS01} with $ \chi v_t $. After integration by parts in the resulting, similar arguments as before would yield,
	\begin{equation*}
	\begin{aligned}
		& \dfrac{d}{dt} \dfrac{1}{2} \int \chi \left( \dfrac{x}{r}\right)^2 \rho_0 v_t^2 \,dx + \int \chi \left( r_x \dfrac{v_t^2}{r^2} + \dfrac{v_{xt}^2}{r_x} \right) \,dx \\
		& ~~~~ \leq P(M(t)) \times \left\lbrace \int \chi \left( \dfrac{x}{r}\right)^2 \rho_0 v_t^2 \,dx  + \int \chi r_x \dfrac{v^2}{r^2} \,dx + \int \chi \dfrac{v_x^2}{r_x} \,dx \right. \\
		& ~~~~~~~~ \left.  + \int x^2 \rho_0 v_t^2 \,dx + \int r^2 \dfrac{v_x^2}{r_x} \,dx + \int r_x v^2 \,dx  \right\rbrace.
	\end{aligned}	
	\end{equation*}
	Therefore, together with \eqref{lm1} and \eqref{ee0031}, we shall have
	\begin{equation}\label{local005}
	\begin{aligned}
		& \dfrac{d}{dt} \mathfrak E_{0}(t) + \mathfrak D_{0} (t) \leq P(M(t)) \times \left\lbrace \mathfrak E_0 (t) +  \int r^2 \dfrac{v_x^2}{r_x} \,dx \right. \\
		& ~~~~~~ \left. + \int r_x v^2 \,dx  + \int \chi r_x \dfrac{v^2}{r^2} \,dx + \int \chi \dfrac{v_x^2}{r_x} \,dx  \right\rbrace.
	\end{aligned}	
	\end{equation}
	In the meantime, from \eqref{ee001}, after integration by parts, it should hold
	\begin{equation*}
%		\begin{aligned}
			\int r^2 \dfrac{v_x^2}{r_x} \,dx + \int r_x v^2 \,dx \leq P(M(t)) \times \mathfrak E_0(t). 
%		\end{aligned}
	\end{equation*}
	Similarly, multiplying \eqref{LgNS} with $ \chi v $ would eventually yield
	\begin{equation}\label{local009}
	\begin{aligned}
		& \int \chi r_x \dfrac{v^2}{r^2}\,dx + \int \chi \dfrac{v_x^2}{r_x}\,dx \leq P(M(t)) \times\left\lbrace \int r^2 \dfrac{v_x^2}{r_x}\,dx + \int r_x v^2 \,dx \right. \\
		& ~~~~ \left. + \int r^2 r_x \left( \dfrac{x^2 \rho_0}{r^2r_x}\right)^\gamma \,dx + \int \chi \left( \dfrac{x}{r}\right)^2 \rho_0 v_t^2 \,dx + 1 \right\rbrace.
	\end{aligned}
	\end{equation}
	Then \eqref{local005} can be written as 
	\begin{equation*}\label{local006}
	\dfrac{d}{dt} \mathfrak E_0(t) + \mathfrak D_0(t) \leq P(M(t)) \times\left\lbrace \mathfrak E_0(t) + 1 \right\rbrace. 	
	\end{equation*}
	Together with \eqref{local004}, then it can be shown there is a $ 0 <  T^* < \tilde T $ such that 
	\begin{equation}\label{local007}
		\sup_{0 \leq t\leq T^* }\mathfrak E_0(t) \leq 2 \mathfrak E_{0}(0),
	\end{equation}
	and
	\begin{equation}\label{local008}
	\int_0^{T^*} \mathfrak D_0(t) \,dt \leq 2 \mathfrak E_0(0).
	\end{equation}
	The regularity follows from similar arguments as in Section \ref{sec:regularity}. It should be noted when applying the arguments as in Lemma \ref{lm:D2regu}, the initial data for the ODE of $ \int \mathcal G_x^2 \,dx $ is given by the last integral in \eqref{local:initial}. 
	Now we briefly demonstrate \eqref{local004}. Denote
	$$ L(t) = \sup_{0\leq s\leq t} \left\lbrace \norm{\dfrac{r}{x}}{\supnorm_x}, \norm{r_x}{\supnorm_x} , \norm{\dfrac{x}{r}}{\supnorm_x},\norm{\dfrac{1}{r_x}}{\supnorm_x} \right\rbrace. $$
	Indeed, from \eqref{ueept001} and \eqref{ueept002}, 
	\begin{equation*}
	\begin{aligned}
		& \left| \dfrac{v_x}{r_x} \right|^2 , \left| \dfrac{v}{r} \right|^2 \leq H(L) \times \left\lbrace \int \chi \left(\dfrac{x}{r}\right)^2 \rho_0 v_t^2 \,dx + \int x^2 \rho_0 v_t^2 \,dx + \left. v^2 \right|_{x=1} + 1 \right\rbrace \\
		& \leq H(L) \times \left\lbrace \int \chi \left(\dfrac{x}{r}\right)^2 \rho_0 v_t^2 \,dx + \int x^2 \rho_0 v_t^2 \,dx + \int \left( r^2 \dfrac{v_x^2}{r_x} + r_x v^2 \right) \,dx + 1  \right\rbrace,
	\end{aligned}
	\end{equation*}
	where $ H(L) $ is a function of $ L $ and smooth when $ L \neq 0 $. 
	Meanwhile, from \eqref{ee001} again, after integration by parts, it can be derived,
	\begin{equation*}
		\begin{aligned}
			& \int \left( r^2 \dfrac{v_x^2}{r_x} + r_x v^2 \right) \,dx \leq C \int x^2 \rho_0 v^2 \,dx + C \int x^2 \rho_0 v_t^2 \,dx \\
			& ~~~~ + \delta \left\lbrace \int r^2 \dfrac{v_x^2}{r_x} \,dx +  \int r_x v^2 \,dx \right\rbrace + C_\delta H(L).
		\end{aligned}
	\end{equation*}
	Hence, the following inequality holds, 
	\begin{equation*}
		\left| v_x\right| , \left| \dfrac{v}{x} \right| \leq H(L) \times \left\lbrace \mathfrak E_0(t) + 1 \right\rbrace,
	\end{equation*}
	from which, together with the fact $ \dt r = v, \dt r_x = v_x $, \eqref{local004} follows. Thus we finish the local estimate for the strong solution. As for the smooth solution, multiply \eqref{LgNS02} with $ \chi v_{tt} $ and integrate the resulting. After integration by parts, it shall hold
	\begin{equation*}
		\begin{aligned}
			& \dfrac{d}{dt} \dfrac{1}{2} \int \chi \left( \dfrac{x}{r}\right)^2 \rho_0 v_{tt}^2 \,dx + \int \chi \left( \dfrac{v_{xtt}^2}{r_x} + r_x \dfrac{v_{tt}^2}{r^2} \right) \,dx \\
			& ~~~~ \leq P(M(t)) \times \left\lbrace \mathfrak E_{0} + \mathfrak E_{1} + \mathfrak D_0 + \int \chi \left(r_x \dfrac{v^2}{r^2} + \dfrac{v_x^2}{r_x}  \right) \,dx \right\rbrace.
		\end{aligned}
	\end{equation*}
	Together with \eqref{rg1201}, \eqref{local009}, 
	\begin{equation*}
		\dfrac{d}{dt} \mathfrak E_1 + \mathfrak D_1 \leq P(M(t)) \times \left\lbrace \mathfrak E_0 + \mathfrak E_1 + \mathfrak D_0 + 1 \right\rbrace.  
	\end{equation*}
	Therefore, by noticing \eqref{local004}, \eqref{local007} and \eqref{local008}, $ \exists T^{**} \leq T^* $ such that
	\begin{equation}
		\sup_{0 \leq t\leq T^{**}} \mathfrak E_1(t) \leq 3 \mathfrak E_0(0) + 2 \mathfrak E_1(0).
	\end{equation}
	The regularity would follow from the arguments as in Section \ref{sec:highregularity}. It should be noted when applying the arguments as in Lemma \ref{lm:D3regu}, the initial data for the ODE of $ \int \mathcal G_{xx}^2 \,dx $ is given by the last integral in \eqref{local:initial2}.  In particular, \eqref{local:continuity2} would be justified. 
\end{spf}

%We omit the proof of Lemma \ref{lm:local}. The regularity follows from the a prior estimate in Section \ref{sec:estimates}.

\subsection{Global Existence of Smooth solutions}

{\bf \large Proof of the Global Existence of Smooth Solutions, }

Given initial data $ (\rho_0, u_0) = (\rho_0, u_0)(x) $ satisfying \eqref{InitialAssum}, \eqref{InitialAssum2}, \eqref{UpBdInitial}, \eqref{InitialBoundary}, \eqref{InitialAssum3} and 
\begin{displaymath}
	\mathcal E_0, \mathcal E_1, \mathcal E_2 < \bar\epsilon,~ \mathcal E_3, \mathcal E_4 < \infty,
\end{displaymath}
from \eqref{lm906}, $ \exists \bar\alpha \in ( 1, \infty) $ % and $ \bar\beta > \max_{x\in(0,1)} \left( \left| u_{0,x} \right|, \left|\dfrac{u_0}{x}\right| \right) $
 such that $ \mathcal E_0, \mathcal E_1, \mathcal E_2 < \epsilon(\bar \alpha) $. Moreover, let $ \bar{\beta}_0 = \beta_0(\bar\alpha,\mathcal E_0, \mathcal E_1, \mathcal E_2, M), \bar{\beta}_1 = \beta_1(\bar\alpha,\mathcal E_0, \mathcal E_1, \mathcal E_2, M) $ be given in Lemma \ref{lm:extendinglemma}.  
From Lemma \ref{lm:local}, there is a smooth solution to \eqref{LgNS} with \eqref{LgNSBC} in $ (x, t) \in (0,1)\times (0,\delta_*) $,  $ \delta_* > 0 $. In particular, the regularity listed in \eqref{regularity2} holds. 
Moreover, the continuity in \eqref{local:continuity2} implies
\begin{equation}\label{ge:001}
	0 < \dfrac{x^2}{r^2r_x}\leq {\bar \alpha}^3, ~~ 0 \leq \left| \dfrac{v}{r} \right|, \left|\dfrac{v_x}{r_x}\right| \leq \bar{\beta}_0 \leq \bar{\beta}_1 .
\end{equation}
Denote the existence time as $$ T_* = \sup \left\lbrace T \geq \delta_* | \text{The smooth solution exists in $ t \in (0, T) $ and \eqref{ge:001} holds} \right\rbrace. $$ We claim $ T_* = \infty $. Otherwise, $ T_* < \infty $. Then from \eqref{lm903}, \eqref{lm904} in Section \ref{sec:estimates}, 
\begin{equation}
	\left\Arrowvert\dfrac{x^2}{r^2r_x}\right\Arrowvert_{L^{\infty}_x} (T_*) < \bar\alpha^3 , ~ \left\Arrowvert \dfrac{v_x}{r_x} \right\Arrowvert_{L_x^\infty}(T_*) < \bar{\beta}_0, ~\left\Arrowvert \dfrac{v}{r}\right\Arrowvert_{L^{\infty}_x}(T_*) < \bar{\beta}_0.
\end{equation}
Moreover, by choosing $ r(x, T_*), u(x, T_*)$ as initial data in \eqref
{local001} to the problem \eqref{LgNS} with the boundary condition \eqref{LgNSBC}, from the a priori estimates in Section \ref{sec:estimates}, it would satisfy \eqref{local:initial} and \eqref{local:initial2}. Therefore, Lemma \ref{lm:local} implies there is $ \delta > 0 $ such that there is a smooth solution in $ t \in [T_*, T_* + \delta) $. Additionally, \eqref{local:continuity2} implies that $ (r,v) $ satisfies \eqref{ge:001} for $ \delta $ small enough. This contradicts the definition of $ T_* $. This finishes the proof. \\
\linebreak
{\bf \large Proof of the Global Existence of Strong Solutions}

Consider \eqref{LgNS}, \eqref{LgNSBC}, \eqref{LgNSIN} with $ (\rho_0, u_0 ) $ satisfying \eqref{InitialAssum}, \eqref{InitialAssum2}, \eqref{UpBdInitial}, \eqref{InitialBoundary} and
\begin{equation}
	\mathcal E_0, \mathcal E_1, \mathcal E_2  < \bar\epsilon.
\end{equation}
Construct a smooth sequence $(\rho_{0,\iota} , u_{0,\iota}) \in C^{\infty} $ satisfying \eqref{InitialAssum}, \eqref{InitialAssum2}, \eqref{UpBdInitial}, \eqref{InitialBoundary}, \eqref{InitialAssum4} \eqref{InitialAssum3}, and as $ \iota \rightarrow 0 $, 
\begin{equation}
\begin{aligned}
	& \left\Arrowvert \rho_{0,\iota}^\gamma - \rho_0^\gamma\right\Arrowvert_{H_x^1(0,1)} \rightarrow 0 , &&  \norm{x \rho^{\gamma/2}_{0,l}-x\rho_0^{\gamma/2}}{L^2(0,1)} \rightarrow 0,   \\
	& \left\Arrowvert x \sqrt{\rho_{0,\iota}}u_{0,\iota} - x \sqrt{\rho_0}u_0 \right\Arrowvert_{L_x^2(0,1)} \rightarrow 0, 
	&& \left\Arrowvert x \sqrt{\rho_{0,\iota}}u_{1,\iota} - x \sqrt{\rho_0}u_1 \right\Arrowvert_{L_x^2(0,1)} \rightarrow 0 ,  \\
	&  \left\Arrowvert \sqrt{\rho_{0,\iota}}u_{1,\iota} - \sqrt{\rho_0}u_1 \right\Arrowvert_{L_x^2(0,1)} \rightarrow 0, && 
\end{aligned}
\end{equation}
where $ u_1 $ is defined as \eqref{intut1} and 
\begin{equation}
	u_{1,\iota} = \dfrac{1}{\rho_{0,\iota}}\left\lbrace (2\mu + \lambda) \left( \dfrac{(x^2 u_{0,\iota})_x}{x^2}\right)_x - \left( \rho_{0,\iota}^\gamma \right)_x \right\rbrace.
\end{equation}
Moreover $ (\rho_{0,\iota}, u_{0,\iota}) $ would admit the conditions listed in Lemma \ref{lm:local} for the smooth solution to exist. 
Then for any fixed $ T>0 $, there is a smooth solution $ (r_\iota, v_\iota ) $ for $ t \in (0,T) $ satisfying \eqref{ge:001} with some $ \bar\alpha,\bar{\beta}_0>M $. 
 Notice, from our choice of the approximation sequence, the estimates in Sections \ref{sec:energyestimates}, \ref{sec:uniformestimate}, \ref{sec:regularity} are independent of $ \iota $.
 In particular, the following regularity hold regardless of $ \iota $,
\begin{displaymath}
\begin{aligned}
	& \left\Arrowvert x \sqrt{\rho_{0,\iota}} v_\iota \right\Arrowvert_{L_t^\infty((0,T),L_x^2(0,1))}, \left\Arrowvert x \sqrt{\rho_{0,\iota}} v_{\iota,t} \right\Arrowvert_{L_t^\infty((0,T),L_x^2(0,1))} \leq C_T,\\
	& \left\Arrowvert \sqrt{\rho_{0,\iota}} v_{\iota,t} \right\Arrowvert_{L_t^\infty((0,T),L_x^2(0,1))}, 
	\left\Arrowvert v_\iota \right\Arrowvert_{L_t^\infty((0,T),H_x^2(0,1))},\left\Arrowvert \dfrac{v_\iota}{x} \right\Arrowvert_{L_t^\infty((0,T),H_x^1(0,1))} \leq C_T,\\
%	& \left\Arrowvert \dfrac{v_\iota}{x} \right\Arrowvert_{L_t^\infty((0,T),H_x^1(0,1))} \leq C_T,\\
	& \left\Arrowvert xv_{\iota,x}\right\Arrowvert_{L_t^2((0,T),L_x^2(0,1))}, 
	\left\Arrowvert x v_{\iota,xt}\right\Arrowvert_{L_t^2((0,T),L_x^2(0,1))}, 
	\left\Arrowvert v_\iota \right\Arrowvert_{L_t^2((0,T),L_x^2(0,1))} \leq C_T, \\
	& \left\Arrowvert v_{\iota,t}\right\Arrowvert_{L_t^2((0,T),L_x^2(0,1))}, 
	\left\Arrowvert v_{\iota,xt}\right\Arrowvert_{L_t^2((0,T),L_x^2(0,1))}, 
	\left\Arrowvert \dfrac{v_{\iota,t}}{x} \right\Arrowvert_{L_t^2((0,T),L_x^2(0,1))} \leq C_T,\\
	& 0 < \dfrac{x^2}{r_\iota^2 r_{\iota,x}} \leq \bar{\alpha}^3,~ \left\Arrowvert \dfrac{v_\iota}{r_\iota}\right\Arrowvert_{L_x^\infty(0,1)},  \left\Arrowvert \dfrac{v_{\iota,x}}{r_{\iota,x}}\right\Arrowvert_{L_x^\infty(0,1)} \leq \bar{\beta}_0.
\end{aligned}
\end{displaymath}
Then, by taking $ \iota \rightarrow 0 $, standard compactness arguments would  yield
$$ (r_\iota, v_{\iota}) \rightarrow (r,v) ~~ \text{strongly in} ~ C((0,T), H_x^1(0,1)), $$
and $(r,v) $ is a strong solution to \eqref{LgNS} in $ t \in (0,T) $ and the regularity listed in \eqref{regularity1} holds. This finishes the proof.

\paragraph{\bf Acknowledgements}
This work is part of the doctoral dissertation of the author under the supervision of Professor Zhouping Xin at the Institute of Mathematical Sciences of the Chinese University of Hong Kong, Hong Kong. The author would like to express great gratitude to Prof. Xin for his kindly support and professional guidance.


\begin{thebibliography}{}

\bibitem{Cho2006} Y. Cho, B.J. Jin, Blow-up of viscous heat-conducting compressible flows, J. Math. Anal. Appl. 320(2006) 819-826.

\bibitem{Cho2006a} Y. Cho, H. Kim, Existence results for viscous polytropic fluids with vacuum, J. Differential Equations 228(2006) 377-411.

\bibitem{Cho2006c} Y. Cho, H. Kim, On classical solutions of the compressible Navier-Stokes equations with nonnegative initial densities, Manuscripta Math. 120(2006) 91-129.

\bibitem{Coutand2010} D. Coutand, H. Lindblad, S. Shkoller, A priori estimates for the free-boundary 3D compressible Euler equations in physical vacuum, Commun. Math. Phys. 296(2010) 559-587.

\bibitem{Coutand2011a} D. Coutand, S. Shkoller, Well-posedness in smooth function spaces for moving-boundary 1-D compressible Euler equations in physical vacuum, Commun. Pure Appl. Math. LXIV(2011) 0328-0366.

\bibitem{Coutand2012} D. Coutand, S. Shkoller, Well-posedness in smooth function spaces for the moving-boundary three-dimensional compressible Euler equations in physical vacuum, Arch. Rational Mech. Anal. 206(2012) 515--616.

\bibitem{Gu2012} X. Gu, Z. Lei, Well-posedness of 1-D compressible Euler-Poisson equations with physical vacuum, J. Differential Equations 252(2012) 2160-2188.

\bibitem{Gu2015} X. Gu, Z. Lei, Local well-posedness of the three-dimensional compressible Euler-Poisson equations with physical vacuum, J. Math. Pures Appl. 105(2016) 662-723.

\bibitem{Guo2012a} Z. Guo, H.-L. Li, Z. Xin, Lagrange structure and dynamics for solutions to the spherically symmetric compressible Navier-Stokes equation, Commun. Math. Phys. 309(2012) 371-412

\bibitem{Guo2012b} Z. Guo, Z. Xin, Analytical solutions to the compressible Navier-Stokes equations with density-dependent viscosity coefficients and free boundaries, J. Differential Equations 253(2012) 1-19

\bibitem{HuangLiXin2012} X. Huang, J. Li, Z. Xin, Global well-posedness of classical solutions with large oscillations and vacuum to the three-dimensional isentropic compressible Navier-Stokes equations, Commun. Pure Appl. Math. 65(2012) 0549-0585. 

\bibitem{Itaya1971} N. Itaya, On the Cauchy problem for the system of fundamental equations describing the movement of compressible viscous fluid, Kodai Math. Sem. Rep. 23(1971) 60-120.

\bibitem{Jang2010} J. Jang, Local well-posedness of dynamics of viscous gaseous stars, Arch. Rational Mech. Anal. 195(2010) 797-863.

\bibitem{Jang2009b} J. Jang, N. Masmoudi, Well-posedness for compressible Euler equations with physical vacuum singularity, Commun. Pure Appl. Math. LXII(2009) 1327-1385.

\bibitem{Jang2011} J. Jang, N. Masmoudi, Vacuum in gas and fluid dynamics, in: A. Bressan, G.-Q. G. Chen, M. Lewicka, D. Wang(Eds.), Nonlinear conservation laws and applications, Springer Science+Bussiness Media, LLC 2011, pp. 315-329.

\bibitem{Jang2015} J. Jang, N. Masmoudi, Well-posedness of compressible Euler equations in a physical vacuum, Commun. Pure Appl. Math. LXVIII(2015) 0061-0111.

\bibitem{Li2016} H.-L. Li, X. Zhang, Global strong solutions to radial symmetric compressible Navier-Stokes equations with free boundary, J. Differential Equations 261(2016) 6341-6367.

\bibitem{Liu1996} T.-P. Liu, Compressible flow with damping and vacuum, Japan J. Indust. Appl. Math. 13(1996) 25-32.

\bibitem{Liu1998a} T.-P. Liu, Z. Xin, T. Yang, Vacuum states for compressible flow, Discret. Contin. Dyn. S. 4(1998) 1-32.

\bibitem{Luo2000a} T. Luo, Z. Xin, T. Yang, Interface behavior of compressible Navier-Stokes equations with vacuum, SIAM J. Math. Anal. 31(2000) 1175-1191.

\bibitem{LuoXinZeng2014} T. Luo, Z. Xin, H. Zeng, Well-posedness for the motion of physical vacuum of the three-dimensional compressible Euler equations with or without self-gravitation, Arch. Rational Mech. Anal. 213(2014) 763-831.

\bibitem{LuoXinZeng2015} T. Luo, Z. Xin, H. Zeng, Nonlinear asymptotic stability of the Lane-Emden solutions for the viscous gaseous star problem with degenerate density dependent viscosities, Commun. Math. Phys. 347(2016) 657-702. 

\bibitem{LuoXinZeng2016} T. Luo, Z. Xin, H. Zeng, On nonlinear asymptotic stability of the Lane-Emden solutions for the viscous gaseous star problem, Advances in Mathematics 291(2016) 90-182. 

\bibitem{Matsumura1980} A. Matsumura, T. Nishida, The initial value problem for the equations of motion of viscous and heat-conductive gases, J. Math. Kyoto Univ. (JMKYAZ) 20-1(1980) 67-104.

\bibitem{Matsumura1983} A. Matsumura, T. Nishida, Initial boundary value problems for the equations of motion of compressible viscous and heat-conductive fluids, Commun. Math. Phys. 89(1983) 445-464.

\bibitem{Okada1989} M. Okada, Free boundary value problems for the equation of one-dimensional motion of viscous gas, Japan J. Appl. Math. 6(1989) 161-177.

\bibitem{Okada2004} M. Okada, Free boundary problem for one-dimensional motions of compressible gas and vacuum, Japan J. Indust. Appl. Math. 21(2004) 109-128.

\bibitem{Ou2015} Y. Ou, H. Zeng, Global strong solutions to the vacuum free boundary problem for compressible Navier-Stokes equations with degenerate viscosity and gravity force, J. Differential Equations 259(2015) 6803-6829. 

\bibitem{Serrin1959} J. Serrin, On the uniqueness of compressible fluid motions, Arch. Rational Mech. Anal. 3(1959) 271-288.

\bibitem{Solonnikov1992} V.A. Solonnikov, A. Tani, A problem with free boundary for the Navier-Stokes equations for a compressible fluid in the presence of surface tension, Journal of Soviet Mathematics 62(1992) 2814-2818.

\bibitem{Tani1977} A. Tani, On the first initial-boundary value problem of compressible viscous fluid motion, Publ. RIMS, Kyoto Univ. 13(1977) 193-253. 

\bibitem{Tani1986} A. Tani, The initial value problem for the equations of the motion of compressible viscous fluid with some slip boundary condition, in: T. Nishida, M. Mimura, H. Fujii(Eds.), Studies in Mathematics and Its Applications Volume 18 - Patterns and Waves-Qualitative Analysis of Nonlinear Differential Equations, Kinokuniya Company Ltd. Published by Elsevier B.V., 1986, pp.675-684

\bibitem{zpxin1998} Z. Xin, Blowup of smooth solutions to the compressible Navier-Stokes equation with compact density, Commun. Pure Appl. Math. LI(1998) 0229-0240. 

\bibitem{XinYan2013} Z. Xin, W. Yan, On blowup of classical solutions to the compressible Navier-Stokes equations, Commun. Math. Phys. 321(2013) 529-541.

\bibitem{Yeung2009} Ling Hei Yeung, Yuen Manwai, Analytical solutions to the Navier-Stokes equations with density-dependent viscosity and with pressure, J. Math. Phys. 50, 083101 (2009), doi: 10.1063/1.3197860

\bibitem{Zadrzynska2001} E. Zadrzy\'nska, Evolution free boundary problem for equations of viscous compressible heat-conducting capillary fluids, Math. Meth. Appl. Sci. 24(2001) 713-743.

\bibitem{Zajaczkowski1999} E. Zadrzy\'nska, W.M. Zaj\c{a}czkowski, On nonstationary motion of a fixed mass of a viscous compressible baratropic fluid bounded by a free boundary, Colloquium Mathematicum 79(1999) 283-310.

\bibitem{Zadrzynska1994} E. Zadrzy\'nska, W.M. Zaj\c{a}czkowski, On local motion of a general compressible viscous heat conducting fluid bounded by a free surface, Annales Polonici Mathematici LIX-2(1994) 133-170.

\bibitem{Zajaczkowski1993} W.M. Zaj\c{a}czkowski, On nonstationary motion of a compressible barotropic viscous fluid bounded by a free surface, Dissertationes Mathematicae Rozprawy Matematyczne tom/nr w serii: 324 wydano, 1993.

\bibitem{Zajaczkowski1994} W.M. Zaj\c{a}czkowski, On nonstationary motion of a compressible barotropic viscous capillary fluid bounded by a free surface, SIAM J. Math. Anal. 25(1994) 1-84. 

\bibitem{Zajaczkowski1995} W.M. Zaj\c{a}czkowski, Existence of local solutions for free boundary problems for viscous compressible barotropic fluids, Annales Polonici Mathematici LX-3(1995) 255-287.

\bibitem{ZengHH2015} H. Zeng, Global-in-time smoothness of solutions to the vacuum free boundary problem for compressible isentropic Navier-Stokes equations, Nonlinearity 28(2015) 331-345. 

\bibitem{Zhang2007a} T. Zhang, D. Fang, Existence and uniqueness results for viscous, heat-conducting 3-D fluid with vacuum, arXiv:math/0702170 [math.AP], 2007.

\bibitem{Zhang2009c} T. Zhang, D. Fang, A note on spherically symmetric isentropic compressible flows with density-dependent viscosity coefficients, Nonlinear Analysis: Real World Applications 10(2009) 2272-2285.

\bibitem{Zhang2009} T. Zhang, D. Fang, Global behavior of spherically symmetric Navier-Stokes-Poisson system with degenerate viscosity coefficients, Arch. Rational Mech. Anal. 191(2009) 195-243.

\end{thebibliography}
\end{document}